\documentclass[11pt,final]{article}
\usepackage{a4wide}
\usepackage{amsmath,amssymb,amsfonts}
\usepackage[utf8]{inputenc} 
\usepackage[T1]{fontenc}
\usepackage{enumerate}
\usepackage{pdfpages}
\usepackage{mathrsfs}
\usepackage{graphicx}
\usepackage{tabularx}
\usepackage{xcolor}
\usepackage{url}
\newcommand{\MyBibitem}[2]{\bibitem{#1}}
\newcommand{\MyReference}[3]{\textsc{#1}. \emph{#2}. #3.}
\newcommand{\MyReferenceOhnePunkt}[3]{\textsc{#1}. \emph{#2} #3.}
\newcommand{\MyParagraph}[1]{$\,$\\[-3mm] \noindent \textbf{#1.}\;}
\newcommand{\MyUrl}[1]{#1}
\newcommand{\MyCite}[1]{\cite{#1}}

\newcommand{\ZZ}{\mathbb{Z}}
\newcommand{\RR}{\mathbb{R}}
\newcommand{\CC}{\mathbb{C}}
\newcommand{\dd}{\text{d}}
\newcommand{\ee}{\text{e}}
\newcommand{\ii}{\text{i}}
\newcommand{\abs}[1]{|#1|}

\newcommand{\scp}[3]{(#2\,|\,#3)_{#1}}
\newcommand{\scpbig}[3]{\big(#2\,\big|\,#3\big)_{#1}}
\newcommand{\norm}[2]{\|#2\|_{#1}}
\newcommand{\normbig}[2]{\big\|#2\big\|_{#1}}
\newcommand{\nE}{\mathcal{E}}
\newcommand{\nS}{\mathcal{S}}
\newcommand{\nF}{\mathcal{F}}

\newcommand{\Nablaalpha}{\nabla_{\!\alpha}}
\newcommand{\Nablaalphaone}{\nabla_{\!\alpha_1}}
\newcommand{\Nablaalphatwo}{\nabla_{\!\alpha_2}}
\newcommand{\Nablaalphaj}{\nabla_{\!\alpha_j}}
\newcommand{\Deltaalpha}{\Delta_{\alpha}}
\newcommand{\Deltaalphaj}{\Delta_{\alpha_j}}
\newcommand{\Deltaalphaone}{\Delta_{\alpha_1}}
\newcommand{\Deltaalphatwo}{\Delta_{\alpha_2}}

\newcommand{\Vbetagammadelta}{V_{\beta,\gamma,\delta}}
\newcommand{\Vbetagammadeltaj}{V_{\beta_j,\gamma_j,\delta_j}}

\newcommand{\mualpham}{\mu_{\alpha,m}}
\newcommand{\lambdaalpham}{\lambda_{\alpha,m}}
\newcommand{\wtPsi}{\tilde{\Psi}}
\newcommand{\wtF}{\tilde{F}}
\newcommand{\wtG}{\tilde{G}}
\newcommand{\wtH}{\tilde{H}}
\begin{document}
\title{$\,$\\[-3cm]\textbf{Modified splitting methods for Gross--Pitaevskii systems modelling Bose--Einstein condensates: \\ Time evolution and ground state computation}%
\footnote{
This work has been supported by the Austrian Science Fund FWF through a stand-alone project (grant-doi 10.55776/PAT1281625).
Fruitful discussions during a research stay at the Wolfgang Pauli Institute Vienna have inspired this work; Mechthild Thalhammer is grateful to the director Norbert Mauser and the staff members for their support and hospitality. 
}}
\author{Mechthild Thalhammer, Gregor Thalhammer--Thurner%
\footnote{Addresses:
Mechthild Thalhammer, University of Innsbruck, Department of Mathematics, 6020~Innsbruck, Austria.
Gregor Thalhammer--Thurner, Medical University of Innsbruck, Institute of Biomedical Physics, 6020~Innsbruck, Austria.
Websites: \MyUrl{techmath.uibk.ac.at/mecht}, \MyUrl{www.i-med.ac.at/dpmp/bmp}.
Email addresses: \MyUrl{mechthild.thalhammer@uibk.ac.at}, \MyUrl{gregor.thalhammer@i-med.ac.at}.
}}
\maketitle
\noindent
{\small
The year 2025 marks the 100 and 30 years anniversaries of the discovery of Bose--Einstein condensation and its successful experimental realisation.
Inspired by these important research achievements, a conceptually simple approach is proposed to facilitate reliable and efficient numerical simulations.
The structure of the underlying systems of coupled Gross--Pitaevskii equations suggests the use of optimised high-order operator splitting methods for dynamical evolution and ground state computation.
A second-order barrier, however, prevents the applicability of standard operator splitting methods for both, time evolution as well as imaginary time propagation.
An innovative alternative approach accomplishes the design of novel modified operator splitting methods that remain stable under moderate smallness assumptions on the time increments. 
The core idea is to incorporate commutators of the defining differential and nonlinear multiplication operators, since this permits to fulfill the basic stability requirement of positive method coefficients.
Further improvements with respect to convergence at the targeted precision arise from automatic adjustments of the time stepsizes by an inexpensive local error control.
The presented numerical experiments confirm the favourable performance of a specific fourth-order modified operator splitting method.
Amongst others, it is demonstrated that the excellent mass and energy conservation in long-term evolutions, intrinsic attributes of geometric numerical integrators for Hamiltonian systems, is maintained for a sensible variation of the time stepsizes.
Moreover, the benefits of adaptive higher-order approximations in ground state computations are illustrated.
}
\section{Introduction}
\MyParagraph{Main objective}
In light of the 100 and 30 years anniversaries of the theoretical discovery of Bose--Einstein condensation~\MyCite{Einstein1925} and its successful experimental realisation at ultracold temperatures~\MyCite{AndersonEtAl1905}, a main objective of the present work is the introduction and exemplification of a novel and conceptually simple approach, in order to facilitate the reliable and efficient numerical simulation of multi-species Bose--Einstein condensates.

\MyParagraph{Bose--Einstein condensates}
The creation of Bose--Einstein condensate collective systems occurs when a gas of bosons is cooled to nanokelvin temperatures by laser technologies, since the controlled reduction of atomic motion causes a large number of microscopic particles to form a single macroscopic quantum state at lowest energy. 
Worldwide, various experimental physics groups devote their current research to the study and manipulation of this particular state of matter.
As an instance, we highlight the \emph{Strongly Correlated Quantum Matter Group}%
\footnote{See \MyUrl{quantummatter.at}.}
at the University of Innsbruck, led by \textsc{Hanns--Christoph N{\"a}gerl}, and the recent contribution~\MyCite{HCN2025}. 

\MyParagraph{Gross--Pitaevskii systems}
Fundamental mathematical models that describe the dynamical evolution of the macroscopic wave functions associated with multi-species Bose--Einstein condensates are provided by systems of coupled time-dependent Gross--Pitaevskii equations~\MyCite{Gross1961,Pitaevskii1961}.
We focus on typical instances that are characterised by partial differential equations of Schr{\"o}dinger-type comprising the Laplacian with respect to the spatial coordinates, real-valued space-dependent functions, and non-analytic nonlinearities.
The defining operators reflect kinetic energies, external trapping and lattice potentials, and the interactions between the particles of different species. 
The physically most relevant setting corresponds to three space dimensions and cubic nonlinearities.
In addition, we study the special cases of one and two dimensions with the advantages of accelerating numerical tests and simplifying the visualisations of solution profiles. 

\MyParagraph{Operator splitting approach}
The structures of the defining linear differential and nonlinear multiplication operators suggest the use of standard operator splitting methods combined with fast Fourier techniques for spatial discretisation.
The main idea is to divide the original problem into subproblems and to apply adjusted solvers to each of the subproblems.
Employing a usual time-stepping approach, the respective approximations are composed several times, whereby the time increments are scaled by certain method coefficients.
The particular choice of the method coefficients, linked to the number of compositions per time step, determines the overall quantitative and qualitative properties of a specific splitting scheme.
Detailed information on splitting methods is found in~\MyCite{BlanesCasasMurua2024,McLachlanQuispel2002}. 
A small selection of contributions in the context of linear and nonlinear Schr{\"o}dinger equations is \MyCite{Bao2004,BlanesCasasMurua2015,CaliariZuccher2021,Thalhammer2008,Thalhammer2012,ThalhammerAbhau2012}, see also the references given therein.

\MyParagraph{Geometric numerical integrators}
Optimised high-order splitting methods exhibit a favourable behaviour regarding stability, accuracy, and efficiency.
By reason of their distinguishing attributes of geometric numerical integrators for Hamiltonian systems \MyCite{BlanesCasas2016,HairerLubichWanner2006, IserlesQuispel2018,SanzSernaCalvo2018}, they also show an excellent preservation of conserved quantities such as mass and energy over long times.

\MyParagraph{Limitations for non-reversible systems}
However, the unified applicability of standard operator splitting methods for both, the time evolution and the imaginary time propagation of Gross--Pitaevskii systems, is limited by a second-order barrier for systems that are non-reversible in time.
In other words, due to the fact that higher-order splitting methods necessarily involve negative coefficients~\MyCite{BlanesCasas2005} and thus suffer from severe instabilities, they are excluded from ground and excited state computations.
Hence, the symmetric second-order Strang splitting method is left as canonical choice.

\MyParagraph{Modified operator splitting approach}
In this contribution, we exploit an alternative approach, which we refer to as modified operator splitting, and illustrate its benefits over standard splitting methods.
Within the context of \emph{linear} differential equations, the basic concept is known as modified potential or force-gradient operator splitting, see~\MyCite{AuerKrotscheckChin2001,BlanesCasasGonzalezThalhammer2023,Chin1997,Chin2005,Kieri2015,OmelyanMryglodFolk2002}.
In our recent work~\MyCite{BlanesCasasGonzalezThalhammer2024}, we could accomplish the generalisation to a class of nonlinear partial differential equations by means of the formal calculus of Lie-derivatives~\MyCite{Varadarajan1984}.
A key observation is that the incorporation of double commutators of the defining differential and multiplication operators in the nonlinear subproblems accomplishes for the design of fourth-order schemes with positive coefficients.
In what follows, we establish modified operator splitting methods for the first time for Gross--Pitaevskii systems.
Provided that the time stepsizes are moderately small, we demonstrate that modified operator splitting methods with positive coefficients remain stable for time-dependent Gross--Pitaevskii systems as well as for their parabolic-type counterparts arising in the imaginary time propagation for ground and excited state computations, see~\MyCite{BaoDu2004,BlanesCasasGonzalezThalhammer2023,DanailaProtas2017,LehtovaaraToivanenEloranta2007}.

\MyParagraph{Adaptivity in time}
With regard to the evident advantages of adaptivity in time, we propose an inexpensive local error control. 
According to the requirements of the considered problem and the ranges of the decisive physical parameters, this strategy automatically adapts the time increments to ensure stability and also enhances the targeted precision.
We measure the effort of the resulting adaptive modified splitting methods by the numbers of fast Fourier transforms and their inverses, since these procedures constitute the computationally most costly components.

\MyParagraph{Efficient realisation}
Besides, we clarify the validity of an invariance principle that has a significant impact on the efficient realisation of modified splitting methods for time-dependent Gross--Pitaevskii systems.
In essence, it justifies the replacement of nonlinear terms depending on the time-dependent solutions by expressions involving only the initial solution values.
As a consequence, for the nonlinear subproblems that result from omitting the Laplacians, this invariance property provides the knowledge of the exact solutions on each of the temporal subintervals.  

\MyParagraph{Numerical experiments}
Our numerical experiments for model problems in one, two, and three space dimensions illustrate the benefits of a specific adaptive fourth-order modified operator splitting method in comparison with standard splitting methods based on uniform time grids.
Amongst others, we confirm the favourable energy conservation in long-term evolutions and demonstrate the computation of ground states with high accuracy.

\MyParagraph{Notation}
In the following, we denote by~$L^2(\Omega, \CC)$ the Lebesgue space of square integrable complex-valued functions, which forms a Hilbert space with inner product and associated norm defined by 
\begin{subequations}
\label{eq:ScpNorm}  
\begin{equation}
\scp{L^2}{\varphi}{\psi} = \int_{\Omega} \varphi(x) \, \overline{\psi(x)} \; \dd x\,, \quad 
\norm{L^2}{\psi} = \sqrt{\scp{L^2}{\psi}{\psi}}\,, \quad \varphi, \psi \in L^2(\Omega, \CC)\,.
\end{equation}  
Accordingly, for vector-valued functions, we set 
\begin{equation}
\begin{gathered}
\scp{(L^2)^d}{\varphi}{\psi} = \sum_{j=1}^{d} \scp{L^2}{\varphi_j}{\psi_j}\,, \quad 
\norm{(L^2)^d}{\psi} = \sqrt{\sum_{j=1}^{d} \norm{L^2}{\psi_j}^2}\,, \\
\varphi = (\varphi_1, \dots, \varphi_d) \in L^2(\Omega, \CC^d)\,, \quad \psi = (\psi_1, \dots, \psi_d) \in L^2(\Omega, \CC^d)\,.
\end{gathered}
\end{equation}  
\end{subequations}
The arising derivatives of functionals are determined as G{\^a}teaux derivatives generalising directional derivatives 
\begin{equation}
\label{eq:Gateaux}  
F'(v) \, w = \lim_{\varepsilon \to 0} \tfrac{1}{\varepsilon} \, \big(F(v + \varepsilon \, w) - F(v)\big)\,.
\end{equation}
Throughout, we tacitly assume that for the situations under consideration sufficiently regular solutions are ensured such that the application of high-order space and time discretisation methods is indeed expedient.

\MyParagraph{Outline}
The present work is organised as follows.   
In Sections~\ref{sec:GPE} and~\ref{sec:GPS}, we state Gross--Pitaevskii equations as well as extended systems modelling multi-species Bose--Einstein condensates.
A compact symbolic notation permits to emphasise the common structures of the initial value problems for dynamical evolution and ground state computation.
In Section~\ref{sec:Splitting}, we first give a brief summary of standard operator splitting methods and then detail our novel approach of modified operator splitting for the time evolution and imaginary time propagation of Gross--Pitaevskii systems. 
Comparative numerical experiments regarding the stability and accuracy behaviour of modified operator splitting methods in short-term as well as long-term evolutions for two-component systems are presented in Section~\ref{sec:Experiments}.
Concluding remarks are given in the final section~\ref{sec:Conclusions}.
\section{Gross--Pitaevskii equations}
\label{sec:GPE}
\MyParagraph{Model problem}
As a basic model for the nonlinear dynamics of a single-species Bose--Einstein condensate, we consider the initial value problem for a time-dependent Gross--Pitaevskii equation.
For our purposes, it suffices to study the dimensionless formulation
\begin{subequations}
\label{eq:GPE}  
\begin{equation}
\label{eq:GPE1}  
\begin{cases}
\ii \, \partial_t \Psi(x, t) = \Deltaalpha \Psi(x, t) + \Vbetagammadelta(x) \, \Psi(x, t) + \vartheta \, \abs{\Psi(x, t)}^2 \, \Psi(x, t)\,, \\
\Psi(x, t_0) = \Psi_0(x)\,, \quad (x, t) \in \Omega \times [t_0, T]\,.
\end{cases}
\end{equation}
Here, we denote by $\ii \in \CC$ the imaginary unit, by
\begin{equation}
\Psi: \Omega \times [t_0, T] \longrightarrow \CC: (x, t) \longmapsto \Psi(x, t)
\end{equation}  
the complex-valued macroscopic wave function, defined on a certain space domain $\Omega \subseteq \RR^d$ and time interval $[t_0, T] \subset \RR$, by
\begin{equation}
\begin{gathered}  
\Deltaalpha = \sum_{i=1}^{d} \alpha_i \, \partial_{x_i}^2\,, \\
\alpha_i \in \RR \setminus \{0\}\,, \quad i \in \{1, \dots, d\}\,, \quad x = (x_1, \dots, x_d) \in \Omega\,, 
\end{gathered}  
\end{equation}  
the weighted Laplacian with respect to the spatial coordinates, by
\begin{equation}
\begin{gathered}
\Vbetagammadelta: \Omega \longrightarrow \RR: x \longmapsto \sum_{i=1}^{d} \big(\beta_i \, x_i^2 + \gamma_i \, \sin^2(\delta_i \, x_i)\big)\,, \\
\beta_i > 0\,, \quad \gamma_i, \delta_i \in \RR\,, \quad i \in \{1, \dots, d\}\,,
\end{gathered}
\end{equation}  
\end{subequations}
an external harmonic potential, possibly comprising an additional lattice, and by $\vartheta \in \RR$ the coupling constant.

\MyParagraph{Compact reformulation}
Henceforth, we denote by $\psi: \Omega \to \CC$ a sufficiently regular space-dependent function and employ convenient abbreviations
\begin{subequations}
\label{eq:GPECompact}
\begin{equation}
\begin{gathered}
H(\psi) = H_1(\psi) + H_2(\psi)\,, \\
H_1(\psi) = \Deltaalpha \psi + \Vbetagammadelta \, \psi\,, \quad H_2(\psi) = \vartheta \, \abs{\psi}^2 \, \psi\,, 
\end{gathered}
\end{equation}
for the operators defining the right-hand side of the partial differential equation in~\eqref{eq:GPE1}. 
Consequently, we in particular obtain the compact reformulation
\begin{equation}
\begin{cases}
\partial_t \Psi(x, t) = - \, \ii \, H\big(\Psi(x, t)\big)\,, \\
\Psi(x, t_0) = \Psi_0(x)\,, \quad (x, t) \in \Omega \times [t_0, T]\,.
\end{cases}
\end{equation}
\end{subequations}

\MyParagraph{Conserved quantities}
A fundamental attribute of time-dependent Gross--Pitaevskii equations is the preservation of physical quantities such as total mass (or particle number) and energy 
\begin{subequations}
\label{eq:MassEnergy}
\begin{equation}
\begin{gathered}
M(\psi) = \norm{L^2}{\psi}^2\,, \quad E(\psi) = E_1(\psi) + \tfrac{1}{2} \, E_2(\psi)\,, \\
E_1(\psi) = \scpbig{L^2}{H_1(\psi)}{\psi}\,, \quad E_2(\psi) = \scpbig{L^2}{H_2(\psi)}{\psi}\,, 
\end{gathered}
\end{equation}
see~\eqref{eq:ScpNorm} and~\eqref{eq:GPE}--\eqref{eq:GPECompact}.
Integration-by-parts under the requirement that the boundary terms vanish yields 
\begin{equation}
E_1(\psi) = - \sum_{i=1}^{d} \alpha_i \, \big(\partial_{x_i}^2 \psi\big)^2 + \scpbig{L^2}{\Vbetagammadelta}{\abs{\psi}^2}\,, \quad E_2(\psi) = \vartheta \, \scpbig{L^2}{\abs{\psi}^2}{\abs{\psi}^2}\,,
\end{equation}
which reveals that the energy functional is real-valued and positive for non-trivial solutions, provided that the weights in the Laplacian are negative and that the potential as well as the coupling constant are non-negative.
Straightforward calculations show that the G{\^a}teaux derivatives are given by 
\begin{equation}
\label{eq:MassEnergyGateaux}
\begin{gathered}
M'(\psi) \, \phi = 2 \, \Re \scp{L^2}{\psi}{\phi} = 2 \, \Im \big(\ii \, \scp{L^2}{\psi}{\phi}\big)\,, \\ 
E'(\psi) \, \phi = 2 \, \Re \scp{L^2}{H(\psi)}{\phi} = 2 \, \Im\big(\ii \, \scp{L^2}{H(\psi)}{\phi}\big)\,,
\end{gathered}
\end{equation}
see also~\eqref{eq:Gateaux}.
As a consequence, for regular solutions to the Gross--Pitaevskii equation~\eqref{eq:GPE}--\eqref{eq:GPECompact}, this implies 
\begin{equation*}
\begin{gathered}
\tfrac{\dd}{\dd t} M\big(\Psi(\cdot, t)\big) = M'\big(\Psi(\cdot, t)\big) \, \partial_t \Psi(\cdot, t) = - \, 2 \, \Im\Big(E_1\big(\Psi(\cdot, t)\big) + E_2\big(\Psi(\cdot, t)\big)\Big) = 0\,, \\ 
\tfrac{\dd}{\dd t} E\big(\Psi(\cdot, t)\big) = E'\big(\Psi(\cdot, t)\big) \, \partial_t \Psi(\cdot, t) = - \, 2 \, \Im\Big(\normbig{L^2}{H\big(\Psi(\cdot, t)\big)}^2\Big) = 0\,, \\
t \in [t_0, T]\,,
\end{gathered}
\end{equation*}
and thus confirms the characterisations of mass and energy conservation as 
\begin{equation}
M\big(\Psi(\cdot, t)\big) = M\big(\Psi(\cdot, t_0)\big)\,, \quad E\big(\Psi(\cdot, t)\big) = E\big(\Psi(\cdot, t_0)\big)\,, \quad t \in [t_0, T]\,.
\end{equation}
\end{subequations}

\MyParagraph{Regularity and normalisation conditions}
For the time-dependent Gross--Pitaevskii equation~\eqref{eq:GPE}--\eqref{eq:GPECompact}, due to mass conservation, it is justified to impose the normalisation condition 
\begin{equation*}
\norm{L^2}{\Psi(\cdot, t)} = \sqrt{N_0}\,, \quad t \in [t_0, T]\,, 
\end{equation*}  
with a strictly positive real number $N_0 > 0$.
Besides, we tacitely suppose that the prescribed initial states and hence the associated solutions are sufficiently regular such that in particular the energy functional is well-defined.
Stronger regularity requirements are needed in the context of higher-order operator splitting methods to retain the full orders of convergence. 

\MyParagraph{Ground and excited states}
Ground and excited states are solutions to Schr{\"o}dinger equations with separate dependencies on space and time variables.
More precisely, inserting the ansatz 
\begin{equation*}
\begin{gathered}
\Psi(x, t) = e(t) \, \Phi(x)\,, \quad (x, t) \in \Omega \times [t_0, T]\,, \\ 
e: [t_0, T] \longrightarrow \CC\,, \quad \abs{e(t)} = 1\,, \quad t \in [t_0, T]\,, \\
\Phi: \Omega \longrightarrow \RR\,, \quad \norm{L^2}{\Phi} = \sqrt{N_0}\,, 
\end{gathered}
\end{equation*}
into the Gross--Pitaevskii equation~\eqref{eq:GPE}--\eqref{eq:GPECompact}, implies
\begin{equation*}
\ii \, e'(t) \, \Phi(x) = e(t) \, H\big(\Phi(x)\big)\,, \quad (x, t) \in \Omega \times [t_0, T]\,, 
\end{equation*}
and hence $e(t) = \ee^{- \, \ii \, \mu \, t}$ for $t \in [t_0, T]$ with $\mu \in \RR$ as well as
\begin{equation*}
H\big(\Phi(x)\big) = \mu \, \Phi(x)\,, \quad x \in \Omega\,.
\end{equation*}
Taking the inner product and performing integration-by-parts shows the close connection of the so-called chemical potential with the energy 
\begin{equation*}
\begin{gathered}
N_0 \, \mu = \scp{L^2}{H(\Phi)}{\Phi} = E(\Phi) + \tfrac{1}{2} \, E_2(\Phi) = E\big(\Psi(\cdot, t)\big) + \tfrac{1}{2} \, E_2\big(\Psi(\cdot, t)\big)\,, \\
t \in [t_0, T]\,,  
\end{gathered}
\end{equation*}
and in particular confirms that it is a real-valued quantity.
Summarised, the ground state solution at minimal energy level is given by a pair $(\Phi, \mu)$ that satisfies the relations 
\begin{equation}
\label{eq:GroundExcitedState}  
\begin{gathered}
\Psi(x, t) = \ee^{- \, \ii \, \mu \, t} \, \Phi(x)\,, \quad (x, t) \in \Omega \times [t_0, T]\,, \quad \norm{L^2}{\Phi} = \sqrt{N_0}\,, \\
H(\Phi) = \mu \, \Phi\,, \quad \mu = \tfrac{1}{N_0} \, \big(E(\Phi) + \tfrac{1}{2} \, E_2(\Phi)\big)\,. 
\end{gathered}
\end{equation}  
Accordingly, solutions of this form at higher energy levels are referred to as excited states. 

\MyParagraph{Imaginary time propagation}
The imaginary time propagation is a heuristic approach for the computation of ground and first excited state solutions.
It relies on the time integration of the parabolic counterpart of the Gross--Pitaevskii equation
\begin{subequations}
\label{eq:GPEImaginary}  
\begin{equation}
\label{eq:GPEImaginary1}  
\begin{cases}
\partial_t \wtPsi(x, t) = - \, H\big(\wtPsi(x, t)\big)\,, \\
\wtPsi(x, t_0) = \wtPsi_0(x)\,, \quad (x, t) \in \Omega \times [t_0, T]\,,
\end{cases}
\end{equation}
formally retained from~\eqref{eq:GPE}--\eqref{eq:GPECompact} by multiplying the time variable by the imaginary unit and using that $\ii \, \partial_t = - \, \partial_{\, \ii \, t}$.
For well-definedness, it is essential that the weights in the Laplacian are negative. 
In order to ensure convergence towards a non-trivial stationary solution, an additional normalisation condition is imposed, which is often realised by a straightforward scaling of the numerical solution values at certain intermediate times. 
For suitably chosen real-valued initial states and sufficiently large final times, it is expected that this approach yields appropriate approximations to ground and excited state solutions through 
\begin{equation}
\sqrt{N_0} \, \frac{\wtPsi(x, T)}{\norm{L^2}{\wtPsi(\cdot, T)}} \approx \Phi(x)\,, \quad x \in \Omega\,.
\end{equation}
\end{subequations}
More precisely, the resulting space-dependent function~$\Phi: \Omega \to \RR$ determines~$\mu \in \RR$ and hence $\Psi: \Omega \times [t_0, T] \to \CC$ as stated in~\eqref{eq:GroundExcitedState}.  
In brief, the imaginary time propagation is substantiated by the application of the gradient descent method for the computation of the solution with minimal energy
\begin{equation*}
\wtPsi_{\text{new}} = \wtPsi_{\text{old}} - \tau E'\big(\wtPsi_{\text{old}}\big) = \wtPsi_{\text{old}} - 2 \, \tau H\big(\wtPsi_{\text{old}}\big)
\end{equation*}
and its interpretation as explicit Euler approximation to~\eqref{eq:GPEImaginary1}.
Hereby, the connection between the directional derivative of the energy functional and the operator defining the right-hand side of the partial differential equation in~\eqref{eq:GPEImaginary1} is used, see~\eqref{eq:MassEnergyGateaux}.
Further enhancements of the procedure incorporate the knowledge of approximate ground and first excited state solutions in simplified settings and stepwise adaptations to the problem under consideration.
Specifically, for harmonic potentials, the first Hermite function 
\begin{equation}
\label{eq:Hermite0}
\begin{gathered}
\gamma_i = 0\,, \quad i \in \{1, \dots, d\}\,, \quad \vartheta = 0\,, \\
\Psi(x, t) = \tfrac{\sqrt{N_0}}{\sqrt[4]{\pi^d}} \, 
\prod_{i=1}^{d} \Big(\sqrt[8]{\tfrac{\beta_i}{- \, \alpha_i}} \, \ee^{- \, \ii \, t \, \sqrt{- \, \alpha_i \, \beta_i} - \frac{1}{2} \, \sqrt{\frac{\beta_i}{- \, \alpha_i}} \, x_i^2}\Big)\,, \quad x = (x_1, \dots, x_d) \in \RR^d\,, 
\end{gathered}
\end{equation}
and the Thomas--Fermi approximation
\begin{equation}
\label{eq:ThomasFermi}
\begin{gathered}
\alpha_i = 0\,, \quad \gamma_i = 0\,, \quad i \in \{1, \dots, d\}\,, \\
\Phi(x) = \begin{cases} \sqrt{W(x)}\,, & W(x) = \tfrac{1}{\vartheta} \, \big(\mu - \Vbetagammadelta(x)\big) > 0\,, \\ 0\,, & W(x) \leq 0\,, \end{cases} \\
\mu = \begin{cases}
\sqrt[3]{\tfrac{9}{16} \,  \, \beta_1 \, \vartheta^2 N_0^2}\,, &d = 1\,, \\
\sqrt{\tfrac{2}{\pi} \, \sqrt{\beta_1 \, \beta_2} \, \vartheta N_0}\,, &d = 2\,, \\
\sqrt[5]{\tfrac{225}{8 \, \pi} \, \beta_1 \, \beta_2 \, \beta_3 \, \vartheta^2 N_0^2}\,, &d = 3\,, \\
\end{cases} \\
\Psi(x, t) = \ee^{- \, \ii \, \mu \, t} \, \Phi(x)\,, \quad x \in \RR^d\,, 
\end{gathered}
\end{equation}
yield ground state solutions for the linear Schr{\"o}dinger equation and the limiting case of a negligible kinetic part, respectively.
Another option to improve the convergence properties of the algorithm is the application of the gradient descent with momentum, where~\eqref{eq:GPEImaginary1} is replaced by 
\begin{subequations}
\label{eq:GPEImaginaryMomentum} 
\begin{equation}
\begin{cases}
\partial_{tt} \wtPsi(x, t) = - \, c \, \partial_t \wtPsi(x, t) - H\big(\wtPsi(x, t)\big)\,, \\
\wtPsi(x, t_0) = \wtPsi_{10}(x)\,, \quad \partial_t \wtPsi(x, t_0) = \wtPsi_{20}(x)\,, \quad (x, t) \in \Omega \times [t_0, T]\,, 
\end{cases}
\end{equation}
or an associated first-order system such as
\begin{equation}
\begin{cases}
\partial_t \wtPsi_1(x, t) = \wtPsi_2(x, t)\,, \\
\partial_t \wtPsi_2(x, t) =- \, c \, \wtPsi_2(x, t) - H\big(\wtPsi_1(x, t)\big)\,, \\
\wtPsi_1(x, t_0) = \wtPsi_{10}(x)\,, \quad \wtPsi_2(x, t_0) = \wtPsi_{20}(x)\,, \quad (x, t) \in \Omega \times [t_0, T]\,, 
\end{cases}
\end{equation}
\end{subequations}
with suitably chosen initial conditions for the solution and its time derivative. 
\section{Gross--Pitaevskii systems}
\label{sec:GPS}
\MyParagraph{Model problem}
In order to exemplify systems of coupled time-dependent Gross--Pitaevskii equations describing the dynamical evolution of multi-species Bose--Einstein condensates, we consider the model problem  
\begin{subequations}
\label{eq:GPS}
\begin{equation}
\label{eq:GPS1}
\begin{cases}
&\ii \, \partial_t \Psi_j(x, t) = \Deltaalphaj \Psi_j(x, t) + \Vbetagammadeltaj(x) \, \Psi_j(x, t) \\
&\displaystyle \qquad\qquad\qquad + \sum_{k=1}^J \vartheta_{jk} \, \abs{\Psi_k(x, t)}^2 \, \Psi_j(x, t)\,, \\
&\Psi_j(x, t_0) = \Psi_{0j}(x)\,, \quad (x, t) \in \Omega \times [t_0, T]\,, \quad j \in \{1, \dots, J\}\,. 
\end{cases}
\end{equation}
By analogy to~\eqref{eq:GPE}, the vector-valued macroscopic wave function is defined on a certain space domain and time interval
\begin{equation}
\Psi: \Omega \times [t_0, T] \longrightarrow \CC^d: (x, t) \longmapsto \Psi(x, t) = \big(\Psi_1(x, t), \dots, \Psi_d(x, t)\big)^T\,, 
\end{equation}  
and the weighted Laplace operators and external potentials are given by
\begin{equation}
\begin{gathered}
\Deltaalphaj = \sum_{i=1}^{d} \alpha_{ji} \, \partial_{x_i}^2\,, \quad \alpha_j = (\alpha_{j1}, \dots, \alpha_{jd})\,, \\
\Vbetagammadeltaj: \Omega \longrightarrow \RR: x \longmapsto \sum_{i=1}^{d} \big(\beta_{ji} \, x_i^2 + \gamma_{ji} \, \sin^2(\delta_{ji} \, x_i) \Big)\,, \\
\beta_j = (\beta_{j1}, \dots, \beta_{jd})\,, \quad \gamma_j = (\gamma_{j1}, \dots, \gamma_{jd})\,, \quad \delta_j = (\delta_{j1}, \dots, \delta_{jd})\,, \\
\alpha_{ji}, \beta_{ji} > 0\,, \quad \gamma_{ji}, \delta_{ji} \in \RR\,, \quad i \in \{1, \dots, d\}\,, \quad j \in \{1, \dots, J\}\,.
\end{gathered}
\end{equation}  
\end{subequations}

\MyParagraph{Compact reformulation}
With regard to~\eqref{eq:GPECompact}, we rewrite~\eqref{eq:GPS1} in compact form as 
\begin{subequations}
\label{eq:GPSCompact}
\begin{equation}
\begin{cases}
&\ii \, \partial_t \Psi(x, t) = H\big(\Psi(x, t)\big) \\
&\Psi(x, t_0) = \Psi_0(x)\,, \quad (x, t) \in \Omega \times [t_0, T]\,, 
\end{cases}
\end{equation}
where we employ the convenient short notation
\begin{equation}
\begin{gathered}
H(\psi) = H_1(\psi) + H_2(\psi)\,, \\
\psi = (\psi_1, \dots, \psi_J)^T\,, \quad H_{\ell}(\psi) = \big(H_{\ell 1}(\psi), \dots, H_{\ell J}(\psi)\big)^T\,, \quad \ell \in \{1, 2\}\,, \\ 
H_{1j}(\psi) = \Deltaalphaj \psi_j + \Vbetagammadeltaj \, \psi_j\,, \quad H_{2j}(\psi) = \sum_{k=1}^J \vartheta_{jk} \, \abs{\psi_k}^2 \, \psi_j\,, \\
j \in \{1, \dots, J\}\,.
\end{gathered}
\end{equation}
\end{subequations}

\MyParagraph{Conserved quantities and ground states}
By means of these abbreviations, it is straightforward to extend our considerations for a single Gross--Pitaevskii equation to systems. 
That is, for regular solutions to~\eqref{eq:GPS}--\eqref{eq:GPSCompact}, total mass and energy 
\begin{equation*}
\begin{gathered}
M(\psi) = \norm{(L^2)^d}{\psi}^2\,, \quad E(\psi) = E_1(\psi) + \tfrac{1}{2} \, E_2(\psi)\,, \\
E_{\ell j}(\psi) = \scpbig{L^2}{H_{\ell j}(\psi)}{\psi_j}\,, \quad E_{\ell}(\psi) = \scpbig{(L^2)^d}{H_{\ell}(\psi)}{\psi} = \sum_{j=1}^{J} E_{{\ell}j}(\psi)\,, \\
\ell \in \{1, 2\}\,, \quad j \in \{1, \dots, J\}\,, 
\end{gathered}
\end{equation*}
are preserved over time 
\begin{equation*}
M\big(\Psi(\cdot, t)\big) = M\big(\Psi(\cdot, t_0)\big)\,, \quad E\big(\Psi(\cdot, t)\big) = E\big(\Psi(\cdot, t_0)\big)\,, \quad t \in [t_0, T]\,,
\end{equation*}
which justifies the normalisation condition
\begin{equation*}
\norm{(L^2)^d}{\Psi(\cdot, t)} = \sqrt{N_0}\,, \quad t \in [t_0, T]\,, 
\end{equation*}  
with a strictly positive real number $N_0 > 0$, see also~\eqref{eq:ScpNorm} and~\eqref{eq:MassEnergy}.
Ground and first excited state solutions are characterised by the relations
\begin{equation}
\label{eq:GPSGroundstate}  
\begin{gathered}
\Psi_j(x, t) = \ee^{- \, \ii \, \mu_j \, t} \, \Phi_j(x)\,, \quad (x, t) \in \Omega \times [t_0, T]\,, \quad \norm{(L^2)^d}{\Phi_j} = \sqrt{N_{0j}}\,, \\
H_{1j}(\Phi) + H_{2j}(\Phi) = \mu_j \, \Phi_j\,, \quad \mu_j = \tfrac{1}{N_{0j}} \, \big(E_{1j}(\Phi_j) + E_{2j}(\Phi_j)\big)\,, \\
N_{0j} > 0\,, \quad \sum_{j=1}^{j} N_{0j}^2 = N_0\,, \quad j \in \{1, \dots, J\}\,,
\end{gathered}
\end{equation}  
see~\eqref{eq:GroundExcitedState}.   
Accordingly, their numerical computation relies on the time integration of the non-reversible system 
\begin{subequations}
\label{eq:GPSImaginary}  
\begin{equation}
\begin{cases}
\partial_t \wtPsi(x, t) = - \, H\big(\wtPsi(x, t)\big)\,, \\
\wtPsi(x, t_0) = \wtPsi_0(x)\,, \quad (x, t) \in \Omega \times [t_0, T]\,,
\end{cases}
\end{equation}
under the normalisation condition
\begin{equation}
\sqrt{N_0} \, \frac{\wtPsi(x, T)}{\norm{(L^2)^d}{\wtPsi(\cdot, T)}} \approx \Phi(x)\,, \quad x \in \Omega\,,
\end{equation}
\end{subequations}
or a related second-order in time system, respectively, see~\eqref{eq:GPEImaginary} and~\eqref{eq:GPEImaginaryMomentum}. 
\section{Operator splitting methods}
\label{sec:Splitting}
\MyParagraph{General initial value problem}
In this section, to emphasise common features of our approaches for different situations, we introduce standard modified operator splitting methods for the general initial value problem 
\begin{equation}
\label{eq:IVPGeneral}
\begin{cases}
\partial_t U(x, t) = F\big(U(x, t)\big) = F_1\big(U(x, t)\big) + F_2\big(U(x, t)\big)\,, \\
U(x, t_0) \text{ given}\,, \quad (x, t) \in \Omega \times [t_0, T]\,. 
\end{cases}
\end{equation}
Subsequently, we include detailed information concerning their realisation for Gross--Pitaevskii equations and extended systems, see also~\eqref{eq:GPE}--\eqref{eq:GPECompact}, \eqref{eq:GPEImaginary}, \eqref{eq:GPEImaginaryMomentum}, \eqref{eq:GPS}--\eqref{eq:GPSCompact}, and~\eqref{eq:GPSImaginary}.

\MyParagraph{Time-stepping approach}
Due to the fact that the memory capacities needed in the time integration and imaginary time propagation of Gross--Pitaevskii systems in three space dimensions is limited, we find it most practicable to employ a time-stepping approach and to store a relatively low number of numerical solution values at once.
That is, we subdivide the considered time interval in appropriately chosen subintervals and determine numerical approximation values at certain intermediate times through a recursive procedure
\begin{equation*}
t_0 < t_1 < \dots < t_N = T\,, \quad U_n \approx U(\cdot, t_n)\,, \quad n \in \{0, 1, \dots, N\}\,.
\end{equation*}
Incorporating the evident advantages of an adaptive strategy based on a local error control, we permit non-uniform temporal grid points and denote the associated time increments by  
\begin{equation*}
\tau_n = t_{n+1} - t_n\,, \quad n \in \{0, 1, \dots, N-1\}\,.
\end{equation*}
In view of the statement of standard and modified operator splitting methods, we express the exact solutions to~\eqref{eq:IVPGeneral} and the associated subproblems 
\begin{equation}
\label{eq:IVPGeneralSubproblems}  
\begin{cases}
\partial_t u_j(x, t) = F_j\big(u_j(x, t)\big)\,, \\
u_j(x, t_0) \text{ given}\,, \quad (x, t) \in \Omega \times [t_0, T]\,,
\end{cases} \;
j \in \{1, 2\}\,,
\end{equation}
in terms of the corresponding evolution operators  
\begin{equation*}
\begin{gathered}  
\nE_{t-t_0, F}\big(U(\cdot, t_0)\big) = U(\cdot, t)\,, \\
\nE_{t-t_0, F_j}\big(u_j(\cdot, t_0)\big) = u_j(\cdot, t)\,, \quad j \in \{1, 2\}\,, \\
t \in [t_0, T]\,, 
\end{gathered}  
\end{equation*}
The analogous relations for the numerical solution values in terms of the respective splitting operators read as 
\begin{equation*}
U_{n+1} = \nS_{\tau_n, F}(U_n) \approx U(\cdot, t_{n+1}) = \nE_{\tau_n, F}\big(U(\cdot, t_n)\big)\,, \quad n \in \{0, 1, \dots, N-1\}\,. 
\end{equation*}
For a juxtaposition with modified operator splitting methods, we next include the general format of standard operator splitting methods and specify example methods.

\MyParagraph{Standard splitting methods}
Operator splitting methods for the time integration of nonlinear evolution equations make use of the fact that the defining operators naturally decompose into two or more parts and that the numerical approximation of the associated subproblems is significantly simpler compared to the numerical approximation of the original problem, see~\eqref{eq:IVPGeneral} and~\eqref{eq:IVPGeneralSubproblems}.
Any standard high-order splitting method is given by a composition of the form 
\begin{equation*}
\begin{gathered}
U_{n+1} = \nS_{\tau_n, F}(U_n) \approx U(\cdot, t_{n+1}) = \nE_{\tau_n, F}\big(U(\cdot, t_n)\big)\,, \\
\nS_{\tau_n, F} = \nE_{\tau_n, b_s F_2} \circ \nE_{\tau_n, a_s F_1} \circ \dots \circ \nE_{\tau_n, b_1 F_2} \circ \nE_{\tau_n, a_1 F_1}\,, \\
n \in \{0, 1, \dots, N-1\}\,, 
\end{gathered}
\end{equation*}  
with suitably chosen real coefficients $(a_i, b_i)_{i=1}^{s}$.
Prominant instances are the first-order Lie--Trotter splitting method involving a single composition per step 
\begin{subequations}
\label{eq:Coefficients}
\begin{equation}
\label{eq:SchemeOrder1}
s = 1\,, \quad a_1 = 1\,, \quad b_1 = 1\,, 
\end{equation}
the second-order Strang splitting method based on a symmetric composition
\begin{equation}
\label{eq:SchemeOrder2}
s = 2\,, \quad a_1 = 0\,, \quad a_2 = 1\,, \quad b_1 = \tfrac{1}{2} = b_2\,, 
\end{equation}
and the symmetric fourth-order splitting method by \textsc{Yoshida}~\MyCite{Yoshida1990} including four compositions  
\begin{equation}
\label{eq:CoefficientsYoshida}
\begin{gathered}
s = 4\,, \quad a_1 = 0\,, \quad a_2 = 1 - 2 \, b_2 = a_4\,, \quad a_3 = 4 \, b_2 - 1\,, \\
b_1 = \tfrac{1}{2} - b_2 = b_4\,, \quad b_2 = \tfrac{1}{6} \, \big(1 - \sqrt[3]{2} - \tfrac{1}{2} \sqrt[3]{4}\big) = b_3\,.
\end{gathered}
\end{equation}
\end{subequations}
A variety of higher-order operator splitting methods are found in the literature, see for instance~\MyCite{BlanesCasasMurua2024,McLachlanQuispel2002}.
For our numerical experiments, we select an optimised fourth-order splitting method proposed in~\MyCite{BlanesMoan2002} as a renowned instance with excellent properties regarding stability, accuracy, efficiency, and energy preservation over long times.

\MyParagraph{Modified operator splitting methods}
As already indicated in the introduction, higher-order operator splitting methods such as~\eqref{eq:CoefficientsYoshida} necessarily involve negative coefficients.
Because of severe stability issues in the context of non-reversible systems, they are thus inappropriate for ground and excited state computations.
A remedy to these shortcomings of standard operator splitting methods is the application of modified operator splitting methods, which can be cast into the general form 
\begin{subequations}
\label{eq:ModifiedSplitting}
\begin{equation}
\begin{gathered}
U_{n+1} = \nS_{\tau_n, F}(U_n) \approx U(\cdot, t_{n+1}) = \nE_{\tau_n, F}\big(U(\cdot, t_n)\big)\,, \\
\nS_{\tau_n, F} = \nE_{\tau_n, b_s F_2 + c_s \tau_n^2 \, G} \circ \nE_{\tau_n, a_s F_1} \circ \dots \circ \nE_{\tau_n, b_1 F_2 + c_1 \tau_n^2 \, G} \circ \nE_{\tau_n, a_1 F_1}\,, \\
n \in \{0, 1, \dots, N-1\}\,,
\end{gathered}
\end{equation}  
with certain real coefficients $(a_i, b_i, c_i)_{i=1}^{s}$.
More precisely, this kind of splitting methods relies on the efficient numerical approximation of the subproblems 
\begin{equation}
\label{eq:ModifiedSplitting2}
\begin{gathered}
\partial_t u_1(x, t) = a_i \, F_1\big(u_1(x, t)\big)\,, \quad u_1(x, t_n) \text{ given}\,, \\
\partial_t u_2(x, t) = b_i \, F_2\big(u_2(x, t)\big) + c_i \, \tau_n^2 \, G\big(u_2(x, t)\big)\,, \quad u_2(x, t_n) \text{ given}\,, \\
(x, t) \in \Omega \times [t_n, t_n + \tau_n]\,, \quad i \in \{1, \dots, s\}\,.
\end{gathered}
\end{equation}
Hereby, it is remarkable that the incorporation of the operator 
\begin{equation}
\label{eq:IteratedCommutator}
\begin{split}
G(v) &= F_1''(v) \, F_2(v) \, F_2(v) + F_1'(v) \, F_2'(v) \, F_2(v) + F_2'(v) \, F_2'(v) \, F_1(v) \\
&\quad\, - F_2''(v) \, F_1(v) \, F_2(v) - 2 \, F_2'(v) \, F_1'(v) \, F_2(v)
\end{split}
\end{equation}
\end{subequations}
permits the design of schemes such that the principal coefficients $(a_i, b_i)_{i=1}^{s}$ are non-negative.
In conjunction with the additional factor~$\tau_n^2$, this ensures stability for reversible as well as non-reversible systems under moderate smallness assumptions on the time increments. 
The arising derivatives are again understood as G{\^a}teaux derivatives, see~\eqref{eq:Gateaux}.

\MyParagraph{Historical context and novel aspects}
The alternative approach of modified operator splitting methods goes back to a seminal work of~\textsc{Chin}~\MyCite{Chin1997}, where so-called modified potential or force-gradient operator splitting methods were proposed in the context of \emph{linear} partial differential equations. 
In our recent contribution~\MyCite{BlanesCasasGonzalezThalhammer2024}, we could extend this class of methods within a nonlinear setting such that the operator defined in~\eqref{eq:IteratedCommutator} results as iterated commutator $[D_{F_2}, [D_{F_2}, D_{F_1}]]$ of the associated Lie-derivatives. 
In the present work, we generalise the underlying idea to incorporate commutators of the defining operators to the time integration and imaginary time propagation of Gross--Pitaevskii systems.
We specify~\eqref{eq:IteratedCommutator} for systems comprising two coupled time-dependent Gross--Pitaevskii equations~\eqref{eq:GPS}--\eqref{eq:GPSCompact} and their parabolic counterparts~\eqref{eq:GPSImaginary}, which are of particular relevance in view of actual laboratory set-ups, see for instance~\MyCite{RiboliModugno2002,ThalhammerGEtAl2008}.

\MyParagraph{Iterated commutators ($J = 2$)}
As important objectives regarding the practical implementation of modified operator splitting methods~\eqref{eq:ModifiedSplitting}, we next state the iterated commutators for Gross--Pitaevskii systems involving two coupled equations.
For this purpose, we use the following abbreviations for the component functions and the weighted gradients 
\begin{equation*}
\begin{gathered}
J = 2\,, \quad \wtPsi = (\wtPsi_1, \wtPsi_2)^T\,, \quad \Psi = (\Psi_1, \Psi_2) ^T\,, \quad G = (G_1, G_2) ^T\,, \\
\nabla = (\partial_{x_1}, \dots, \partial_{x_d})^T\,, \quad  \Nablaalphaj = (\alpha_{j1} \, \partial_{x_1}, \dots, \alpha_{jd} \, \partial_{x_d})^T\,, \quad V_j = \Vbetagammadeltaj\,, \\
j \in \{1, 2\}\,, 
\end{gathered}
\end{equation*}
see also~\eqref{eq:GPS} and \eqref{eq:GPSImaginary}.
\begin{enumerate}[(i)]
\item 
We begin with the study of the parabolic-type case~\eqref{eq:GPSImaginary}, where the first component of the operator in~\eqref{eq:IteratedCommutator} is given by
\begin{subequations}
\label{eq:CommutatorGPSImaginary}  
\begin{equation}
G_1(\wtPsi_1, \wtPsi_2) = 2 \, \big(S_0 + S_1(\wtPsi_1, \wtPsi_2) + S_2(\wtPsi_1, \wtPsi_2) + S_3(\wtPsi_1, \wtPsi_2)\big)\,, 
\end{equation}
\begin{equation}
S_0 = \Nablaalphaone V_1 \cdot \nabla V_1 \, \wtPsi_1\,,
\qquad\qquad\qquad\qquad\qquad\qquad\qquad\qquad\qquad\quad
\end{equation}
\begin{equation}
\begin{split}
S_1(\wtPsi_1, \wtPsi_2) 
&= 6 \, \vartheta_{11} \, V_1 \, \Nablaalphaone \wtPsi_1 \cdot \nabla \wtPsi_1 \, \wtPsi_1 \\
&\quad + 2 \, \vartheta_{12} \, V_2 \, \Big(\nabla \wtPsi_2 \cdot \big(2 \, \Nablaalphaone \wtPsi_1 \, \wtPsi_2 + \Nablaalphaone \wtPsi_2 \, \wtPsi_1\big) \\
&\quad + \big(\Deltaalphaone \wtPsi_2 - \Deltaalphatwo \wtPsi_2\big) \, \wtPsi_1 \, \wtPsi_2\Big) \\
&\quad + 2 \, \Nablaalphaone V_1 \cdot \big(3 \, \vartheta_{11} \, \nabla \wtPsi_1 \, \wtPsi_1 + 2 \, \vartheta_{12} \, \nabla \wtPsi_2 \, \wtPsi_2\big) \, \wtPsi_1 \\
&\quad + 2 \, \vartheta_{12} \, \Big(\Nablaalphaone V_2 \cdot \nabla \wtPsi_1 \, \wtPsi_2 \\
&\quad + 2 \, \big(\Nablaalphaone V_2 - \Nablaalphatwo V_2\big) \cdot \nabla \wtPsi_2 \, \wtPsi_1\Big) \, \wtPsi_2 \\
&\quad - \vartheta_{11} \, \Deltaalphaone V_1 \, \wtPsi_1^3 
+ \vartheta_{12} \, \big(\Deltaalphaone V_2 - 2 \, \Deltaalphatwo V_2\big) \, \wtPsi_1 \, \wtPsi_2^2\,, 
\quad\;\,
\end{split}
\end{equation}
\begin{equation}
\begin{split}
S_2(\wtPsi_1, \wtPsi_2)
&= 2 \, \Big(6 \, \vartheta_{11}^2 \, \Nablaalphaone \wtPsi_1 \cdot \nabla \wtPsi_1 \, \wtPsi_1^3 \\
&\quad + 3 \, \big(\vartheta_{11} \, \vartheta_{12} + \vartheta_{12} \, \vartheta_{21}\big) \, \Nablaalphaone \wtPsi_1 \cdot \nabla \wtPsi_1 \, \wtPsi_1 \, \wtPsi_2^2 \\
&\quad - 2 \, \vartheta_{12} \, \vartheta_{21} \, \Nablaalphatwo \wtPsi_1 \cdot \nabla \wtPsi_1 \, \wtPsi_1 \, \wtPsi_2^2 \\
&\quad
+ 4 \, \vartheta_{12} \, \vartheta_{22} \, \Nablaalphaone \wtPsi_1 \cdot \nabla \wtPsi_2 \, \wtPsi_2^3 \\
&\quad + 6 \, \big(\vartheta_{11} \, \vartheta_{12} + \vartheta_{12} \, \vartheta_{21}\big) \, \Nablaalphaone \wtPsi_1 \cdot \nabla \wtPsi_2 \, \wtPsi_1^2 \, \wtPsi_2 \\
&\quad - 4 \, \vartheta_{12} \, \vartheta_{21} \, \Nablaalphatwo \wtPsi_1 \cdot \nabla \wtPsi_2 \, \wtPsi_1^2 \, \wtPsi_2 \\
&\quad + \big(\vartheta_{12} \, \vartheta_{21} - \vartheta_{11} \, \vartheta_{12}\big) \, \Nablaalphaone \wtPsi_2 \cdot \nabla \wtPsi_2 \, \wtPsi_1^3 \\
&\quad + 2 \, \vartheta_{12}^2 \, \Nablaalphaone \wtPsi_2 \cdot \nabla \wtPsi_2 \, \wtPsi_1 \, \wtPsi_2^2 \\
&\quad + 6 \, \vartheta_{12} \, \vartheta_{22} \, \big(\Nablaalphaone - \Nablaalphatwo\big) \wtPsi_2 \, \cdot \nabla \wtPsi_2 \, \wtPsi_1 \, \wtPsi_2^2\Big)\,, 
\qquad\quad\;\;\,
\end{split}
\end{equation}
\begin{equation}
\begin{split}
S_3(\wtPsi_1, \wtPsi_2)
&= 2 \, \Big(2 \, \vartheta_{12} \, \vartheta_{21} \, \big(\Deltaalphaone - \Deltaalphatwo\big) \wtPsi_1 \, \wtPsi_1^2 \, \wtPsi_2^2 \\
&\quad + \big(\vartheta_{11} \, \vartheta_{12} - \vartheta_{12} \, \vartheta_{21}\big) \, \big(\Deltaalphatwo - \Deltaalphaone\big) \wtPsi_2 \, \wtPsi_1^3 \, \wtPsi_2 \\
&\quad + 2 \, \vartheta_{12} \, \vartheta_{22} \, \big(\Deltaalphaone - \Deltaalphatwo\big) \wtPsi_2 \, \wtPsi_1 \, \wtPsi_2^3\Big)\,.
\qquad\qquad\qquad\;\;\,
\end{split}
\end{equation}
\end{subequations}
It is remarkable that only~$S_0$ remains in the simplified linear case.
The first and second spatial derivatives of the potentials arising in~$S_1$ are often known explicitly.
The term~$S_3$ comprises second spatial derivatives of the solution components, whenever the weights in the Laplacians are different.
For the second component, we have to exchange the index one by two and two by one, respectively.
\item
The corresponding identity for~\eqref{eq:GPS} reads as 
\begin{subequations}
\label{eq:CommutatorGPS}  
\begin{equation}
G_1(\Psi_1, \Psi_2) = 2 \, \ii \, \big(S_0 + S_1(\Psi_1, \Psi_2) + S_2(\Psi_1, \Psi_2) + S_3(\Psi_1, \Psi_2)\big) \, \Psi_1\,, 
\end{equation}
\begin{equation}
S_0 = \Nablaalphaone V_1 \cdot \nabla V_1\,,
\qquad\qquad\qquad\qquad\qquad\qquad\qquad\qquad\qquad\qquad\;
\end{equation}
\begin{equation}
\begin{split}  
S_1(\Psi_1, \Psi_2)
&= 2 \, \Big(2 \, \vartheta_{12} \, \big(\Nablaalphaone V_1 - \Nablaalphatwo V_2\big) \cdot \Re\big(\nabla \Psi_2 \, \overline{\Psi_2}\big) \\
&\quad - \big(\vartheta_{11} \, \Deltaalphaone V_1 \, \abs{\Psi_1}^2 + \vartheta_{12} \, \Deltaalphatwo V_2 \, \abs{\Psi_2}^2\big)\Big)\,, 
\qquad\qquad\;\;\;
\end{split}
\end{equation}
\begin{equation}
\begin{split}  
S_2(\Psi_1, \Psi_2)
&= - \, 2 \, \Big(\vartheta_{11}^2 \, \Re\big(\Nablaalphaone \Psi_1 \cdot \nabla \Psi_1 \, \overline{\Psi_1}^2\big) \\
&\quad + 3 \, \vartheta_{11}^2 \, \Nablaalphaone \Psi_1 \cdot \overline{\nabla \Psi_1} \, \abs{\Psi_1}^2 \\
&\quad + 2 \, \vartheta_{12} \, \vartheta_{21} \, \Nablaalphatwo \Psi_1 \cdot \overline{\nabla \Psi_1} \, \abs{\Psi_2}^2\,, \\
&\quad + 2 \, \vartheta_{12} \, \vartheta_{21} \, \Re\big(\Nablaalphatwo \Psi_1 \cdot \nabla \Psi_2 \, \overline{\Psi_1} \, \overline{\Psi_2} \\
&\quad + \Nablaalphatwo \Psi_1 \cdot \nabla \overline{\Psi_2} \, \overline{\Psi_1} \, \Psi_2\big)\,, \\
&\quad - \vartheta_{12}^2 \, \Big(\Re\big(\Nablaalphaone \Psi_2 \cdot \nabla \Psi_2 \, \overline{\Psi_2}^2\big) \\
&\quad + \Nablaalphaone \Psi_2 \cdot \nabla \overline{\Psi_2} \, \abs{\Psi_2}^2\Big) \\
&\quad + 2 \, \vartheta_{12} \, \vartheta_{22} \, \Big(\Re\big(\Nablaalphatwo \Psi_2 \cdot \nabla \Psi_2 \, \overline{\Psi_2}^2\big) \\
&\quad + 2 \, \Nablaalphatwo \Psi_2 \, \nabla \overline{\Psi_2} \, \abs{\Psi_2}^2\Big) \\
&\quad + 2 \, \vartheta_{11} \, \vartheta_{12}\, \Nablaalphaone \Psi_2 \, \nabla \overline{\Psi_2} \, \abs{\Psi_1}^2\Big)\,, 
\qquad\qquad\qquad\qquad\;\;\;
\end{split}
\end{equation}
\begin{equation}
\begin{split}  
S_3(\Psi_1, \Psi_2)
&= - \, 2 \, \Big(2 \, \Re\big(\vartheta_{11}^2 \, \Deltaalphaone \Psi_1 \, \overline{\Psi_1} \, \abs{\Psi_1}^2 \\
&\quad + \vartheta_{12} \, \vartheta_{21} \, \Deltaalphatwo \Psi_1 \, \overline{\Psi_1} \, \abs{\Psi_2}^2\big)\,, \\
&\quad + 2 \, \Re\big(\vartheta_{11} \, \vartheta_{12} \, \Deltaalphaone \Psi_2 \, \overline{\Psi_2} \, \abs{\Psi_1}^2 \\
&\quad + \vartheta_{12} \, \vartheta_{22} \, \Deltaalphatwo \Psi_2 \, \overline{\Psi_2} \, \abs{\Psi_2}^2\big)\Big)\,.
\qquad\qquad\qquad\qquad\quad\;\;\;\,
\end{split}
\end{equation}
\end{subequations}
The analogous identity holds for the second component, but with the index one replaced by two and two by one, respectively.
Again, in the simplified linear case, only the term~$S_0$ remains.
Usually, the first and second spatial derivatives of the potentials that are present in~$S_1$ can be calculated explicitly, whereas the first and second spatial derivatives of the solution components have to be computed numerically. 
However, hereby it is most remarkable that in the first component a common factor~$\Psi_1$ and in the second component a common factor~$\Psi_2$ can be extracted.
Together with the validity of an invariance principle stated below, this special structure of the iterated commutator provides the knowledge of the exact solutions to the nonlinear subproblems on each of the temporal subintervals and hence considerably enhances the efficiency of modified operator splitting methods for time-dependent Gross--Pitaevskii systems. 
\end{enumerate}

\MyParagraph{Iterated commutator ($J = 1$)}
By setting $\alpha = (\alpha_{11}, \dots, \alpha_{1d})$, $(\alpha_{21}, \dots, \alpha_{2d}) = 0$, $\vartheta = \vartheta_{11}$, $\vartheta_{12} = \vartheta_{21} = \vartheta_{22} = 0$, $V = V_1 = \Vbetagammadelta$, and $V_2 = 0$ in~\eqref{eq:CommutatorGPSImaginary} as well as~\eqref{eq:CommutatorGPS}, the corresponding identities for a single Gross--Pitaevskii equation result.
\begin{enumerate}[(i)]
\item 
On the one hand, in the context of the non-reversible equation~\eqref{eq:GPEImaginary}, we retain the significant simplification 
\begin{equation*}
\begin{split}
G(\wtPsi)
&= 2 \, \Big(\Nablaalpha V \cdot \nabla V + \vartheta \, \big(6 \, V \, \Nablaalpha \wtPsi \cdot \nabla \wtPsi + 6 \, \Nablaalpha V \cdot \nabla \wtPsi \, \wtPsi - \, \Deltaalpha V \, \wtPsi^2\big) \\
&\qquad\;\; + 12 \, \vartheta^2 \, \Nablaalpha \wtPsi \cdot \nabla \wtPsi \, \wtPsi^2\Big) \, \wtPsi\,.
\end{split}
\end{equation*}
\item 
On the other hand, in the context of the time-dependent Gross--Pitaevskii equation~\eqref{eq:GPE}, we instead retain 
\begin{equation*}
\begin{split}
G(\Psi)
&= 2 \, \ii \, \Big(\Nablaalpha V \cdot \nabla V - 2 \, \vartheta \, \Deltaalphaone V \, \abs{\Psi}^2 - 2 \, \vartheta^2 \, \Re\big(\Nablaalpha \Psi \cdot \nabla \Psi \, \overline{\Psi}^2\big) \\ 
&\qquad\quad 
- 6 \, \vartheta^2 \, \Nablaalpha \Psi \cdot \overline{\nabla \Psi} \, \abs{\Psi}^2 - 4 \, \vartheta^2 \, \Re\big(\Deltaalpha \Psi \, \overline{\Psi}\big) \, \abs{\Psi}^2\Big) \, \Psi\,.
\end{split}
\end{equation*}
\end{enumerate}
For the special case $\alpha = (1, \dots, 1) \in \RR^d$, both relations are consistent with equations~(14) and~(15) in~\MyCite{BlanesCasasGonzalezThalhammer2024}.

\MyParagraph{Invariance principle}
As indicated in the introduction and shortly before, an invariance principle permits the replacement of the nonlinear subproblems associated with time-dependent Gross--Pitaevskii systems by simplified linear systems with known exact solutions.
In view of a compact formulation of this fundamental result, we make use of the fact that the defining operators can be rewritten as 
\begin{equation*}
H_2(\psi) = \wtH_2(\psi) \, \psi\,, \quad G(\psi) = \wtG(\psi) \, \psi\,. 
\end{equation*}
Moreover, neglecting for notational simplicity the dependencies on the method coefficients $b_i, c_i \in \RR$ and the time increment $\tau_n > 0$, we set
\begin{equation*}
\begin{gathered}
\wtF_3(\psi) = - \, \ii \, \big(b_i \, \wtH_2(\psi) + \tau_n^2 \, c_i \, \wtG(\psi)\big)\,, \quad i \in \{1, \dots, s\}\,, 
\end{gathered}
\end{equation*}
see also~\eqref{eq:ModifiedSplitting2}.
Our finding is that the solution to the initial value problem
\begin{equation*}
\begin{cases}
&\displaystyle     
\partial_t \Psi(x, t) = \wtF_3\big(\Psi(x, t)\big) \, \Psi(x, t)\,, \\ 
&\Psi(x, t_n)  \text{ given}\,, \quad (x, t) \in \Omega \times [t_n, t_n + \tau_n]\,, 
\end{cases}
\end{equation*}
fulfills the simplified linear problem
\begin{equation*}
\begin{cases}
&\displaystyle     
\partial_t \Psi(x, t) = \wtF_3\big(\Psi(x, t_n)\big) \, \Psi(x, t)\,, \\ 
&\Psi(x, t_n)  \text{ given}\,, \quad (x, t) \in \Omega \times [t_n, t_n + \tau_n]\,,
\end{cases}
\end{equation*}
and hence is given by the explicit representation
\begin{equation}
\label{eq:Invariance}  
\Psi(x, t) = \ee^{(t - t_n) \, \wtF_3(\Psi(x, t_n))} \, \Psi(x, t_n)\,, \quad (x, t) \in \Omega \times [t_n, t_n + \tau_n]\,. 
\end{equation}
A rigorous derivation of this relation relies on the calculation of the G{\^a}teaux derivative and the verification of the resulting identity 
\begin{equation*}
\begin{split}
&\tfrac{\dd}{\dd t} \wtF_3\big(\Psi(\cdot, t)\big)
= \wtF_3'\big(\Psi(\cdot, t)\big) \, \partial_t \Psi(\cdot, t) \\
&\quad = \wtF_3'\big(\Psi(\cdot, t)\big) \, \Big(\wtF_3\big(\Psi(x, t)\big) \, \Psi(x, t)\Big) = 0\,, \quad t \in [t_n, t_n + \tau_n]\,,
\end{split}
\end{equation*}
see also~\eqref{eq:Gateaux} and~\eqref{eq:MassEnergy}.
We here omit the technical details and instead refer to confirming numerical experiments.

\MyParagraph{Practical implementation}
Modified operator splitting methods for the time integration and imaginary time propagation of Gross--Pitaevskii systems make use of the fact that both, the numerical approximation of the associated linear subproblems involving weighted Laplacians and the subproblems comprising the potentials and nonlinearities can be realised in a highly efficient manner by Fast Fourier transforms and their inverses, respectively, als well as pointwise multiplications.%
\footnote{An elementary \textsc{Matlab} code that illustrates the realisation of a specific modified operator splitting method for a single Gross--Pitaevskii equation is available through the publicly accessible link \MyUrl{doi.org/10.5281/zenodo.7945624}.
The extension to two-component systems involves a certain amount of technical details.
However, in essence, it requires the implementation of~\eqref{eq:CommutatorGPSImaginary} and~\eqref{eq:CommutatorGPS}.  
}
\begin{enumerate}[(i)]
\item
\emph{Fourier series representations.} \,   
Employing the reasonable assumption that the prescribed initial states and hence the solutions to~\eqref{eq:GPS}--\eqref{eq:GPSCompact} are localised and sufficiently regular, we may replace the underlying domains by Cartesian products of suitably chosen symmetric intervals 
\begin{equation*}
\Omega = \prod_{i=1}^d \, [- \, \omega_i, \omega_i]\,, \quad \omega_i > 0\,, \quad i \in \{1, \dots, d\}\,. 
\end{equation*}
Consequently, denoting by $(\nF_m)_{m \in \ZZ^d}$ the periodic Fourier functions
\begin{equation*}
\begin{gathered}
\nF_m: \RR^d \longrightarrow \CC: x = (x_1, \dots, x_d) \longmapsto \prod_{i=1}^{d} \big(\tfrac{1}{\sqrt{2 \, \omega_i}} \, \ee^{\, \ii \, \pi \, m_i \, (x_i/\omega_i + 1)}\big)\,, \\
m = (m_1, \dots, m_d) \in \ZZ^d\,, 
\end{gathered}  
\end{equation*}
and using that they form a complete orthonormal system of the Hilbert space~$L^2(\Omega, \CC)$, Fourier series representations and Parseval's identities hold true 
\begin{equation*}
\begin{gathered}
\psi = \sum_{m \in \ZZ^d} \psi_m \, \nF_m\,, \quad \psi_m = \int_{\Omega} \psi(x) \, \nF_{-m}(x) \; \dd x\,, \quad m \in \ZZ^d\,, \\
\norm{L^2}{\psi}^2 = \sum_{m \in \ZZ^d} \abs{\psi_m}^2\,, \quad \psi \in L^2(\Omega, \CC)\,.
\end{gathered}  
\end{equation*}
Moreover, by means of the following abbreviations for the purely imaginary eigenvalues associated with first spatial derivatives as well as the corresponding non-positive eigenvalues of second derivatives
\begin{equation*}
\begin{gathered}
\partial_{x_j} \, \nF_m = \mu_{m_j} \, \nF_m\,, \quad \mu_{m_j} = \tfrac{\ii \, \pi \, m_j}{\omega_j} \in \ii \, \RR\,, \\
\mu_m = \big(\mu_{m_j}\big)_{j=1}^d \in \CC^{d \times 1}\,, \quad \mualpham = \big(\alpha_j \, \mu_{m_j}\big)_{j=1}^d \in \CC^{d \times 1}\,, \\
\nabla \nF_m = \mu_m \, \nF_m\,, \quad \Nablaalpha \, \nF_m = \mualpham \, \nF_m\,, \\
\partial_{x_j}^2 \, \nF_m = \lambda_{m_j} \, \nF_m\,, \quad \lambda_{m_j} = \mu_{m_j}^2 = - \, \tfrac{\pi^2 \, m_j^2}{\omega_j^2} \in \RR_{\leq 0}\,, \\
\lambda_m = \sum_{j=1}^{d} \lambda_{m_j} \in \RR_{\leq 0}\,, \quad \lambdaalpham = \sum_{j=1}^{d} \alpha_j \, \lambda_{m_j} \in \RR_{\leq 0}\,, \\
\Delta \, \nF_m = \lambda_m \, \nF_m\,, \quad \Deltaalpha \, \nF_m = \lambdaalpham \, \nF_m\,, \\
m = (m_1, \dots, m_d) \in \ZZ^d\,, \quad x = (x_1, \dots, x_d) \in \Omega\,, 
\end{gathered}  
\end{equation*}  
we obtain series representations such as 
\begin{equation*}
\nabla \psi = \sum_{m \in \ZZ^d} \mu_m \, \psi_m \, \nF_m\,, \quad \Delta \psi = \sum_{m \in \ZZ^d} \lambda_m \, \psi_m \, \nF_m\,,
\end{equation*}  
and their analogues for weighted gradients and Laplacians.
\item 
\emph{Derivatives and linear subproblems.} \,   
The above identities permit to determine on the one hand the spatial derivatives of the current values of the time-discrete solutions arising in the iterated commutators and on the other hand the solutions to the linear subproblems, which can be cast into the form
\begin{equation*}
\begin{gathered}
\tfrac{\dd}{\dd t} \, u(x, t) = C \Deltaalpha u(x, t) \,, \quad u(x, t_n) \text{ given}\,, \\
(x, t) \in \Omega \times [t_n, t_n + \tau_n]\,, \quad C \in \CC\,, \quad \alpha \in \RR^d\,, \\
u(\cdot, t_{n+1}) = \nE_{\tau_n, C \Deltaalpha}\bigg(\sum_{m \in \ZZ^d} u_m(\cdot, t_n) \, \nF_m\bigg) = \sum_{m \in \ZZ^d} \ee^{\, C \, \tau_n \, \lambdaalpham} \, u_m(\cdot, t_n) \, \nF_m\,, \\
\end{gathered}  
\end{equation*}
with high accuracy and efficiency by fast Fourier transforms, pointwise multiplications, and inverse fast Fourier transforms.
In general, these procedures constitute the computationally most expensive components. 
\item
\emph{Nonlinear subproblems (Time integration).} \,   
For Schr{\"o}dinger-type systems as they arise in the time integration of~\eqref{eq:GPS}--\eqref{eq:GPSCompact}, the invariance principle stated above permits to determine the true solutions to the nonlinear subproblems by pointwise multiplications, see~\eqref{eq:Invariance}.  
\item
\emph{Nonlinear subproblems (Imaginary time propagation).} \,   
For parabolic-type systems~\eqref{eq:GPSImaginary}, numerical approximations are obtained by the application of explicit Runge--Kutta methods.
More concretely, in order to enhance the efficiency behaviour of fourth-order modified operator methods, we combine the second-order Strang splitting method yielding 
\begin{equation*}
\nE_{\frac{1}{2} \tau_n, b_i F_2} \circ \nE_{\tau_n, c_i \tau_n^2 G_2} \circ \nE_{\frac{1}{2} \tau_n, b_i F_2} \approx \nE_{\tau_n, b_i F_2 + c_i \tau_n^2 G}\,, \quad i \in \{1, \dots, s\}\,, 
\end{equation*}
with the first-order explicit Euler method for~$\nE_{\tau_n, c_i \tau_n^2 G_2}$ and a fourth-order Runge--Kutta method~$\nE_{\frac{1}{2} \tau_n, b_i F_2}$, respectively. 
Whenever appropriate, for instance in the case of a single Gross--Pitaevskii equation, we may also use explicit solution representations for $\nE_{\frac{1}{2} \tau_n, b_i F_2}$. 
\end{enumerate}

\MyParagraph{Example method}
In our numerical experiments, we focus a fourth-order modified operator splitting method involving a low number of compositions and a single evaluation of the iterated commutator 
\begin{equation}
\label{eq:SchemeOrder4}  
\begin{gathered}
s = 3\,, \quad a_1 = 0\,, \quad a_2 = \tfrac{1}{2} = a_3\,, \quad b_1 = \tfrac{1}{6} = b_3\,, \quad b_2 = \tfrac{2}{3}\,, \\
c_1 = 0 = c_3\,, \quad c_2 = - \, \tfrac{1}{72}\,.
\end{gathered}
\end{equation}
This scheme goes back to a modified potential operator splitting method that was introduced by \textsc{Chin}~\MyCite{Chin1997} in the context of linear problems. 
Here, it is remarkable that the principal coefficients $(a_i, b_i)_{i=1}^{s}$ are non-negative and that the additional factor~$\tau_n^2$ ensures stability of the method for non-reversible systems under moderate smallness assumptions on the time increments. 
It is a common assessment that the validity of nonstiff order conditions deduced for linear ordinary differential equations and the presumption of sufficiently regular initial states and hence solutions is adequate to retain the full order of convergence in substantially more complex settings of nonlinear partial differential equations. 
However, it remains an open question to rigorously justify this implication in the lines of~\MyCite{Thalhammer2008,Thalhammer2012}.

\MyParagraph{Adaptivity in time}
With regard to the evident advantages of adaptivity in time, also observed in our numerical experiments, we combine the fourth-order modified operator splitting method~\eqref{eq:SchemeOrder4} with the second-order Strang splitting method.
Due to the fact that both schemes are suitable for both, the time integration and imaginary time propagation of Gross--Pitaevskii systems, this simple approach for a local error control permits to automatically adapt the time increments and enhances stability, accuracy, and efficiency.
\section{Numerical experiments}
\label{sec:Experiments}
\begin{figure}
\begin{center}
\includegraphics[width=4.2cm]{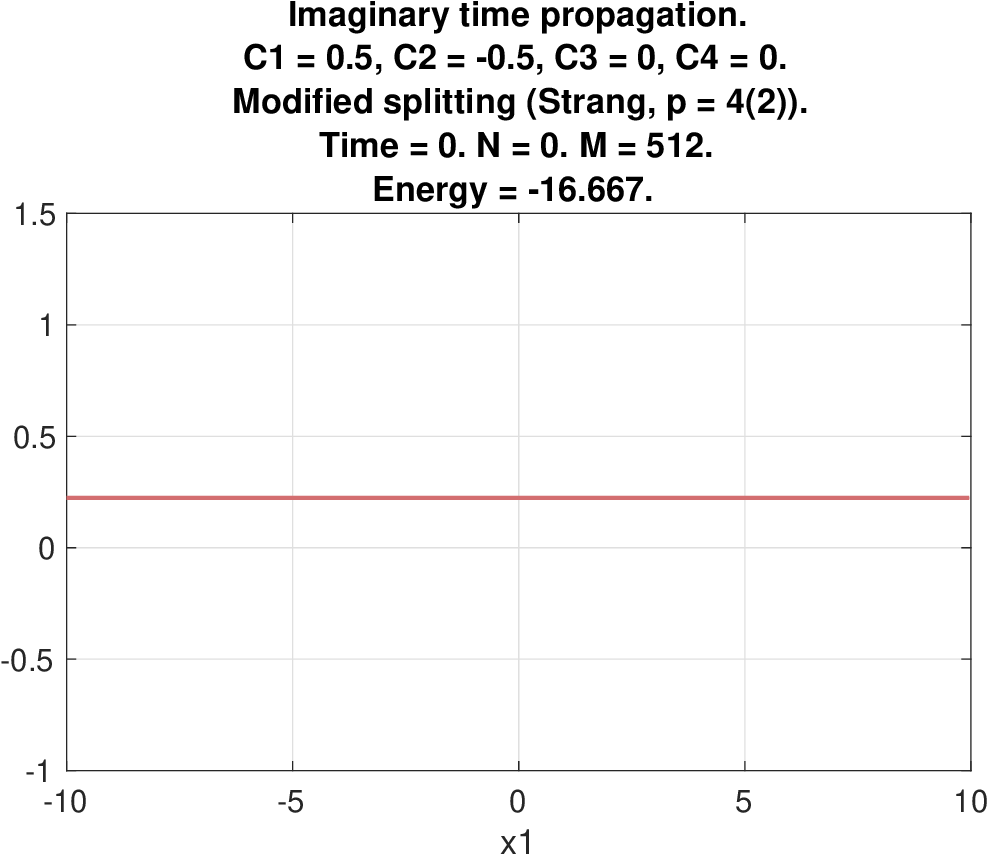} \quad
\includegraphics[width=4.2cm]{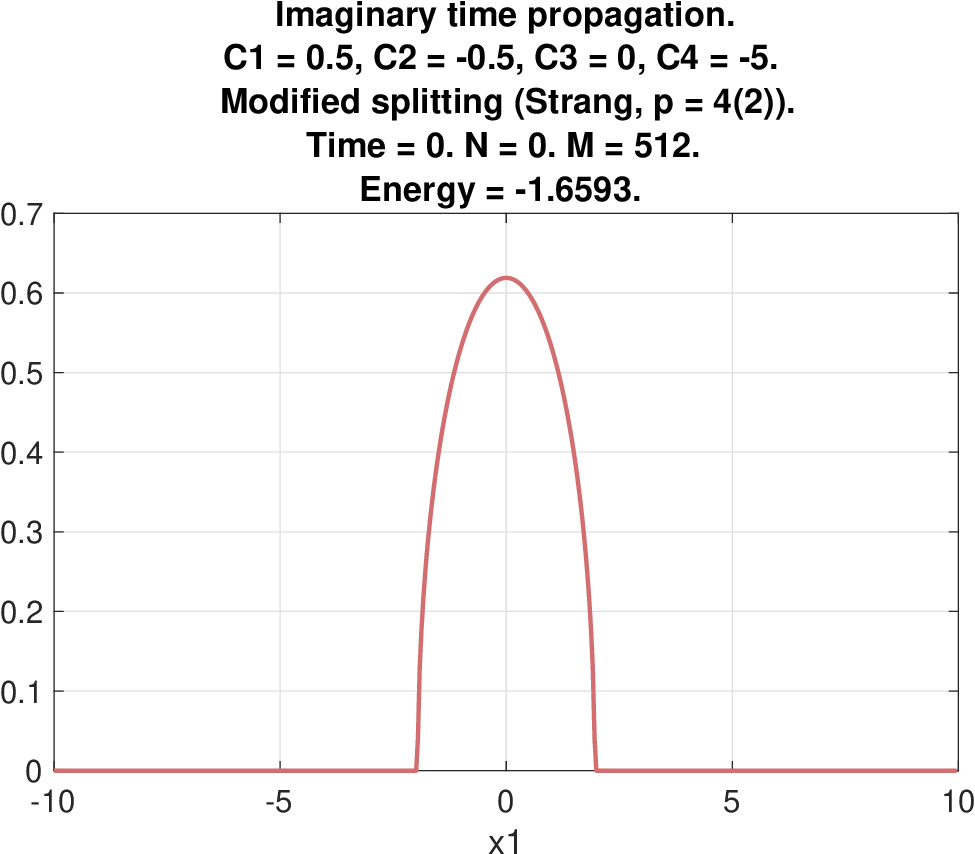} \quad
\includegraphics[width=4.2cm]{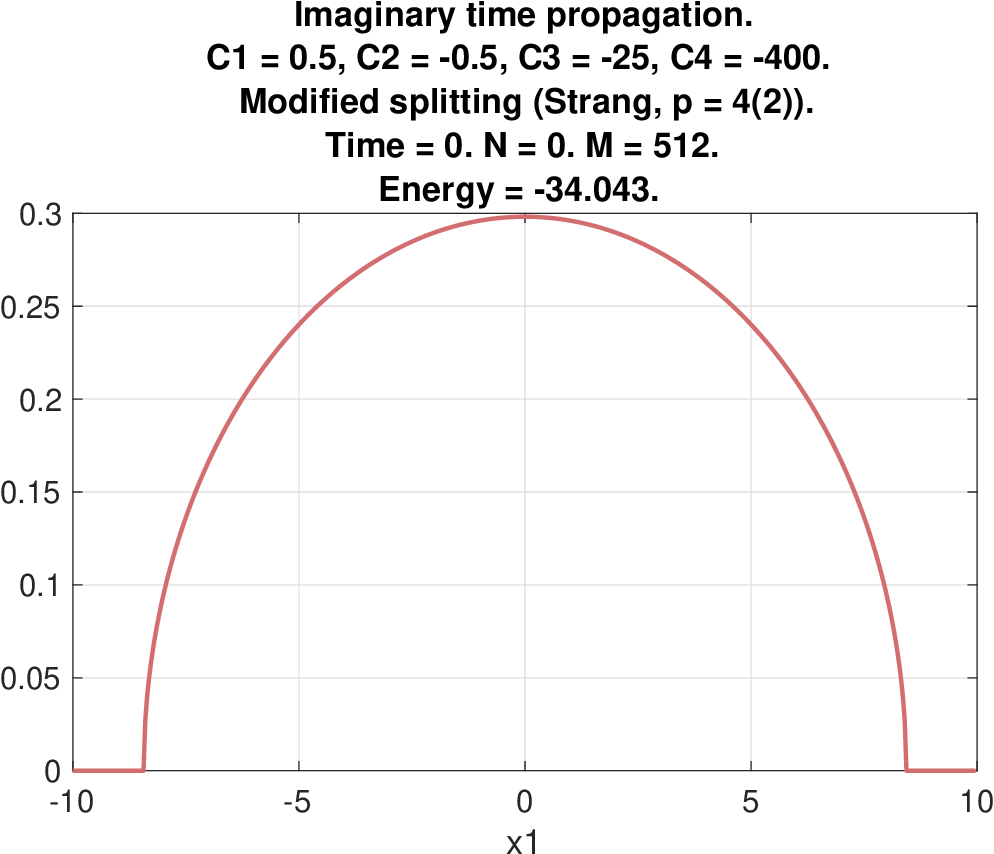} \\[2mm]
\includegraphics[width=4.2cm]{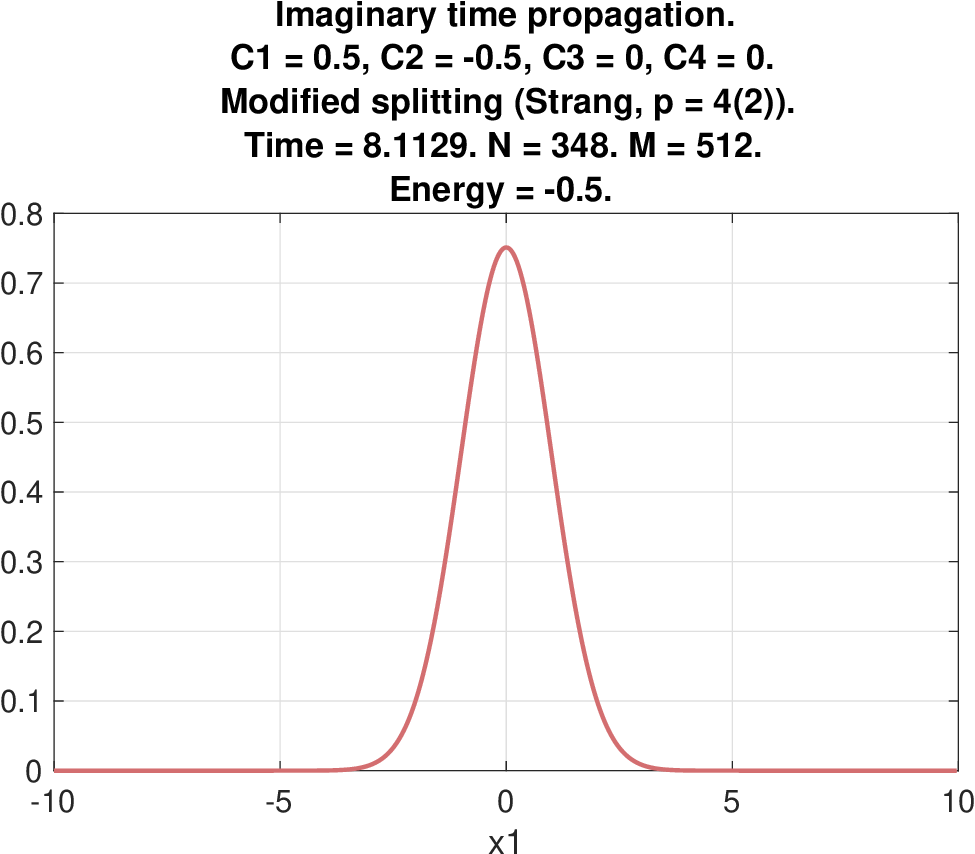} \quad
\includegraphics[width=4.2cm]{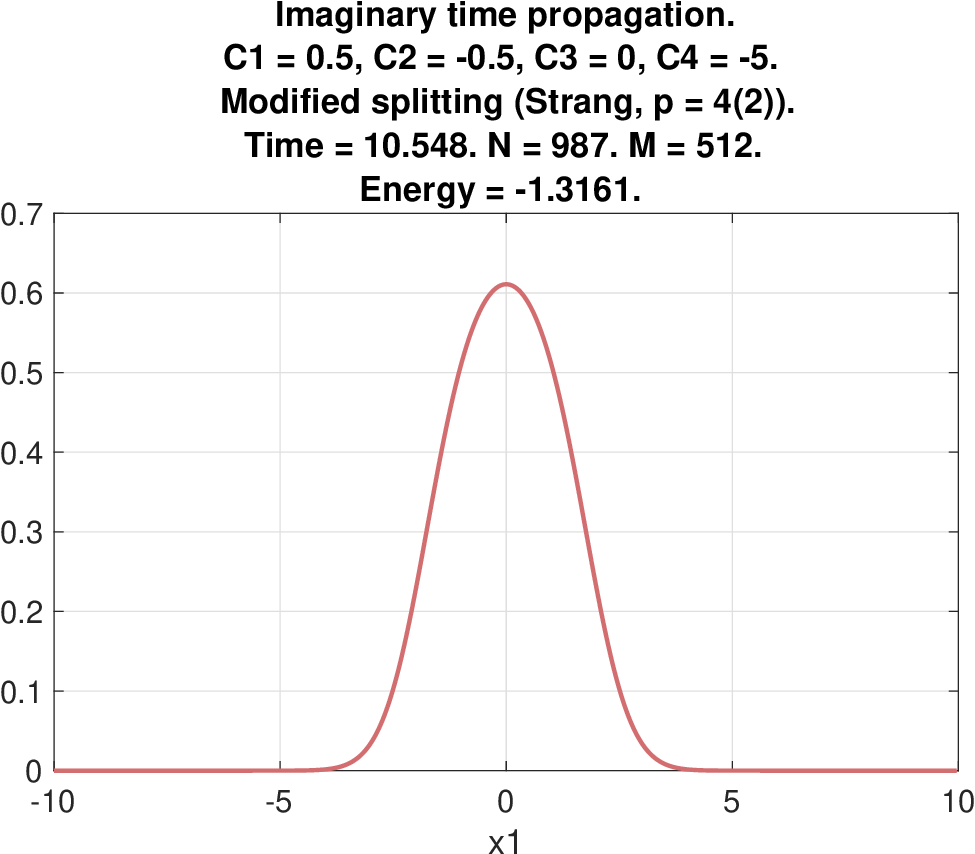} \quad
\includegraphics[width=4.2cm]{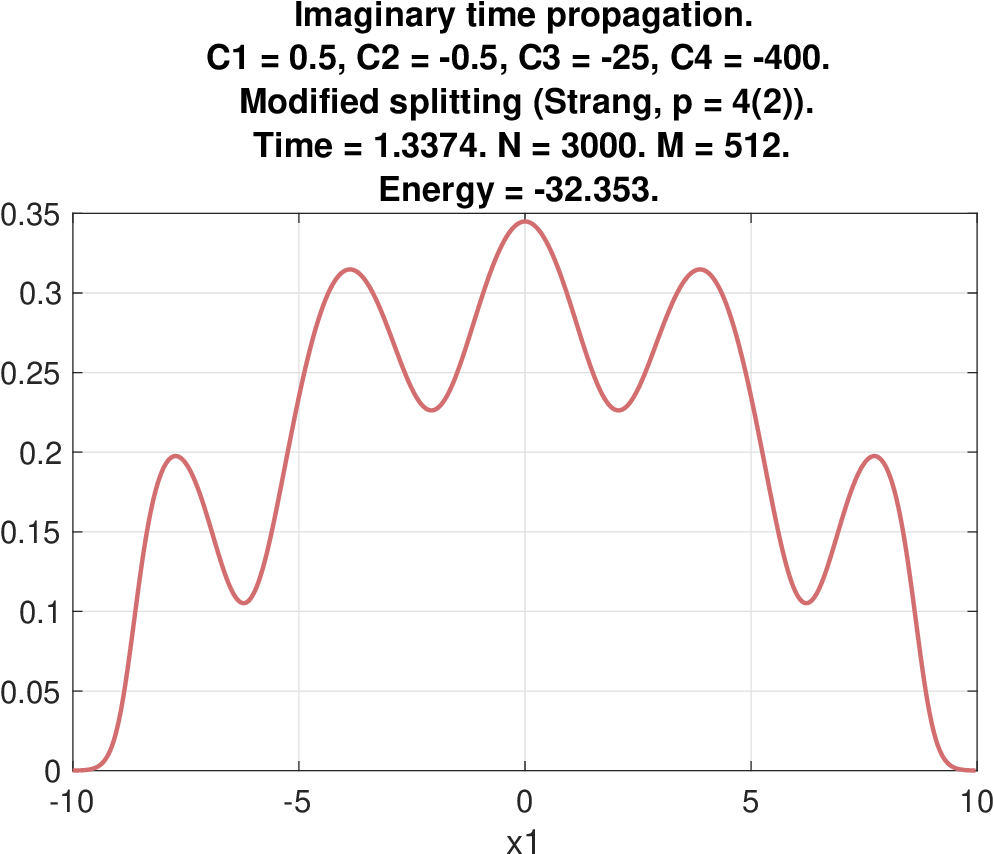} \\[2mm]
\includegraphics[width=4.2cm]{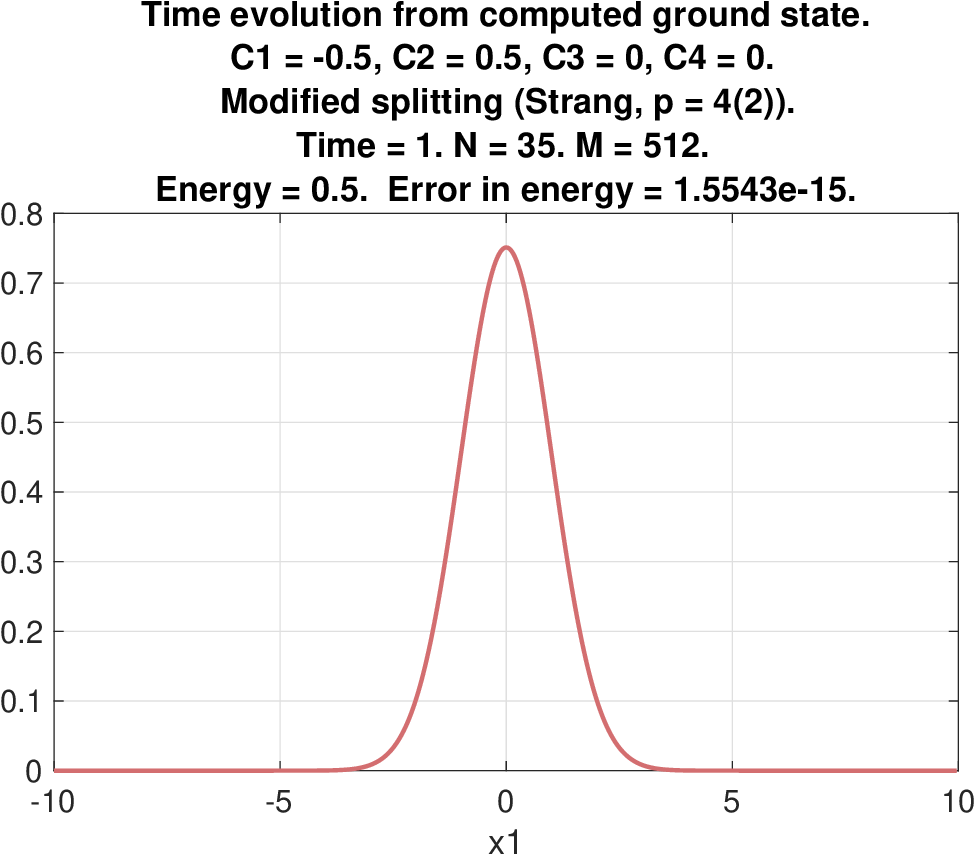} \quad
\includegraphics[width=4.2cm]{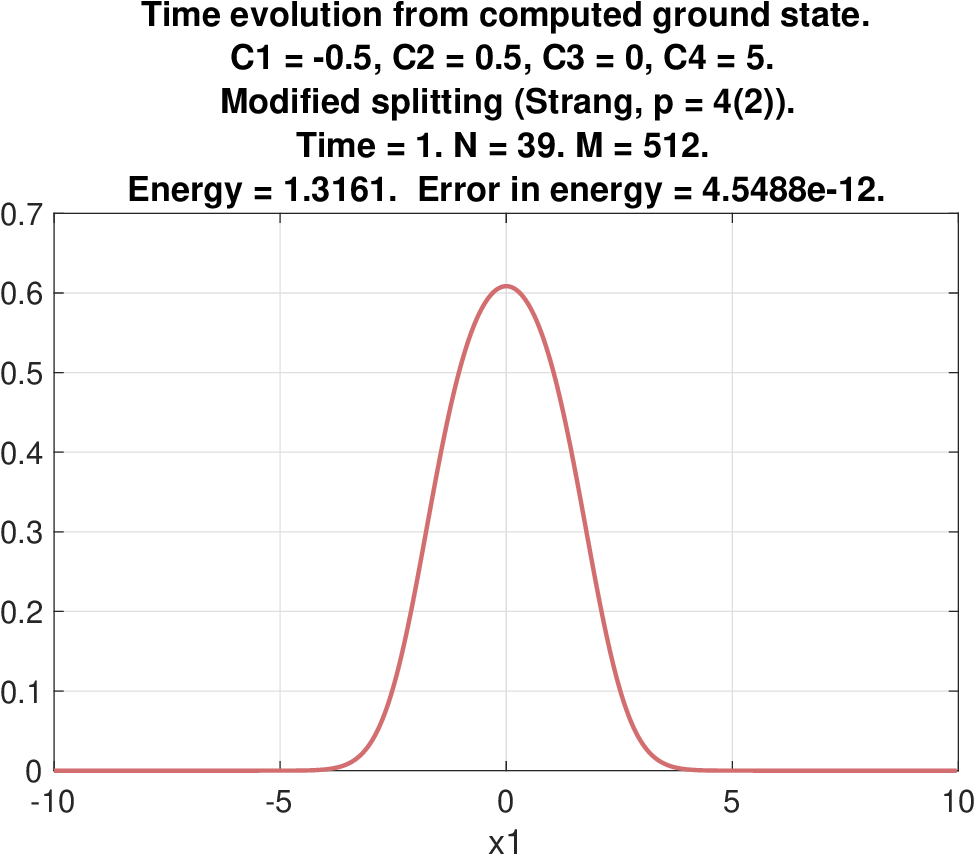} \quad
\includegraphics[width=4.2cm]{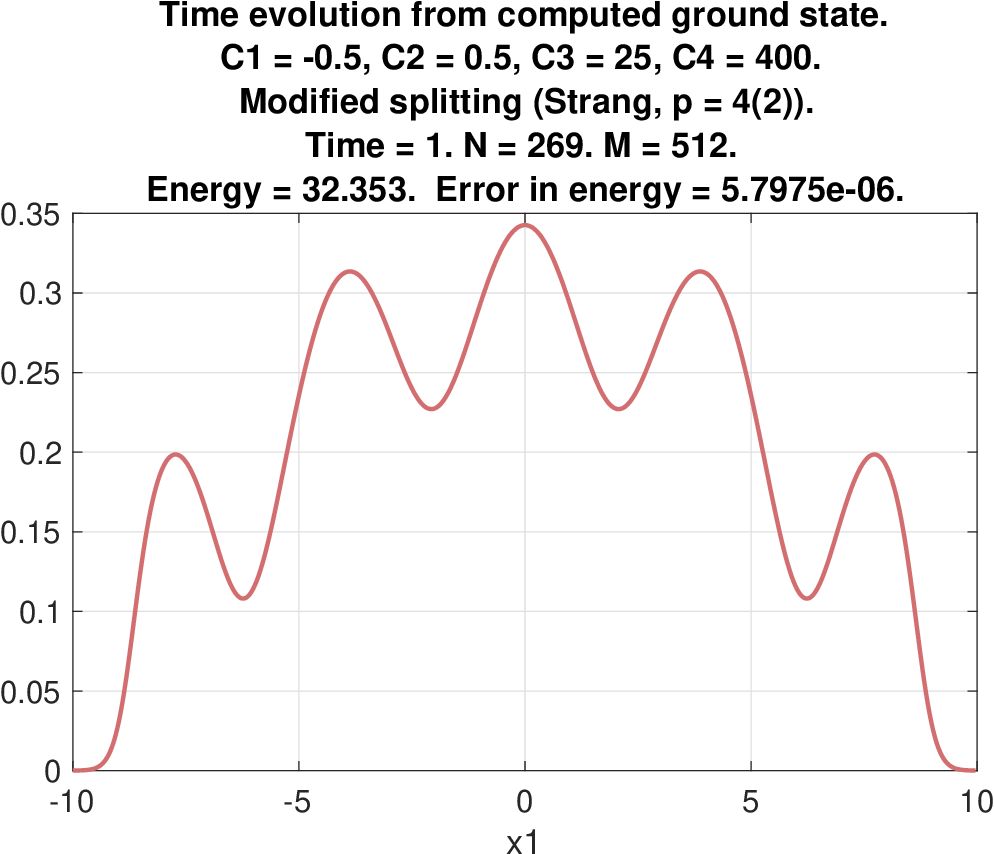} \\[2mm]
\includegraphics[width=4.2cm]{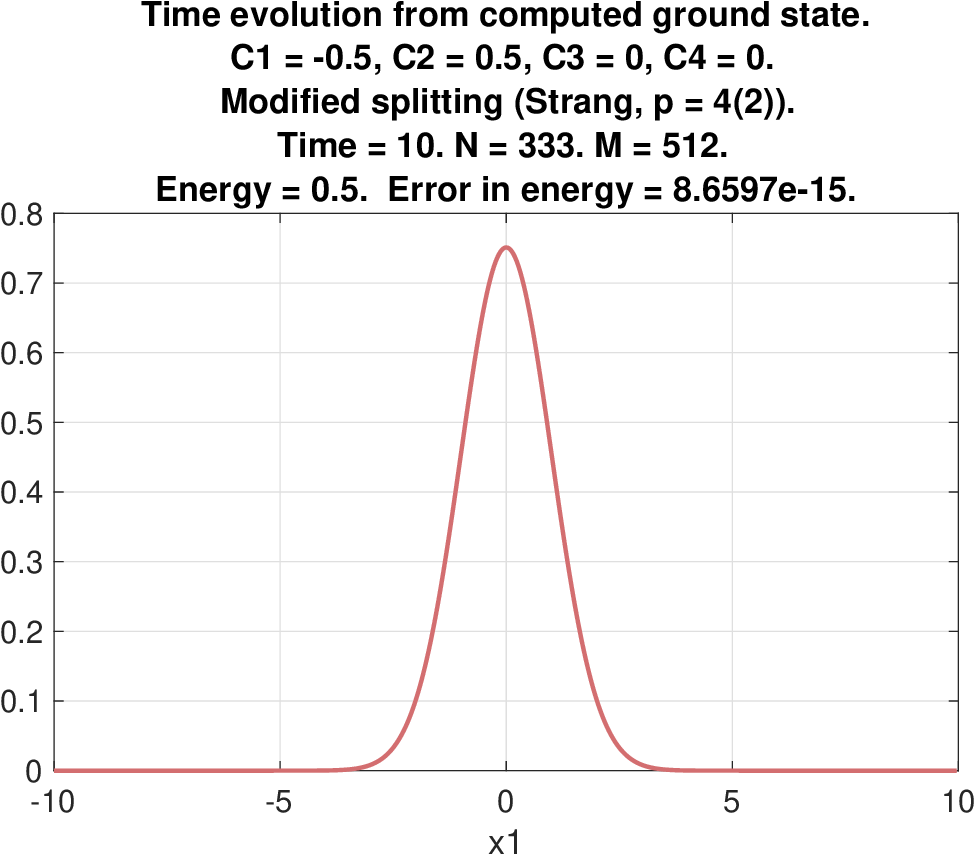} \quad
\includegraphics[width=4.2cm]{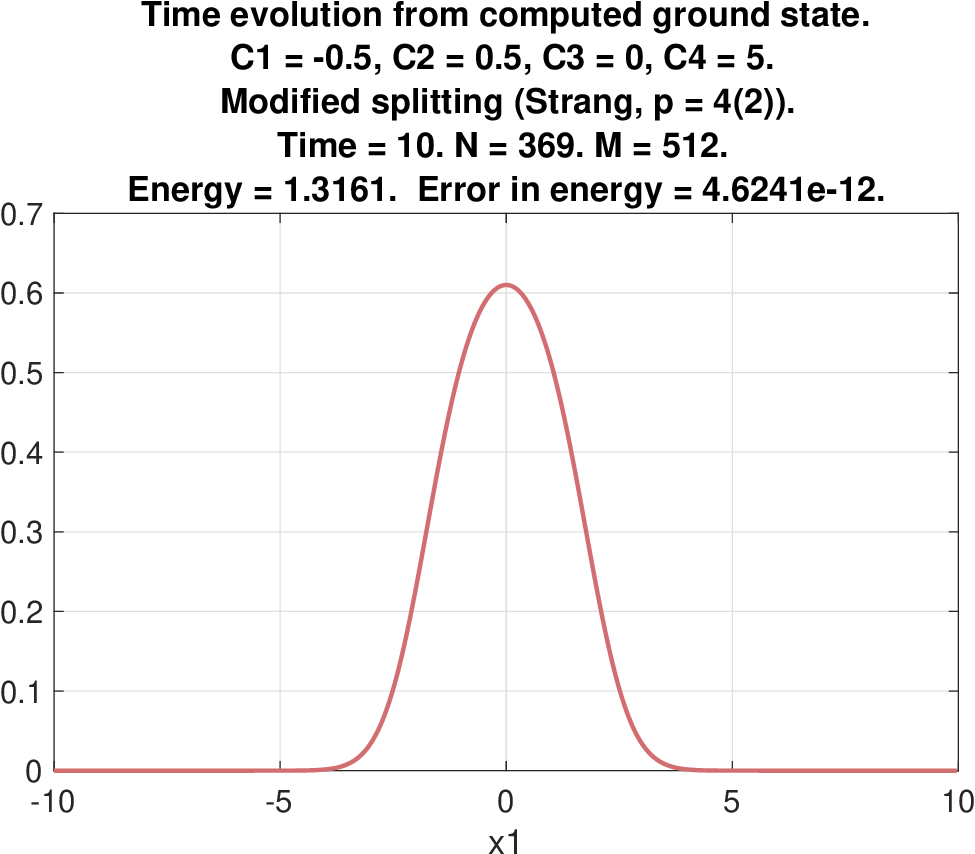} \quad
\includegraphics[width=4.2cm]{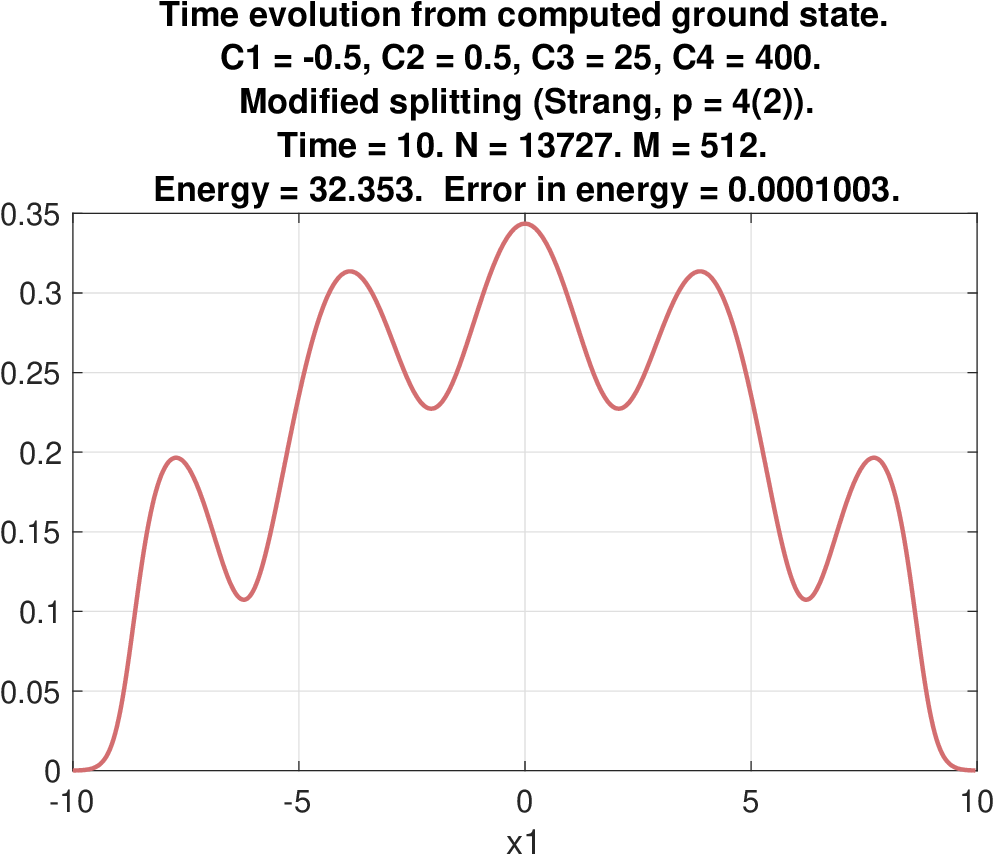} 
\end{center}
\caption{Imaginary time propagation and subsequent time evolution of a one-dimensional linear Schrödinger equation (first column), a one-dimensional Gross--Pitaevskii equation with moderate nonlinearity (second column), and a one-dimensional Gross--Pitaevskii equation with additional multiple lattice potential and large constant in the nonlinearity (third column). 
Profiles of the prescribed initial states for the imaginary time propagation (first row).
The time-dependent solutions at a short intermediate time (third row) and a larger final time (fourth row) coincide with the resulting ground state solutions at the imaginary time~$T$ (second row) in modulus. 
The numbers of iterations, the values of the total energies, and the corresponding errors are included in the headlines of the figures.}
\label{fig:FigureITTE1}
\end{figure}

\begin{figure}
\begin{center}
\includegraphics[width=4.2cm]{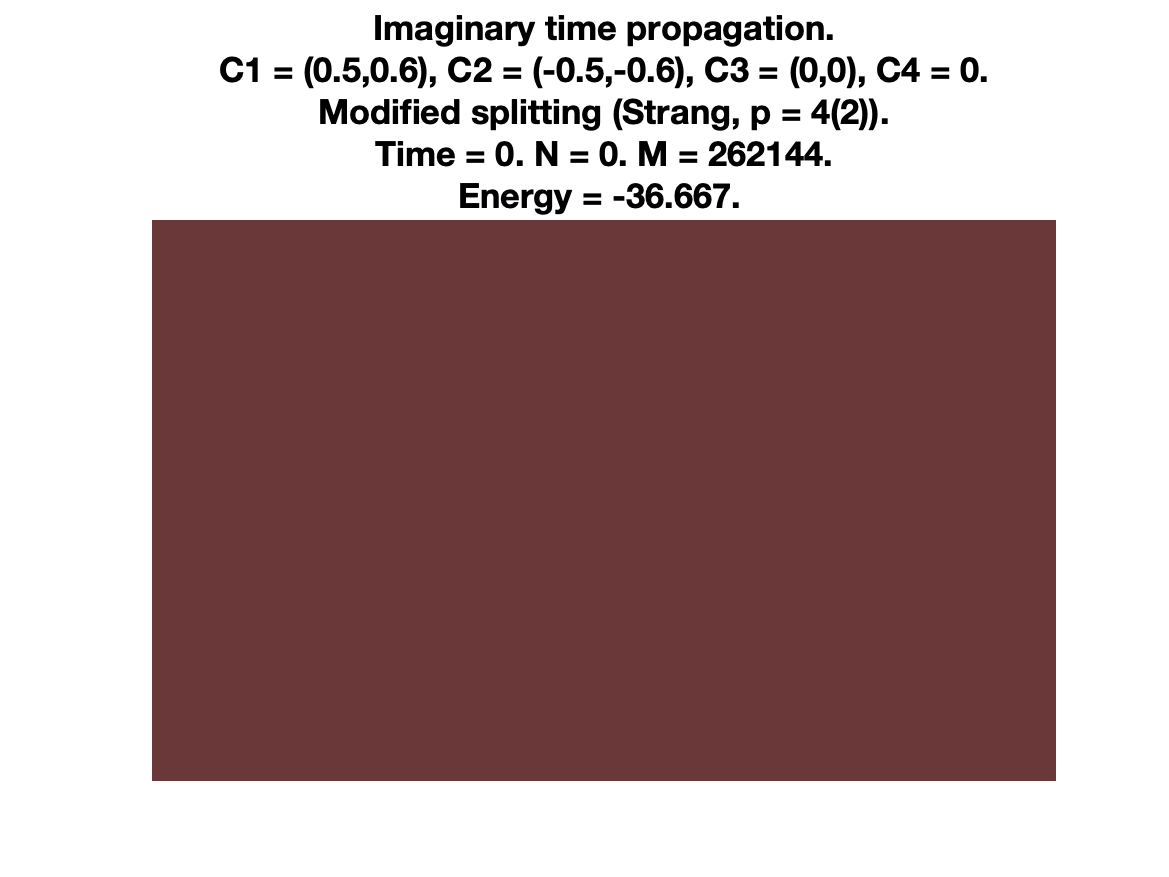} \quad
\includegraphics[width=4.2cm]{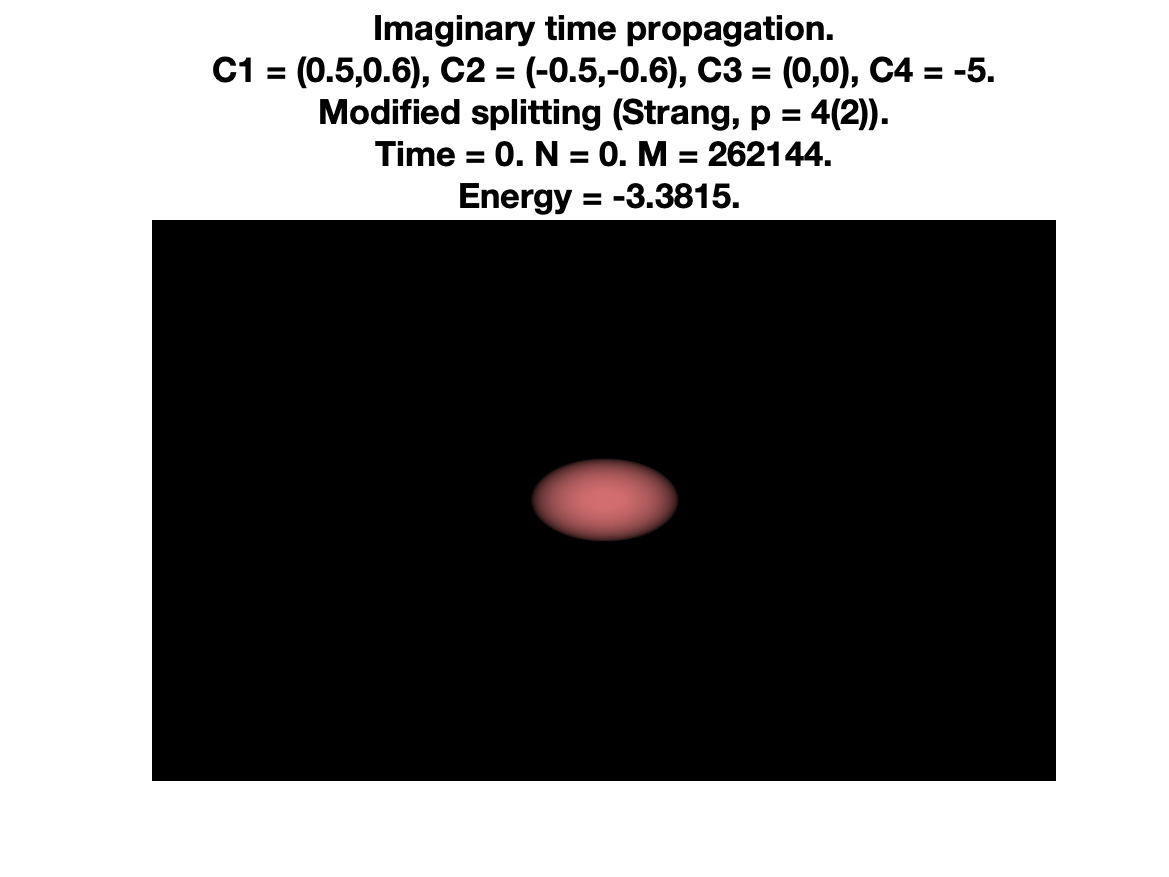} \quad
\includegraphics[width=4.2cm]{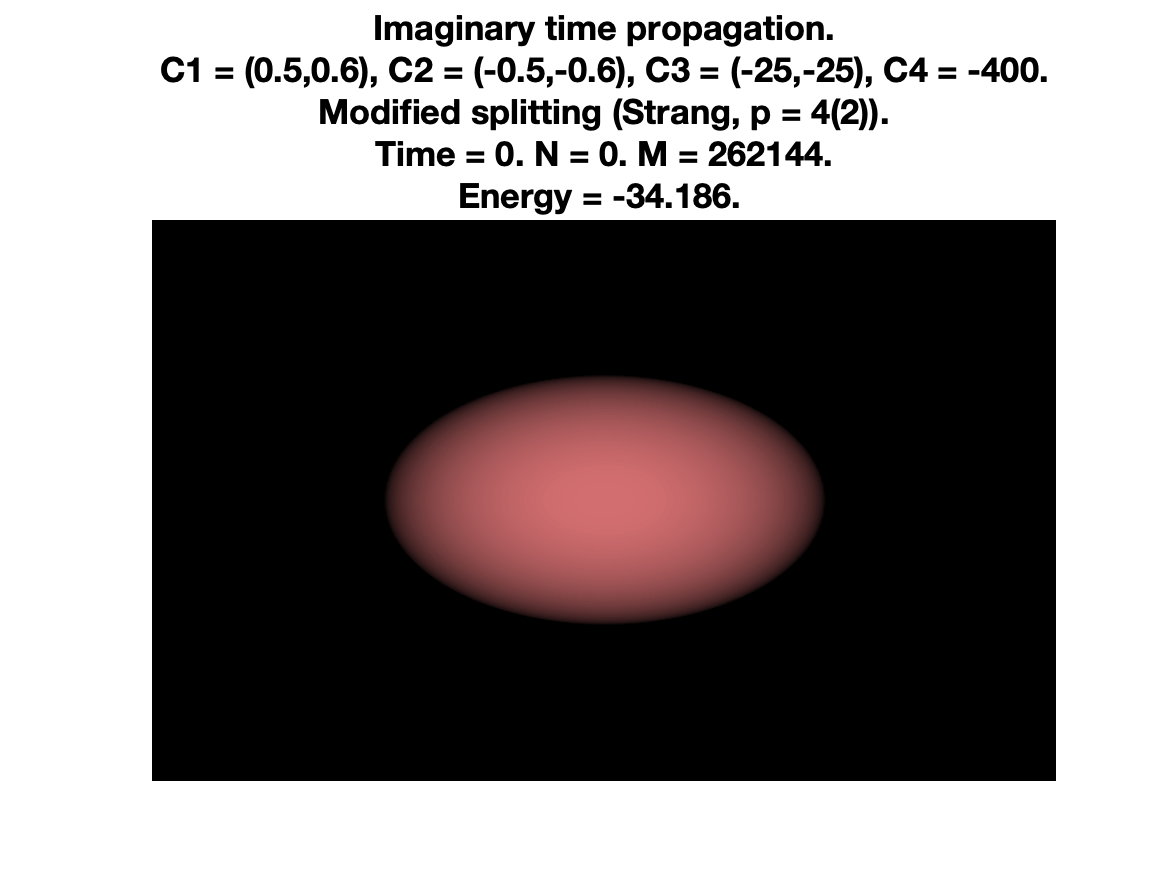} \\[2mm]
\includegraphics[width=4.2cm]{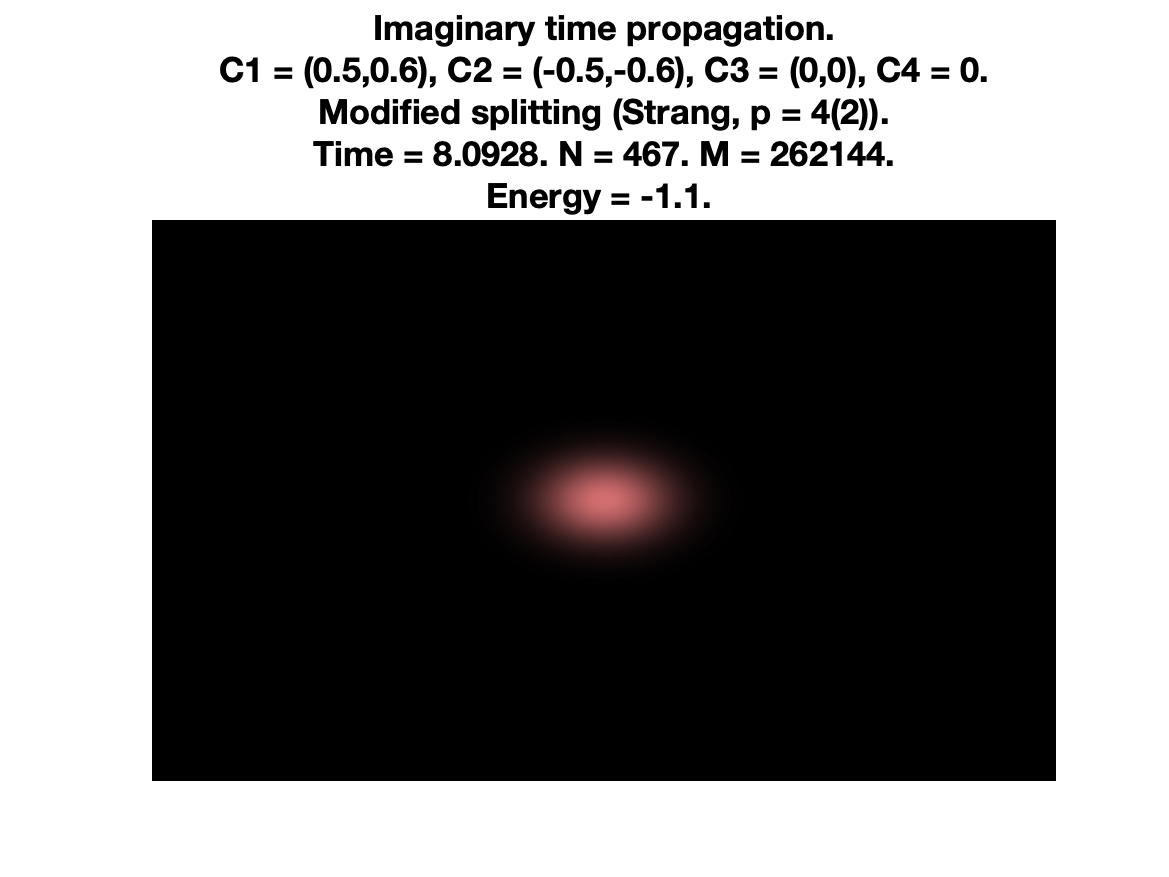} \quad
\includegraphics[width=4.2cm]{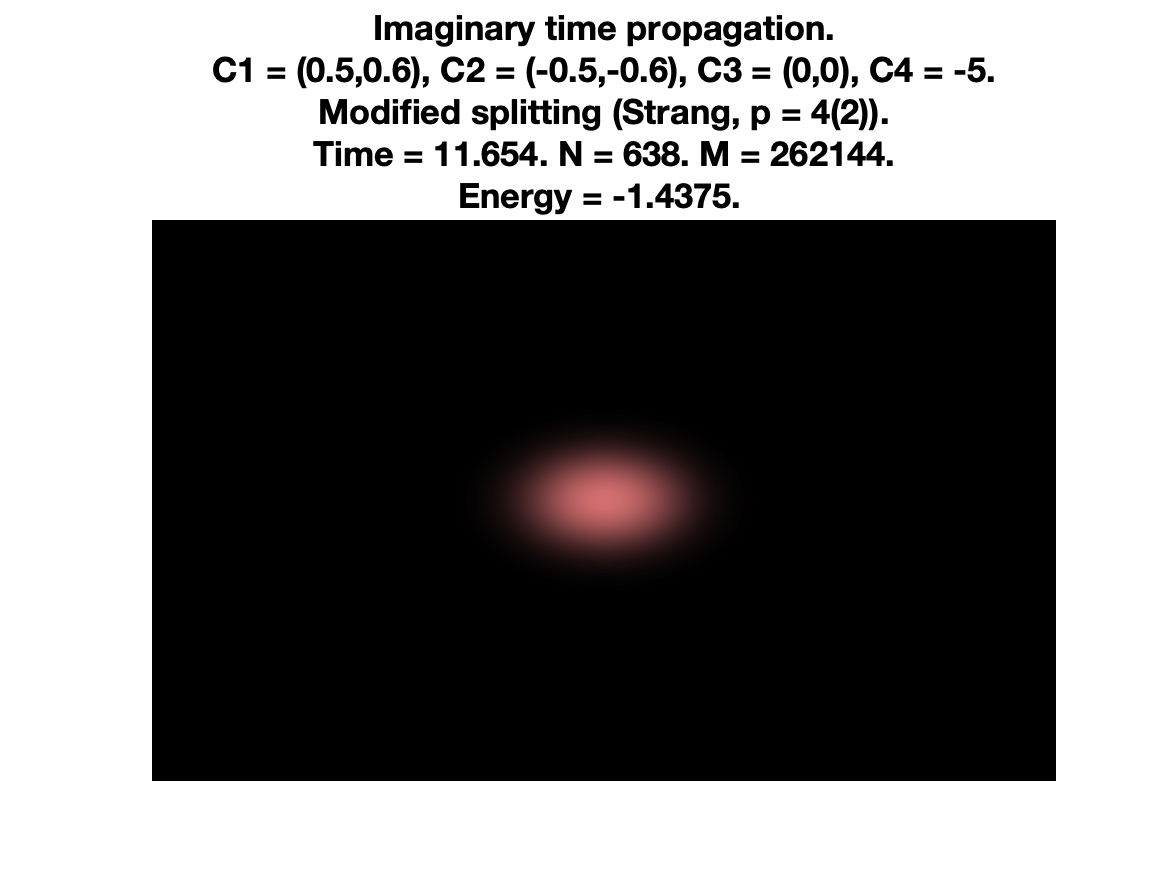} \quad
\includegraphics[width=4.2cm]{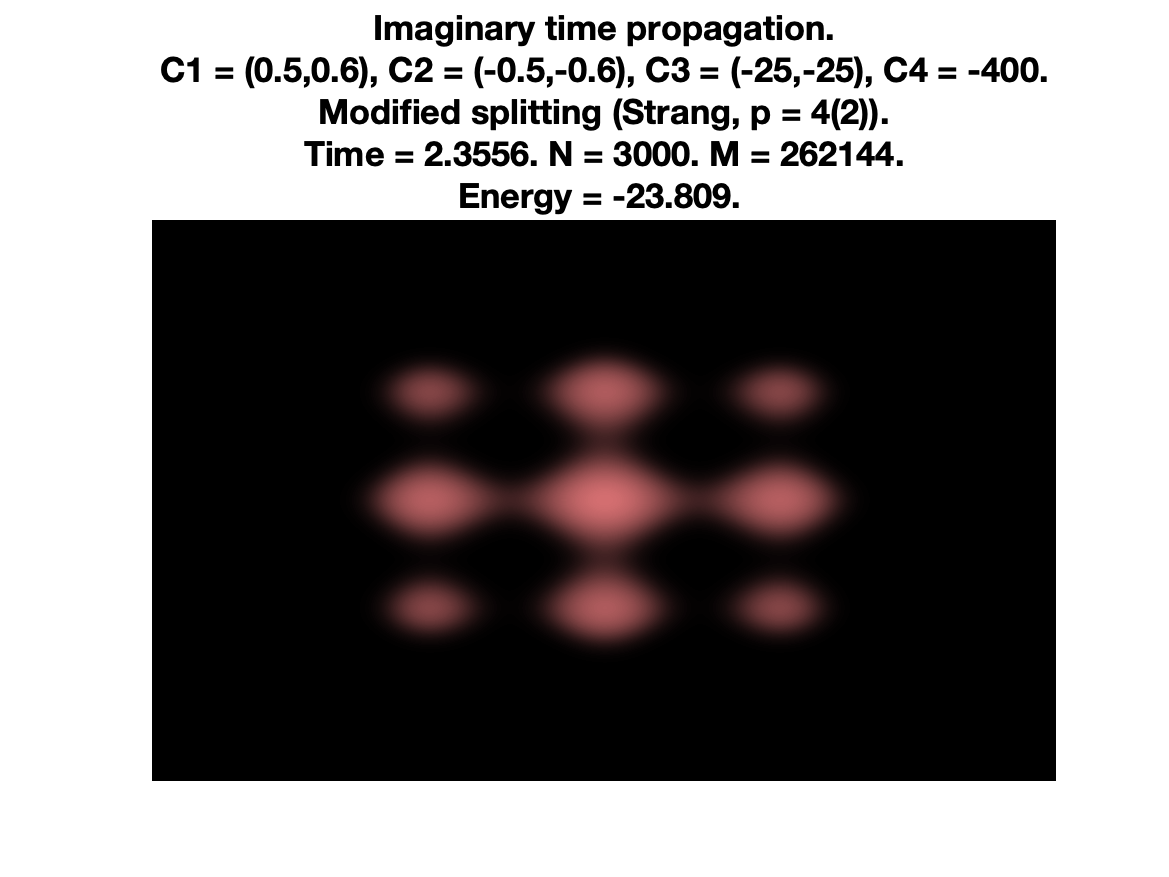} \\[2mm]
\includegraphics[width=4.2cm]{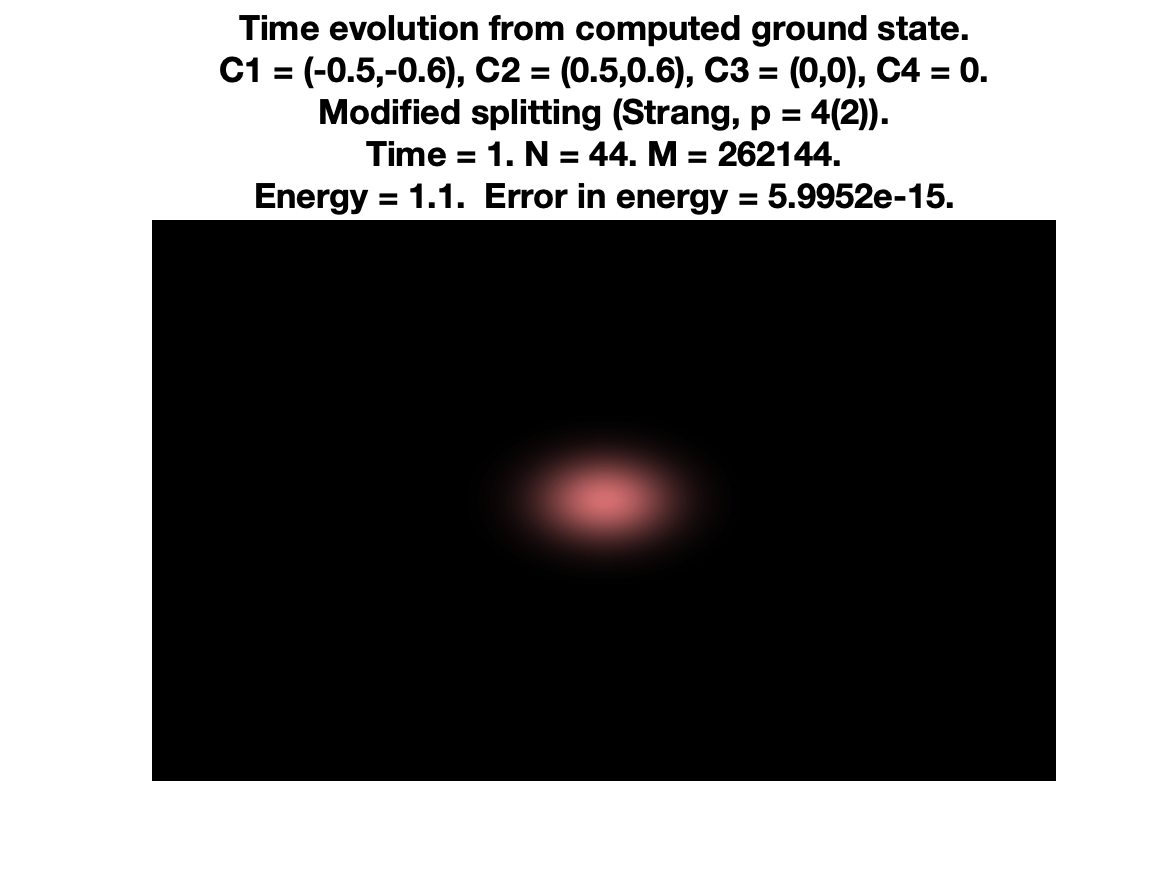} \quad
\includegraphics[width=4.2cm]{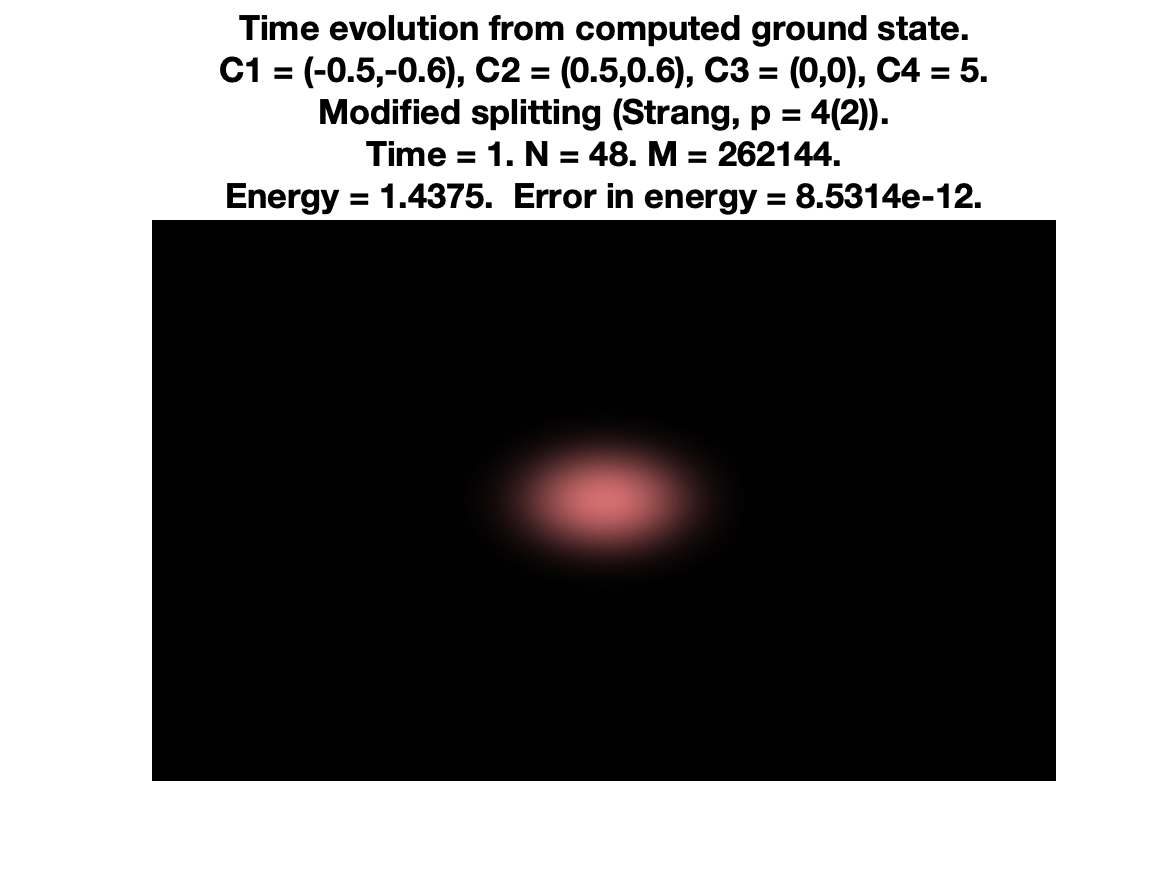} \quad
\includegraphics[width=4.2cm]{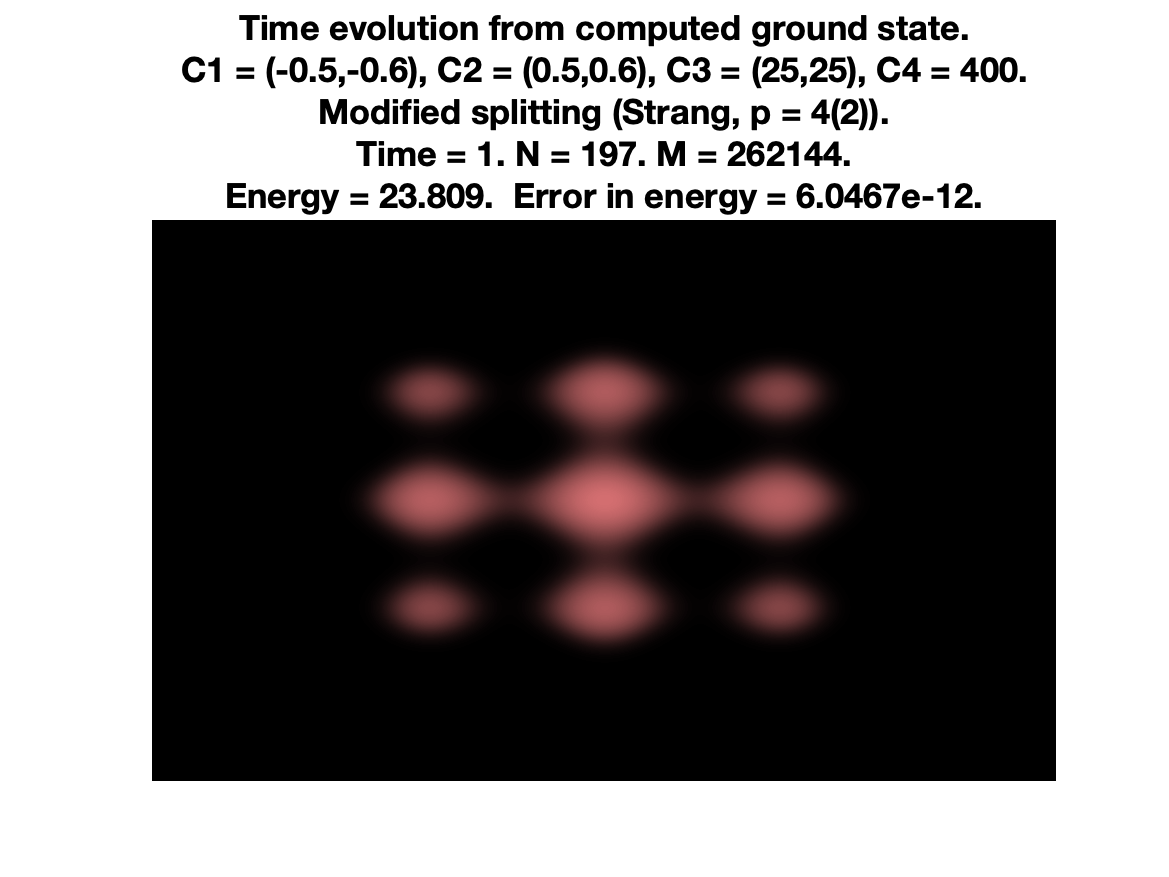} \\[2mm]
\includegraphics[width=4.2cm]{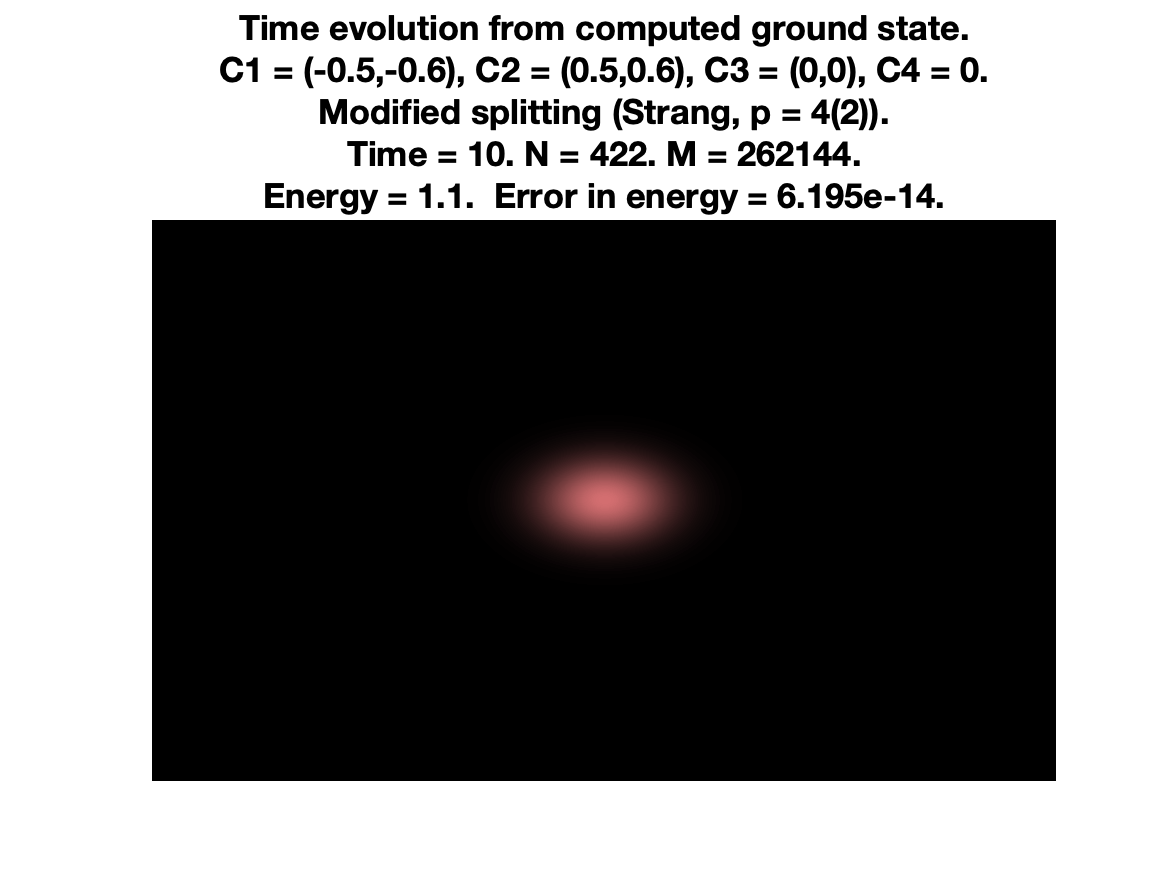} \quad
\includegraphics[width=4.2cm]{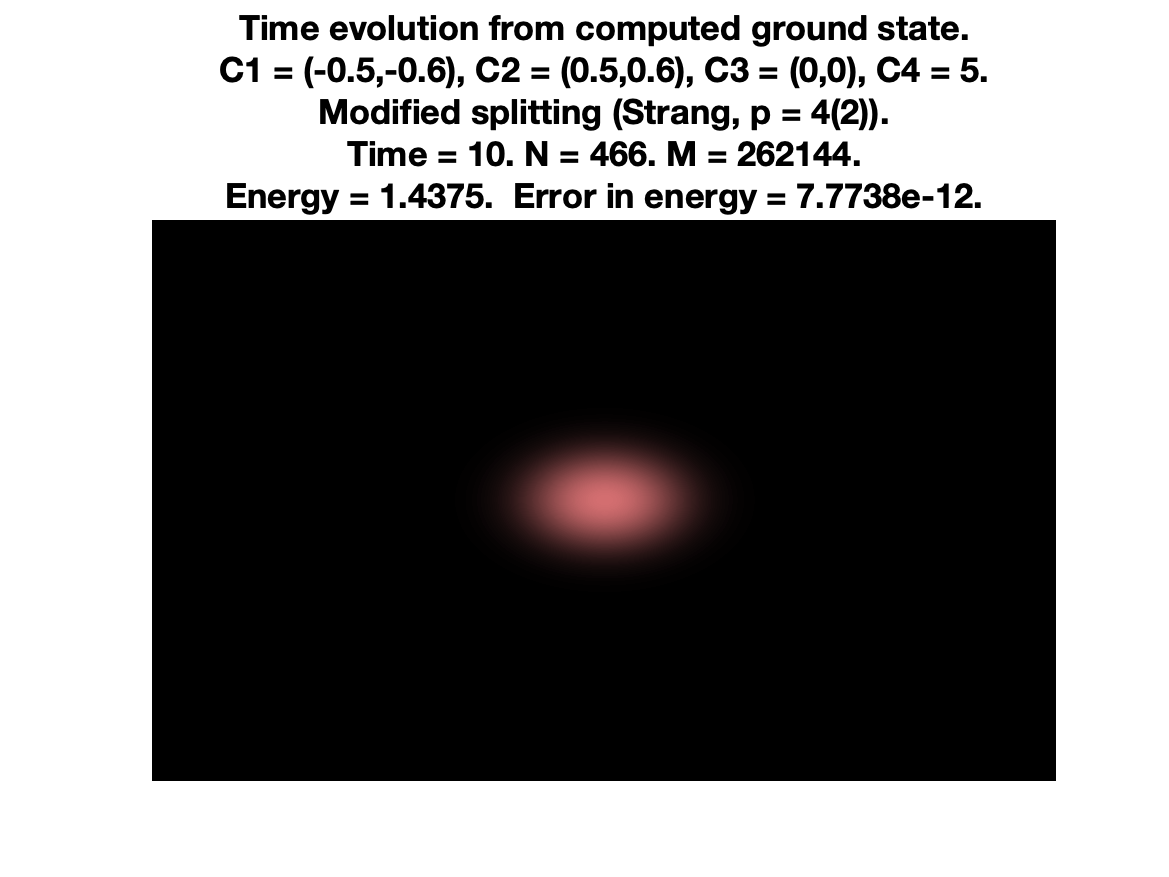} \quad
\includegraphics[width=4.2cm]{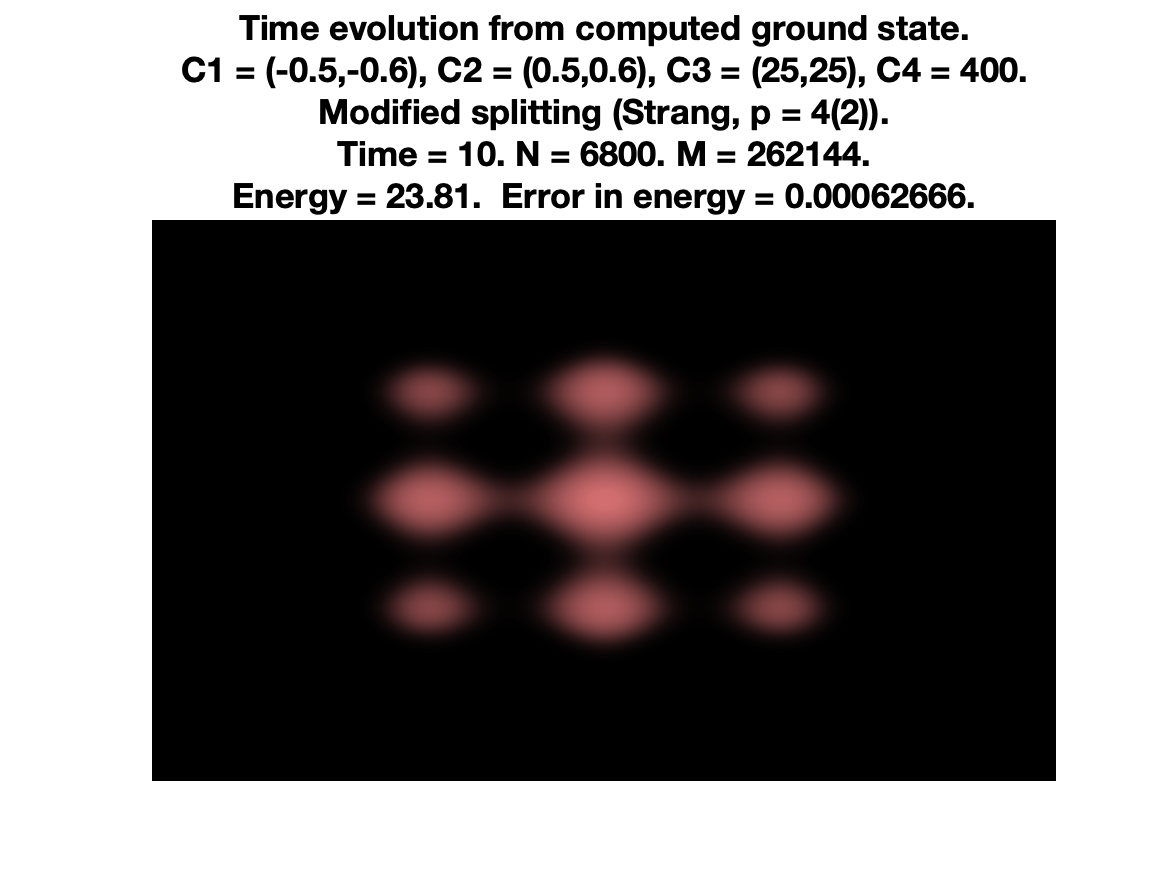} 
\end{center}
\caption{Corresponding results for two space dimensions.}
\label{fig:FigureITTE2}
\end{figure}

\begin{figure}
\begin{center}
\includegraphics[width=4.2cm]{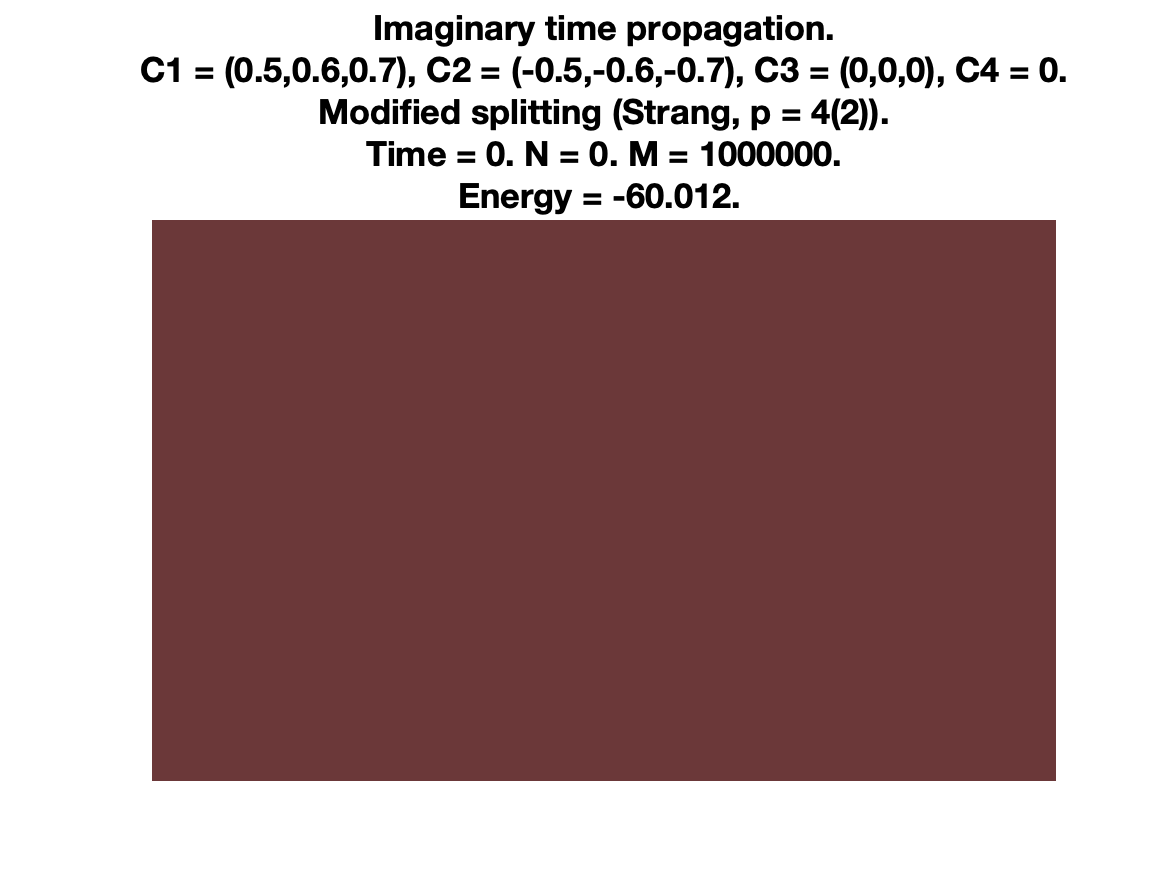} \quad
\includegraphics[width=4.2cm]{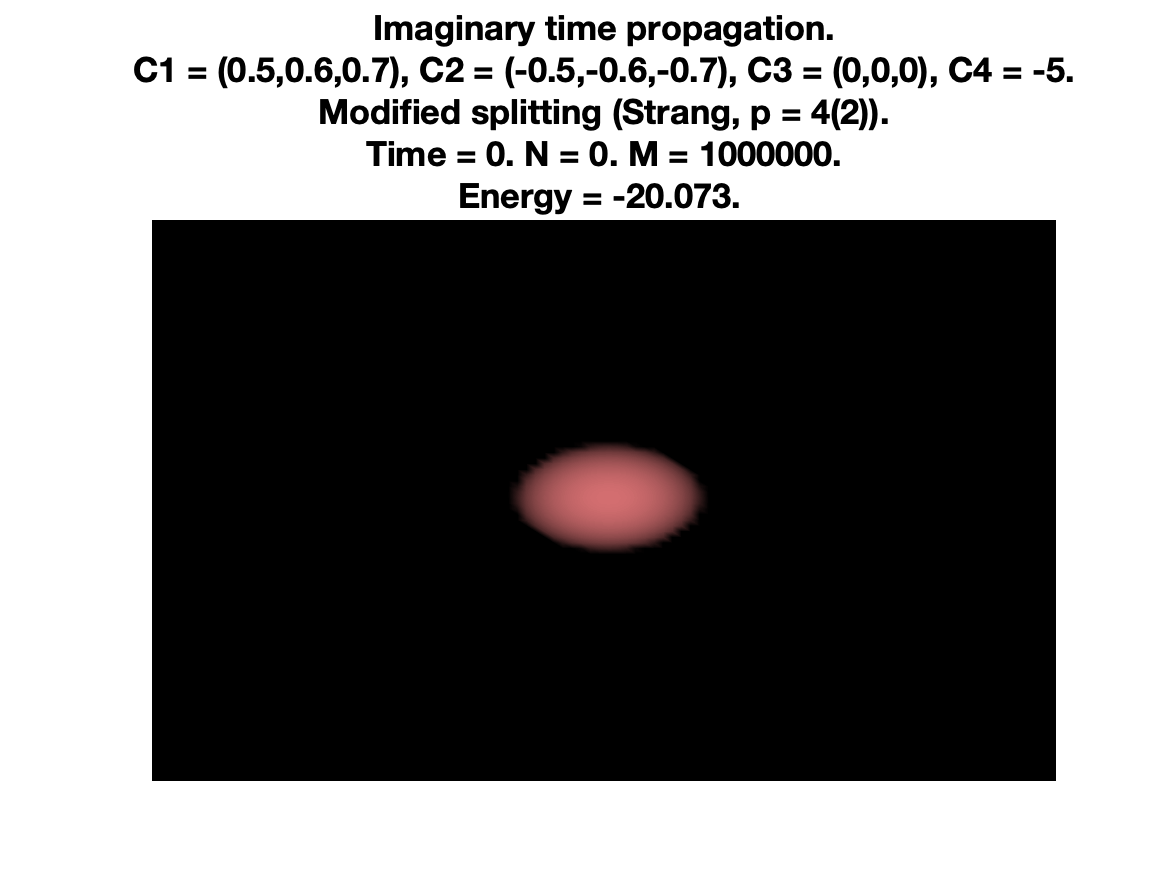} \quad
\includegraphics[width=4.2cm]{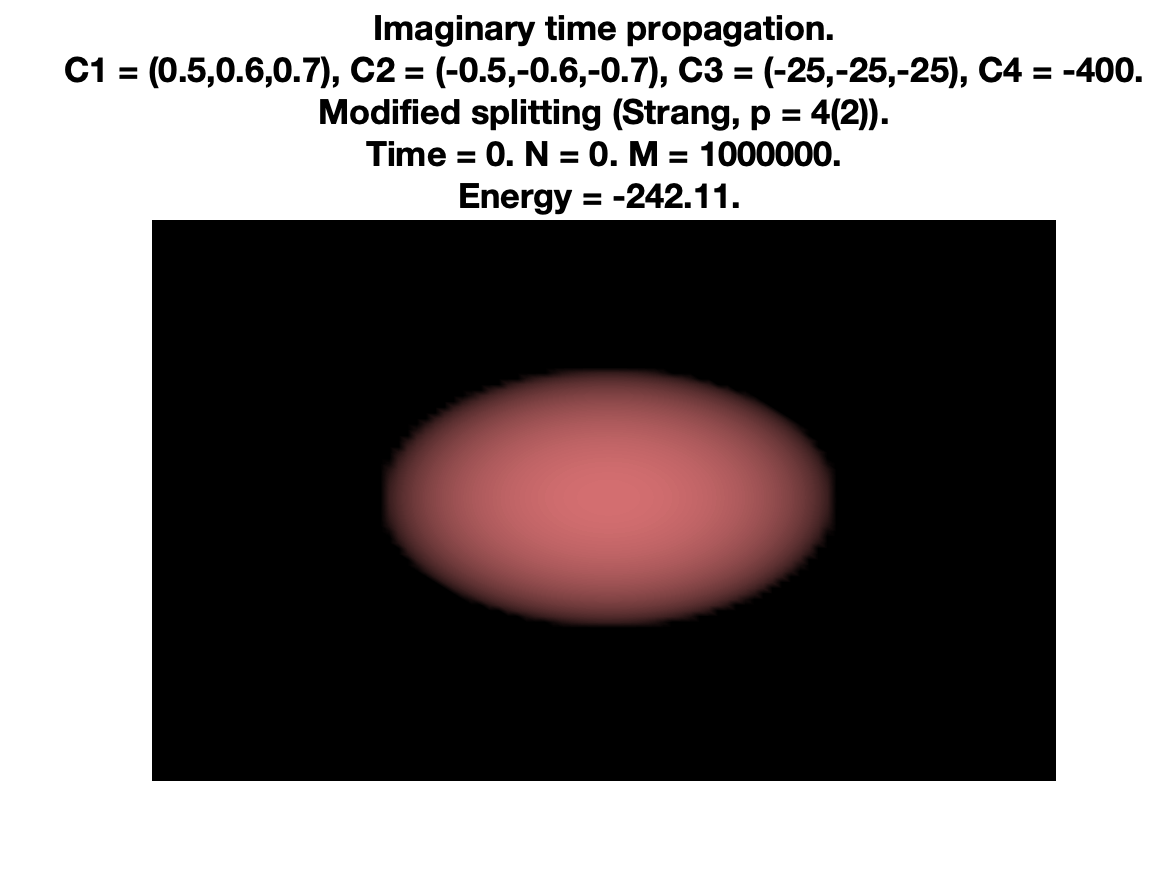} \\[2mm]
\includegraphics[width=4.2cm]{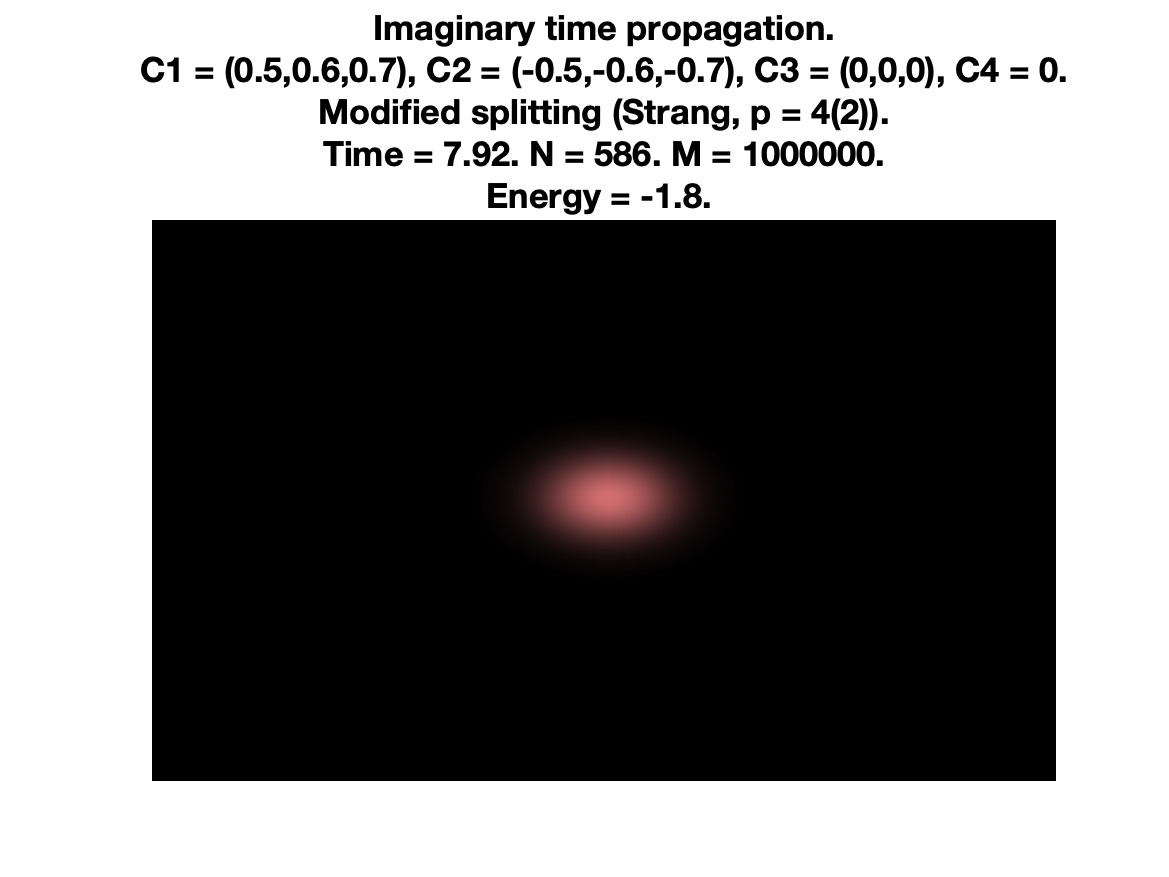} \quad
\includegraphics[width=4.2cm]{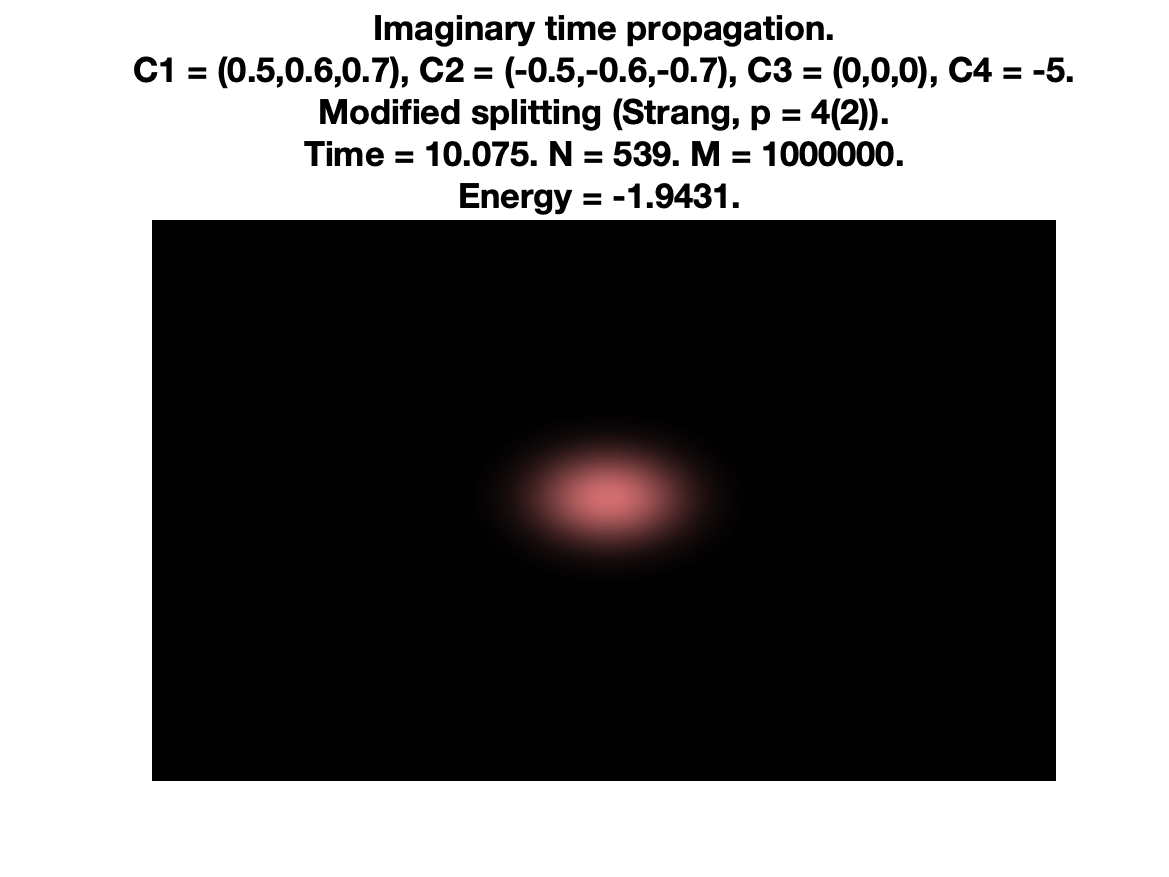} \quad
\includegraphics[width=4.2cm]{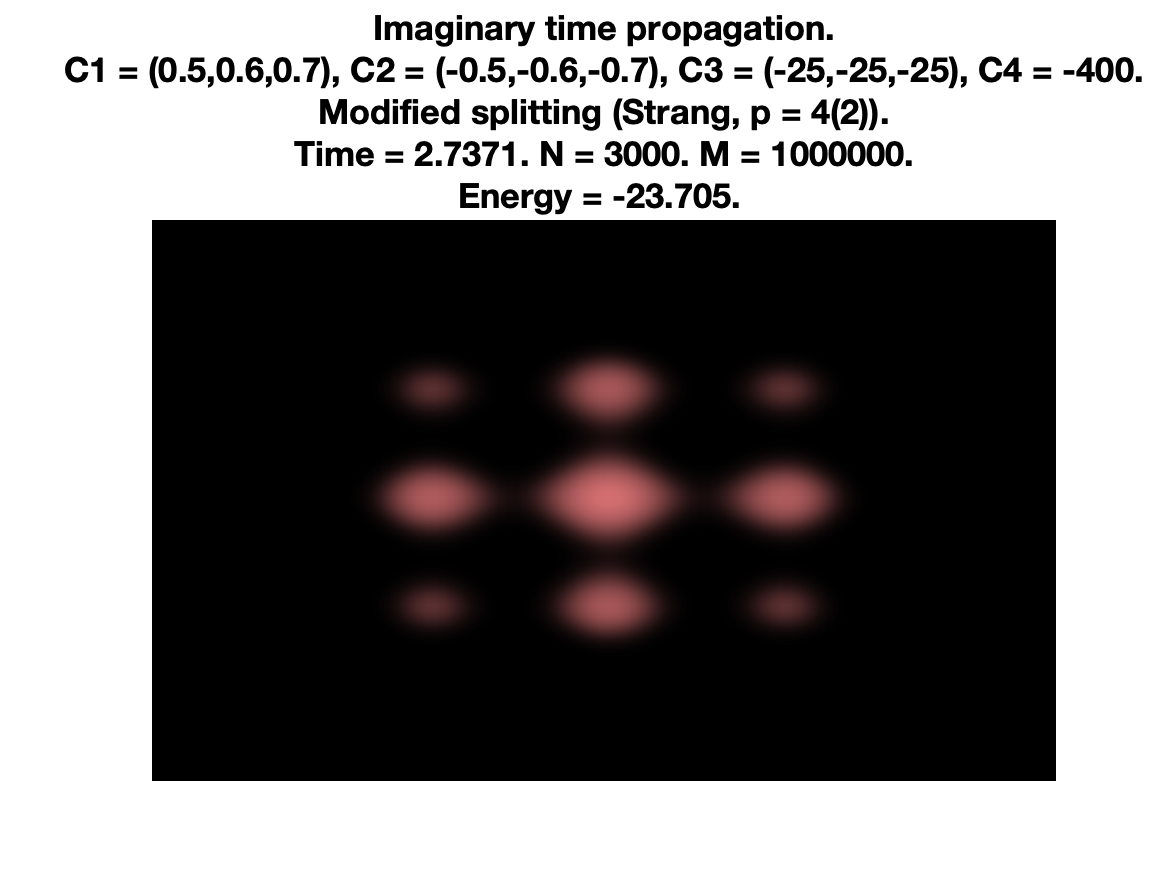} \\[2mm]
\includegraphics[width=4.2cm]{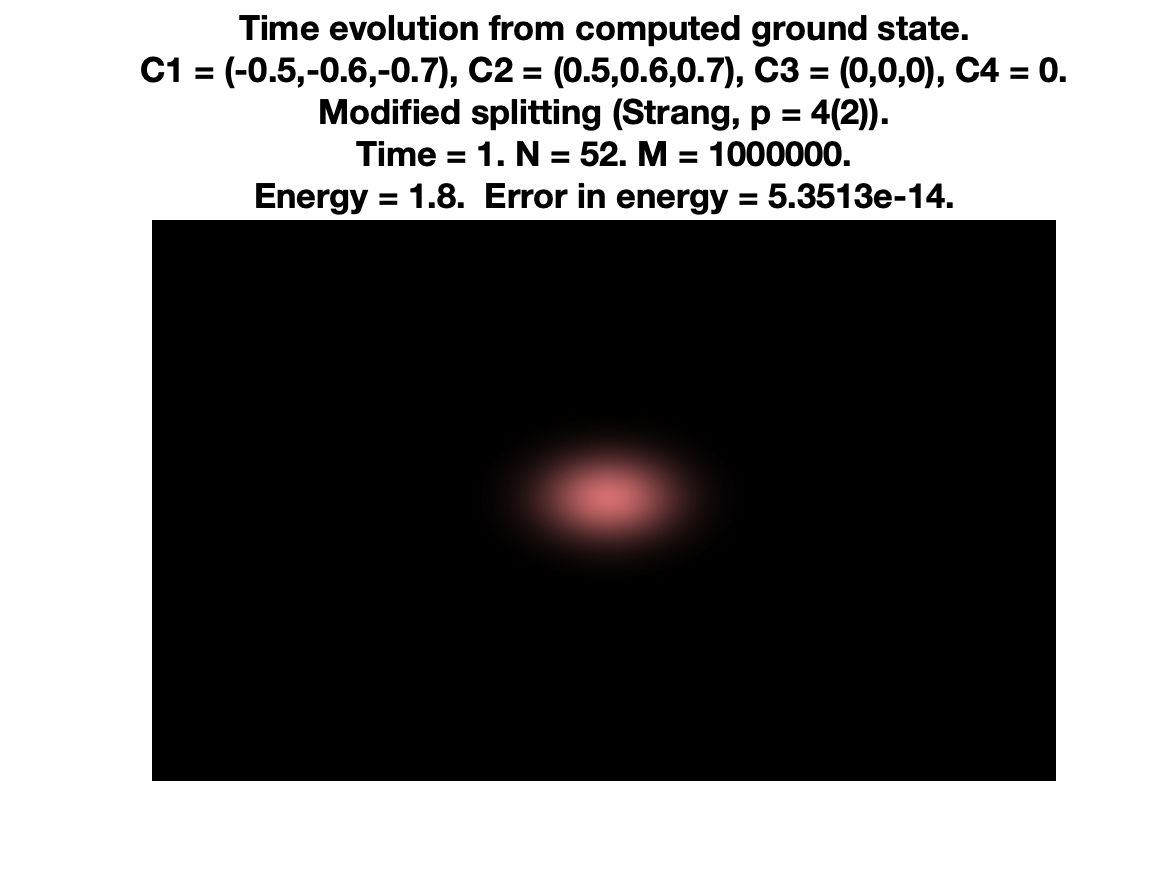} \quad
\includegraphics[width=4.2cm]{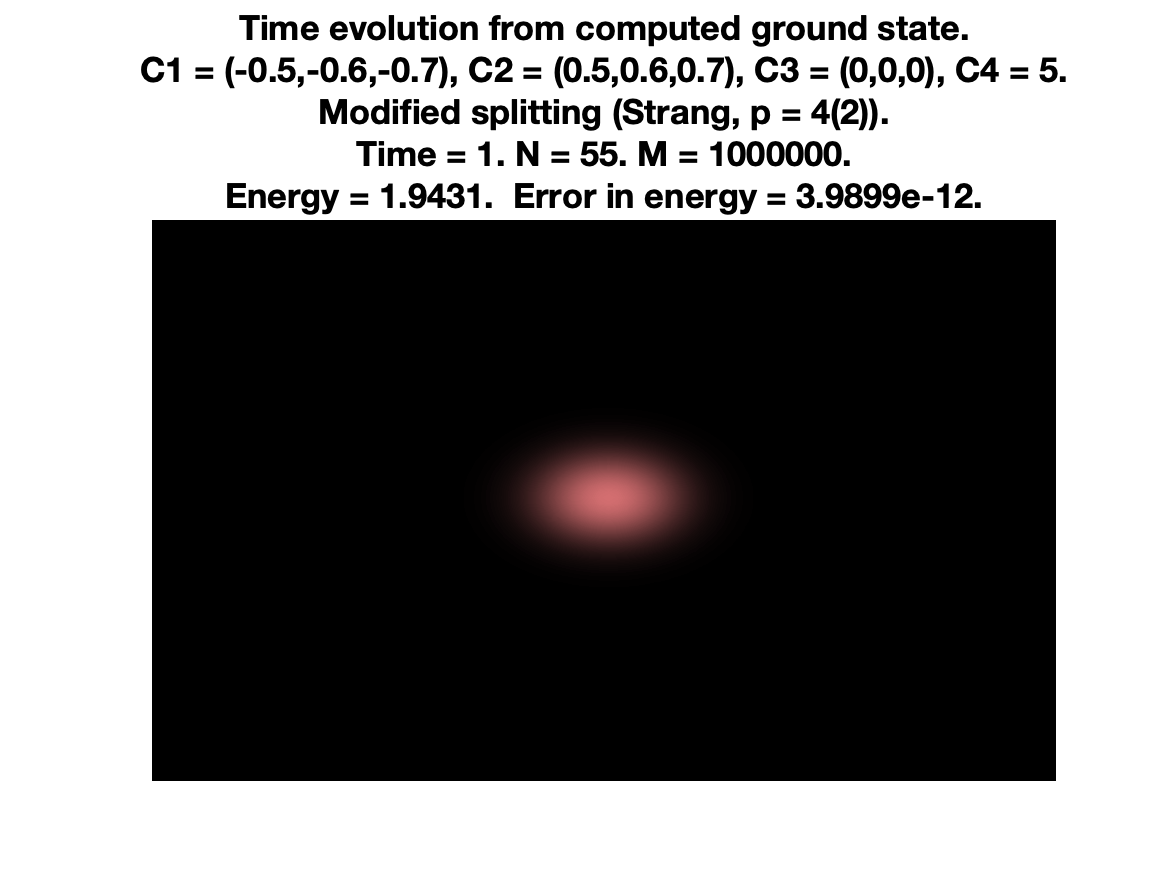} \quad
\includegraphics[width=4.2cm]{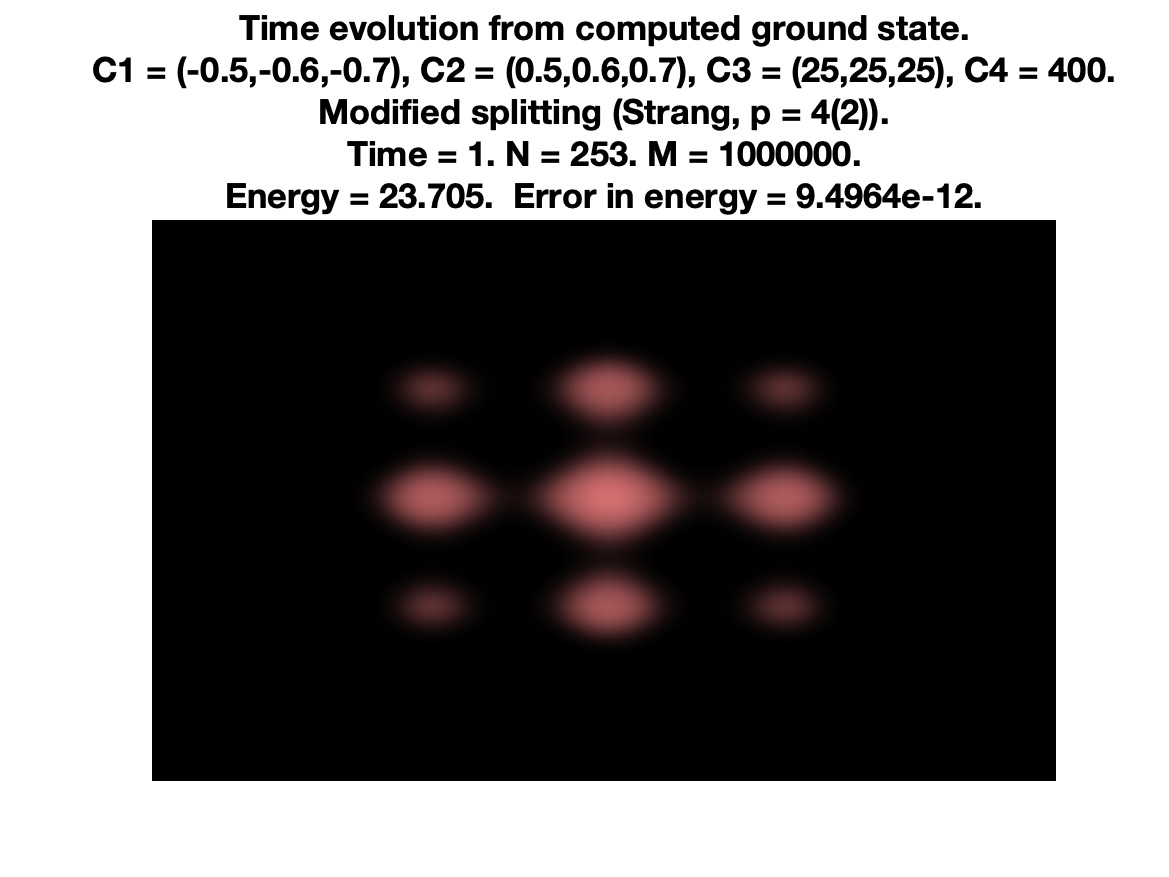} \\[2mm]
\includegraphics[width=4.2cm]{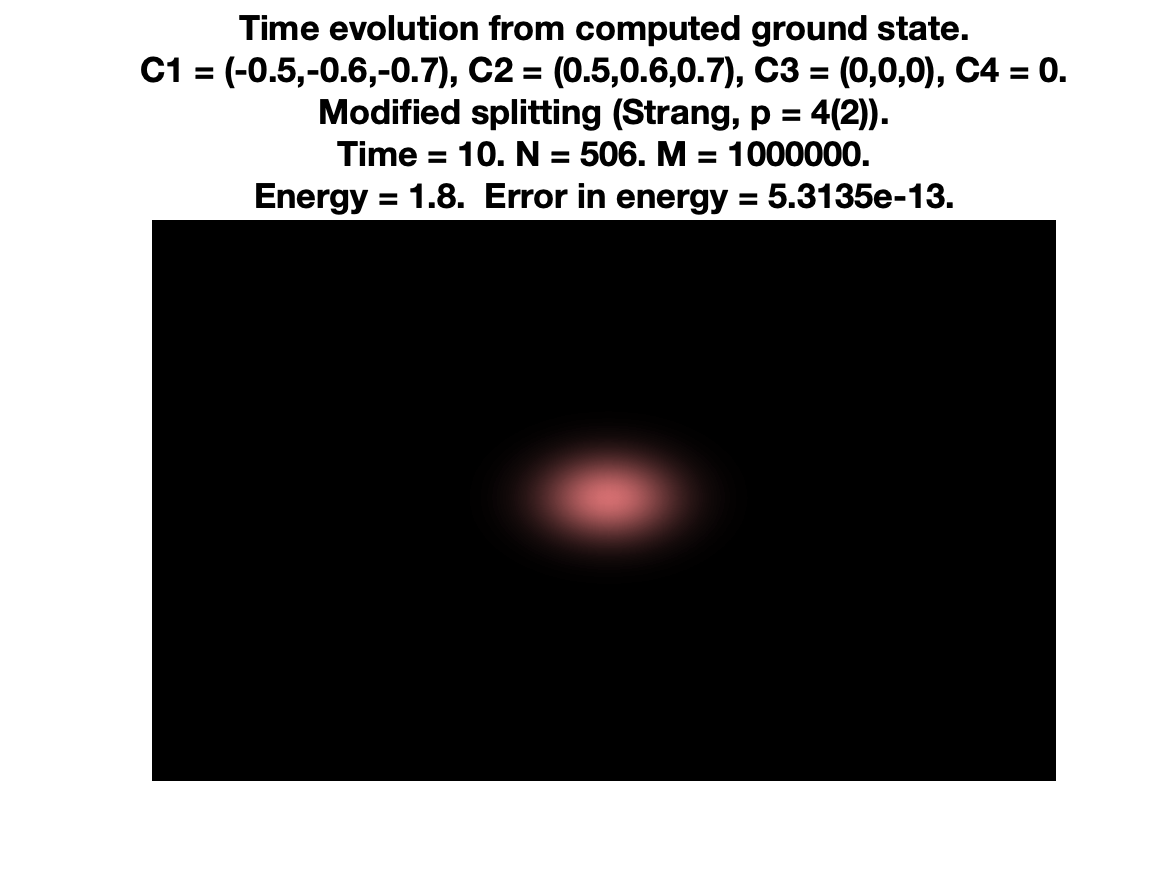} \quad
\includegraphics[width=4.2cm]{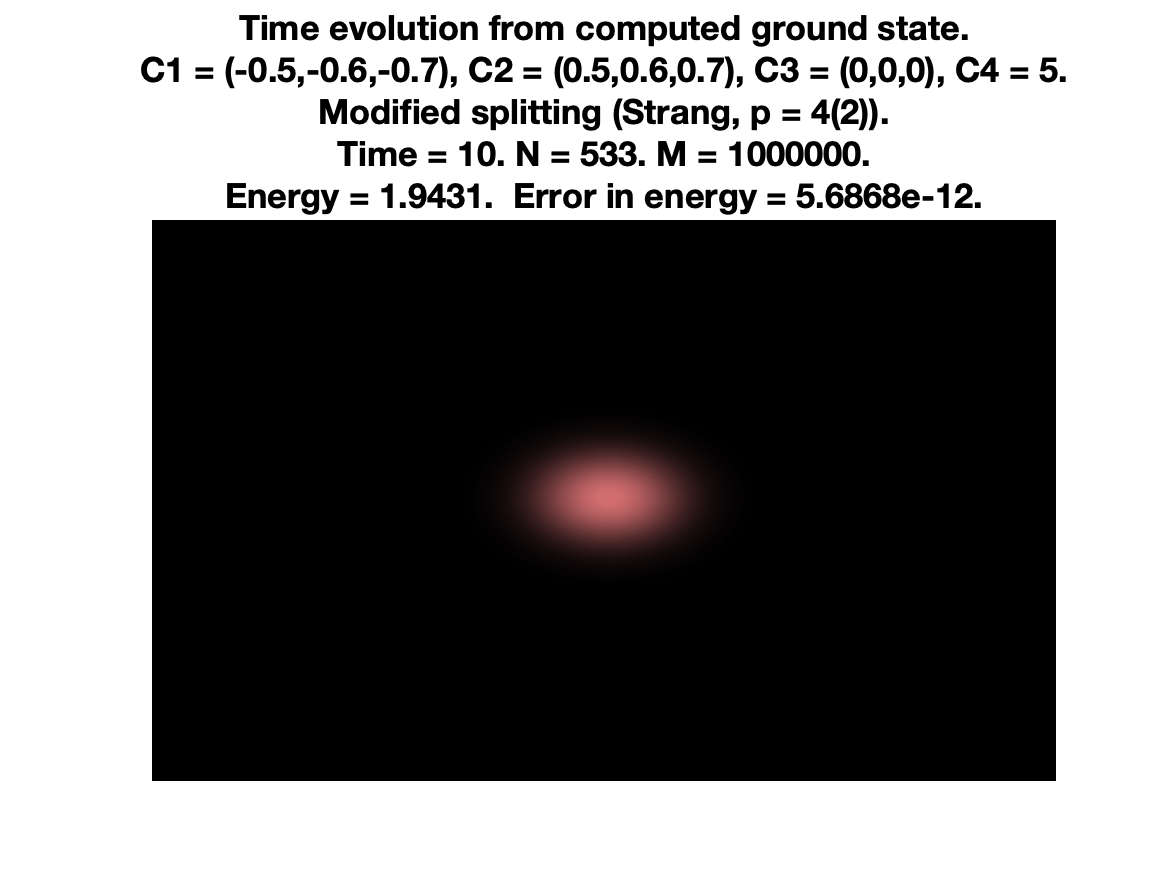} \quad
\includegraphics[width=4.2cm]{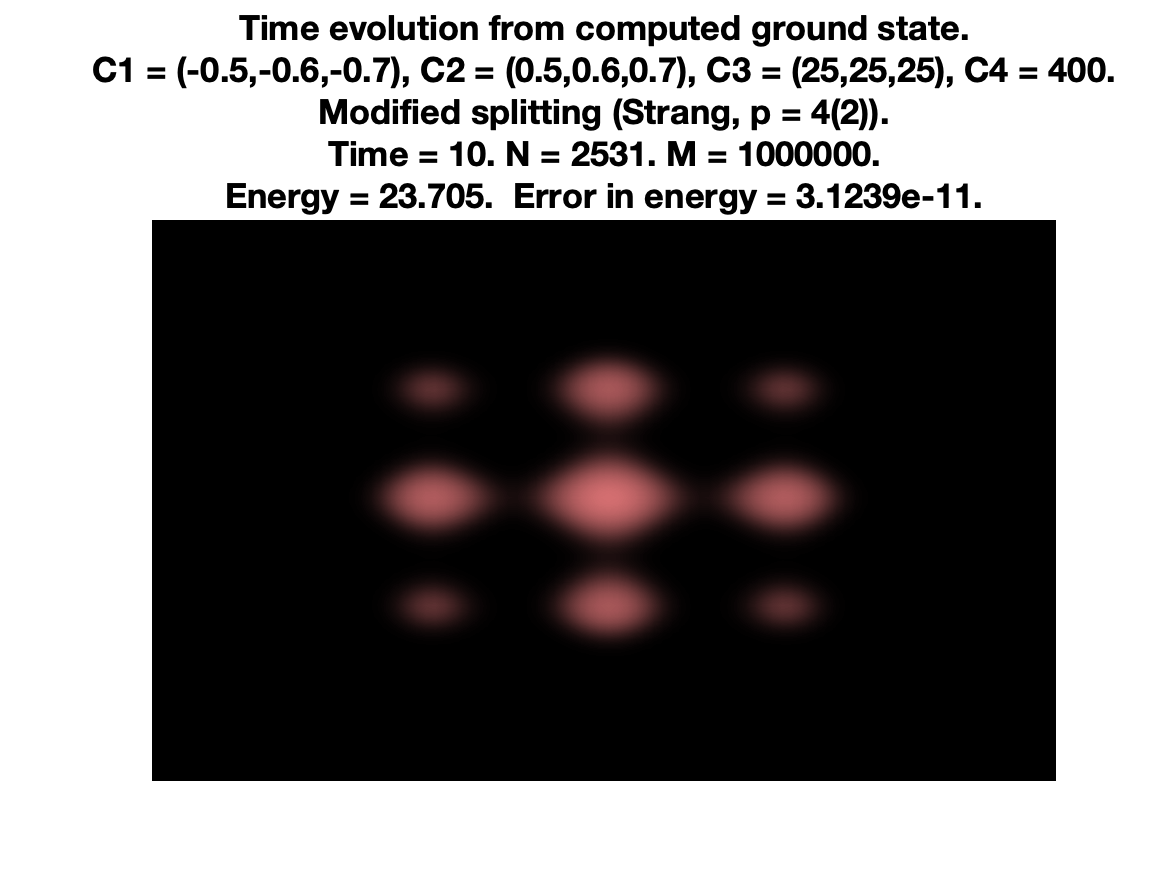} 
\end{center}
\caption{Corresponding results for three space dimensions. The solution profiles along $x_3 = 0$ resemble the results obtained in the two-dimensional case, see Figure~\ref{fig:FigureITTE2}. Differences are in particular encountered in the number of iterations and the values of the total energies.}
\label{fig:FigureITTE3}
\end{figure}

In our numerical experiments, we demonstrate the favourable stability, accuracy, and long-term behaviour of adaptive modified operator splitting methods for Gross--Pitaevskii systems in comparison with standard operator splitting methods.
As renowned instances, we consider the first-order Lie--Trotter splitting method~\eqref{eq:SchemeOrder1}, the second-order Strang splitting method~\eqref{eq:SchemeOrder2}, an optimised fourth-order splitting method proposed in~\MyCite{BlanesMoan2002}, and the fourth-order modified operator splitting method~\eqref{eq:ModifiedSplitting} for the particular choice~\eqref{eq:SchemeOrder4}.   
In the headlines of the graphics, we specify the decisive parameters, the weights in the Laplacians $C_1 = (\alpha_{ji})_{j \in \{1,\dots,J\}, \, i \in \{1, \dots, d\}}$, the weights in the harmonic potentials $C_2 = (\beta_{ji})_{j \in \{1,\dots,J\}, \, i \in \{1, \dots, d\}}$, the weights in the lattice potentials $C_3 = (\gamma_{ji})_{j \in \{1,\dots,J\}, \, i \in \{1, \dots, d\}}$, and the constants in the cubic nonlinearities $C_4 = (\vartheta_{jk})_{j,k \in \{1,\dots,J\}}$, see~\eqref{eq:GPS}. 
We implement the time evolution from the initial time $t_0 = 0$ up to a prescribed final time~$T$.
The imaginary time propagation is stopped when either the change in energy is below the tolerance $\text{Tol} = 10^{-15}$ or the total number of iterations exceeds a certain number.
We begin with numerical experiments for one-dimensional Gross--Pitaevskii equations and two-component systems involving $M = 512$ spatial grid points as well as two-dimensional Gross--Pitaevskii equations with $M = 512^2$ grid points.
The most elaborate and time-consuming three-dimensional case with $M = 100^3$ is treated in final experiments.
The estimation of the width of the Thomas--Fermi approximation yields a first indication for an appropriate choice of the underlying spatial domain. 
In view of our model problems, it is satisfactory to set $[- \, 10, 10]^d$ by default, see Figures~\ref{fig:FigureITTE1} to~\ref{fig:FigureITTE3}.

\begin{figure}[t!]
\begin{center}
\includegraphics[width=6.5cm]{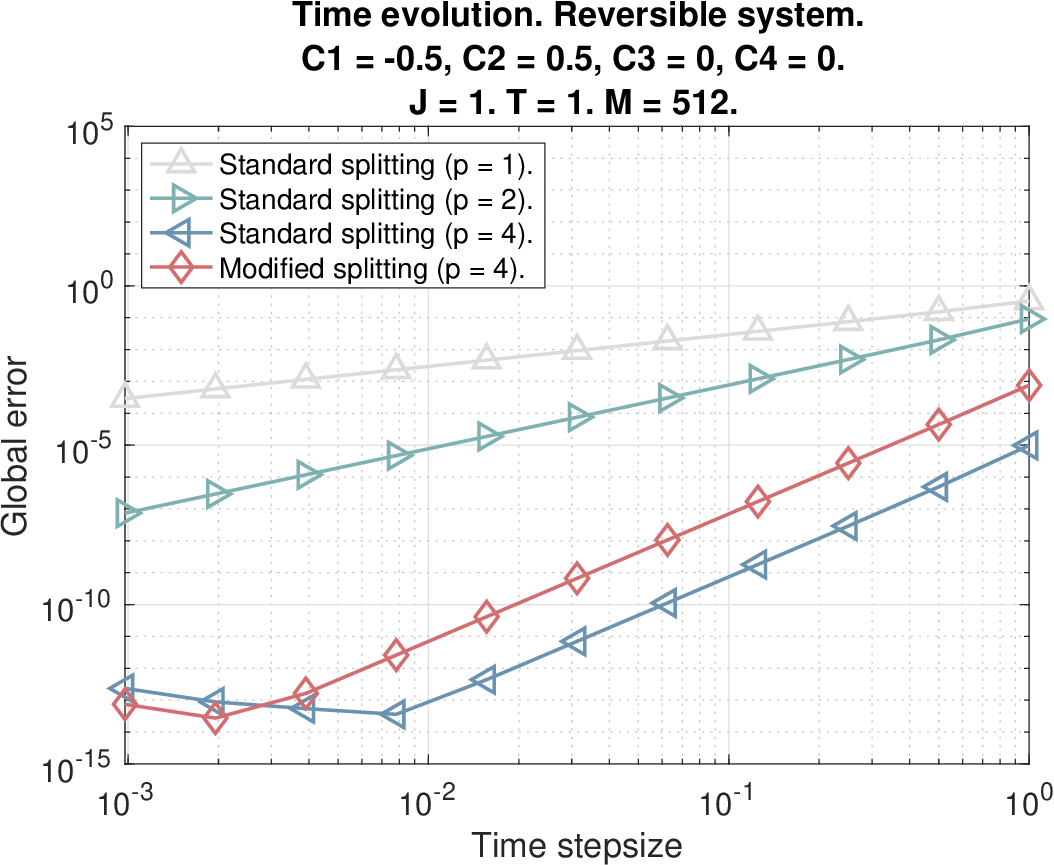} \quad
\includegraphics[width=6.5cm]{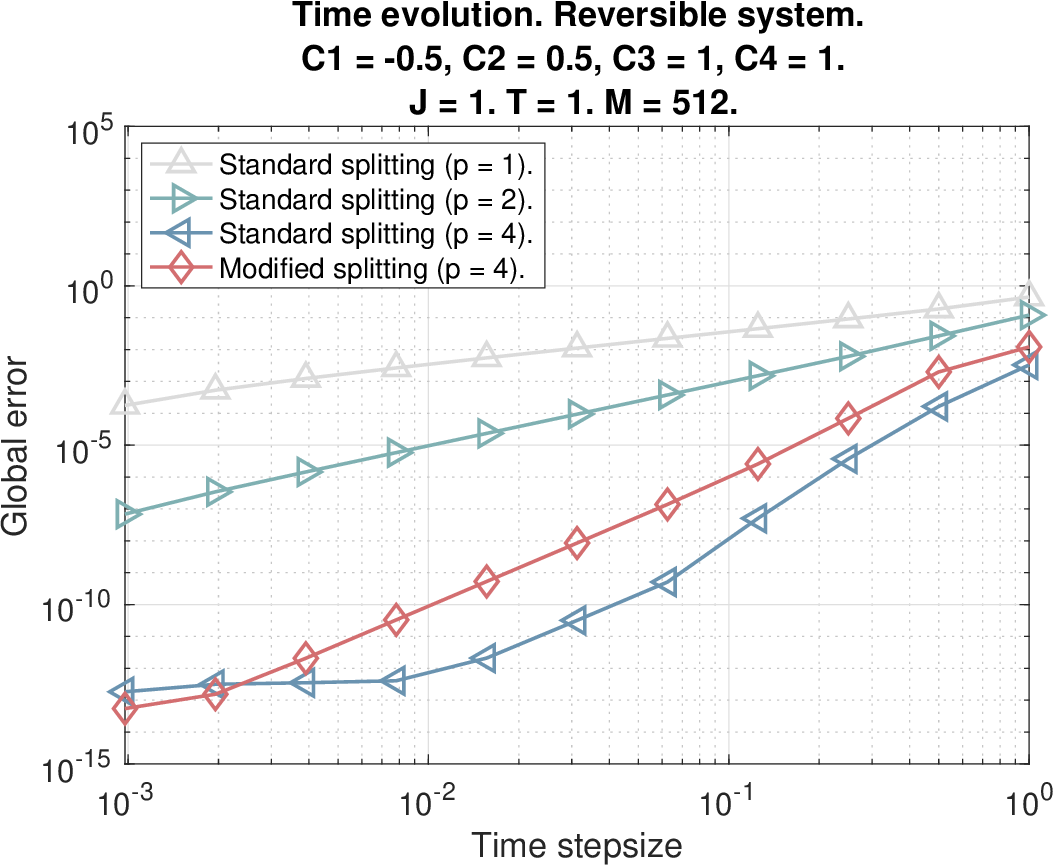} \\[2mm]
\includegraphics[width=6.5cm]{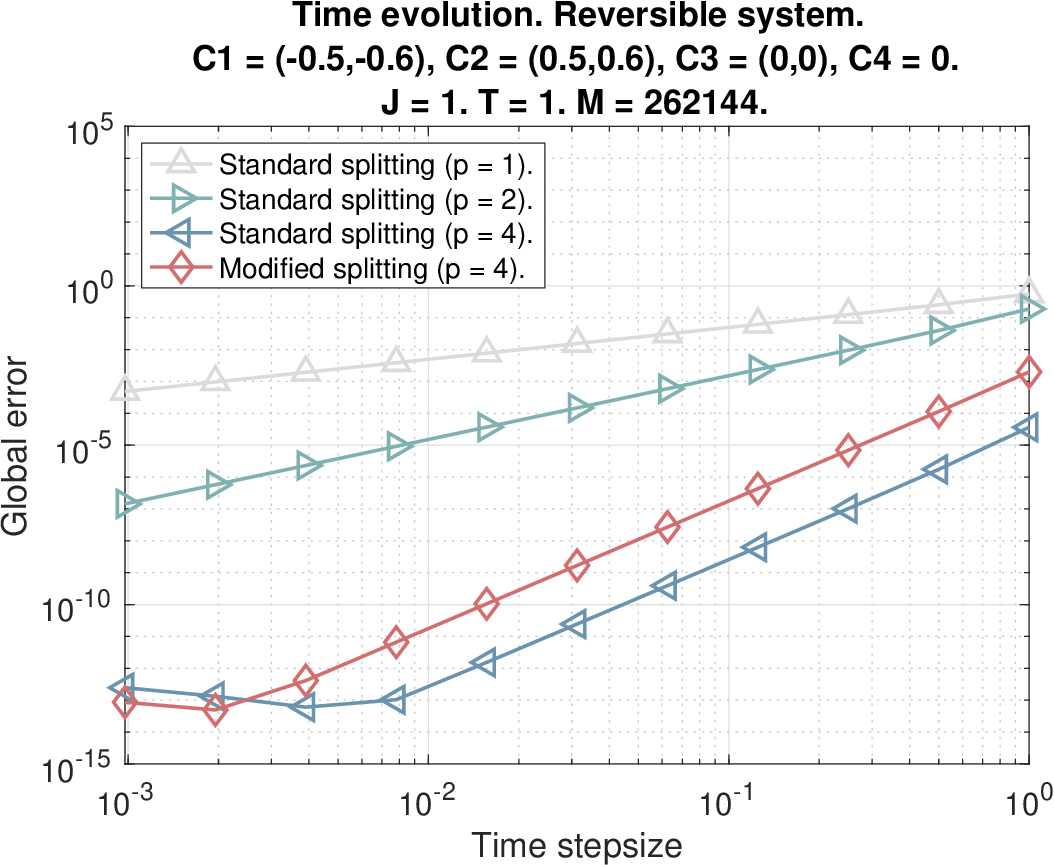} \quad
\includegraphics[width=6.5cm]{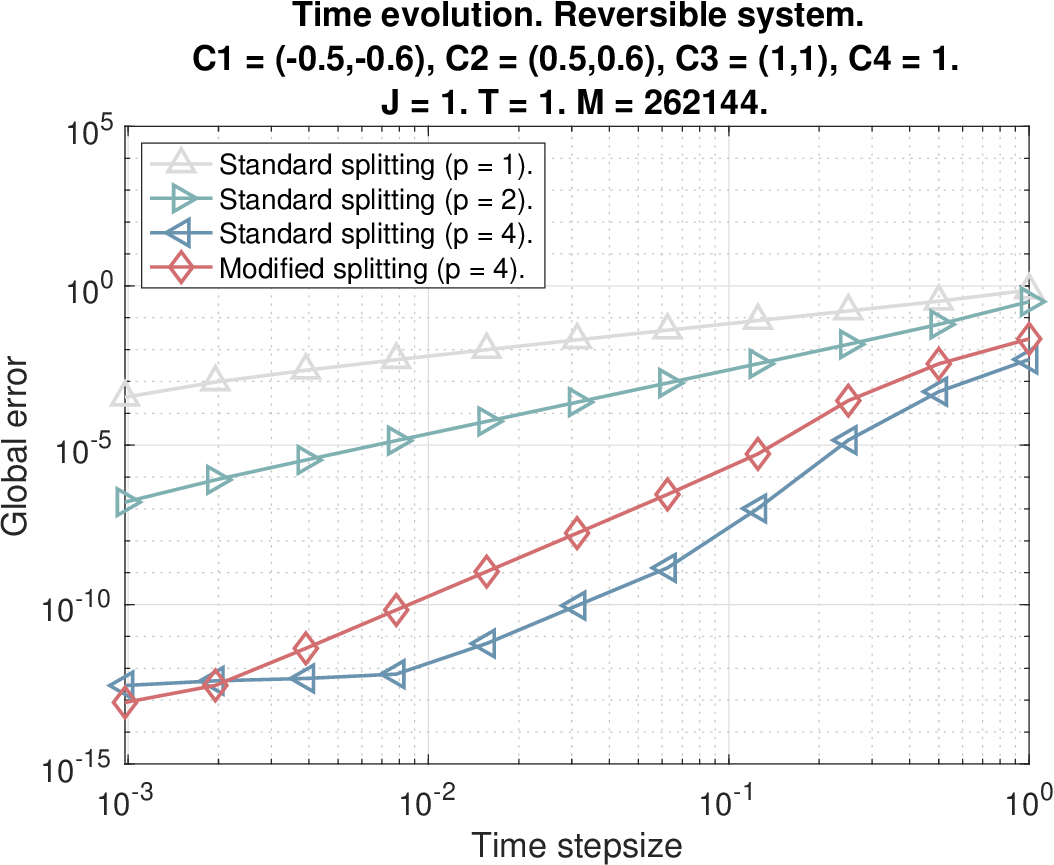} 
\end{center}
\caption{Time evolution of linear (first column) versus nonlinear (second column) reversible-in-time model problems in one (first row) and two (second row) space dimensions.
Application of standard and modified operator splitting methods, specifically the first-order Lie--Trotter splitting~\eqref{eq:SchemeOrder1}, the second-order Strang splitting~\eqref{eq:SchemeOrder2}, the optimised fourth-order splitting by \textsc{Blanes, Moan}~\cite{BlanesMoan2002}, and the fourth-order modified splitting based on~\eqref{eq:SchemeOrder4}.  
The slopes of the lines reflect the temporal orders of convergence and confirm the preservation of the nonstiff orders for problems with sufficiently regular solutions. 
Modified splitting methods remain stable and competitive in accuracy and efficiency in all test cases.}
\label{fig:FigureTO1}
\end{figure}

\begin{figure}[t!]
\begin{center}
\includegraphics[width=6.5cm]{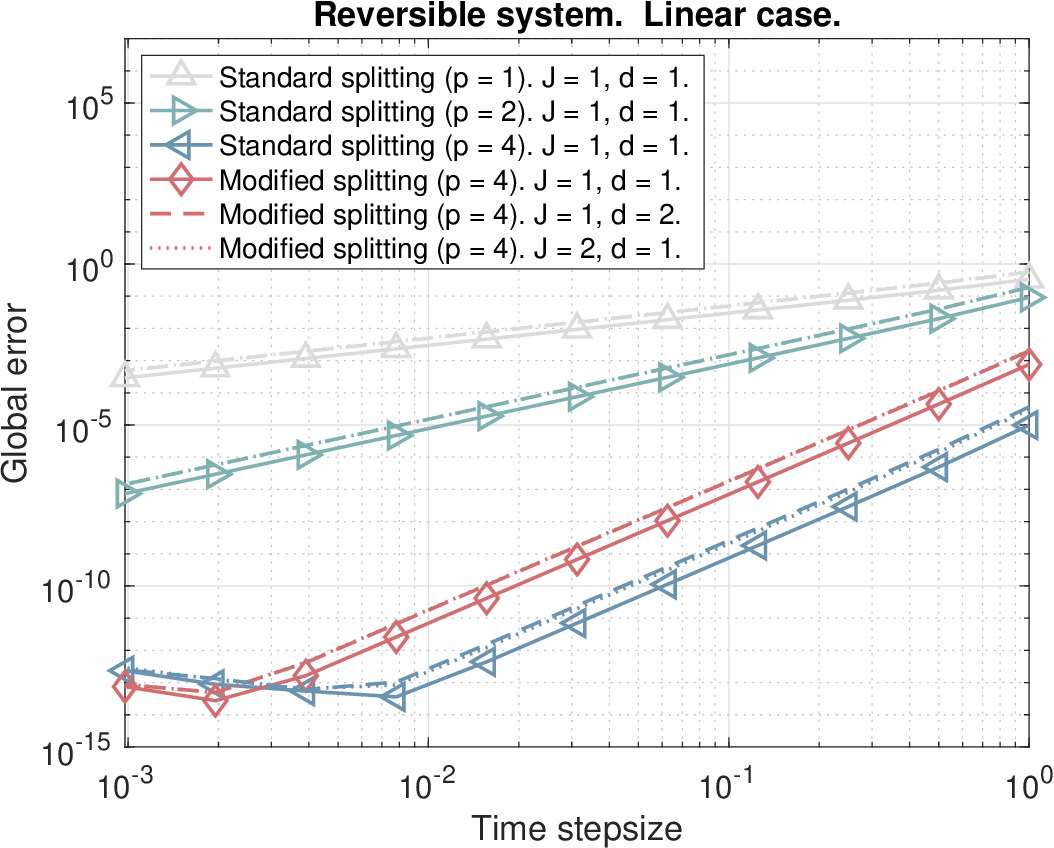} \quad
\includegraphics[width=6.5cm]{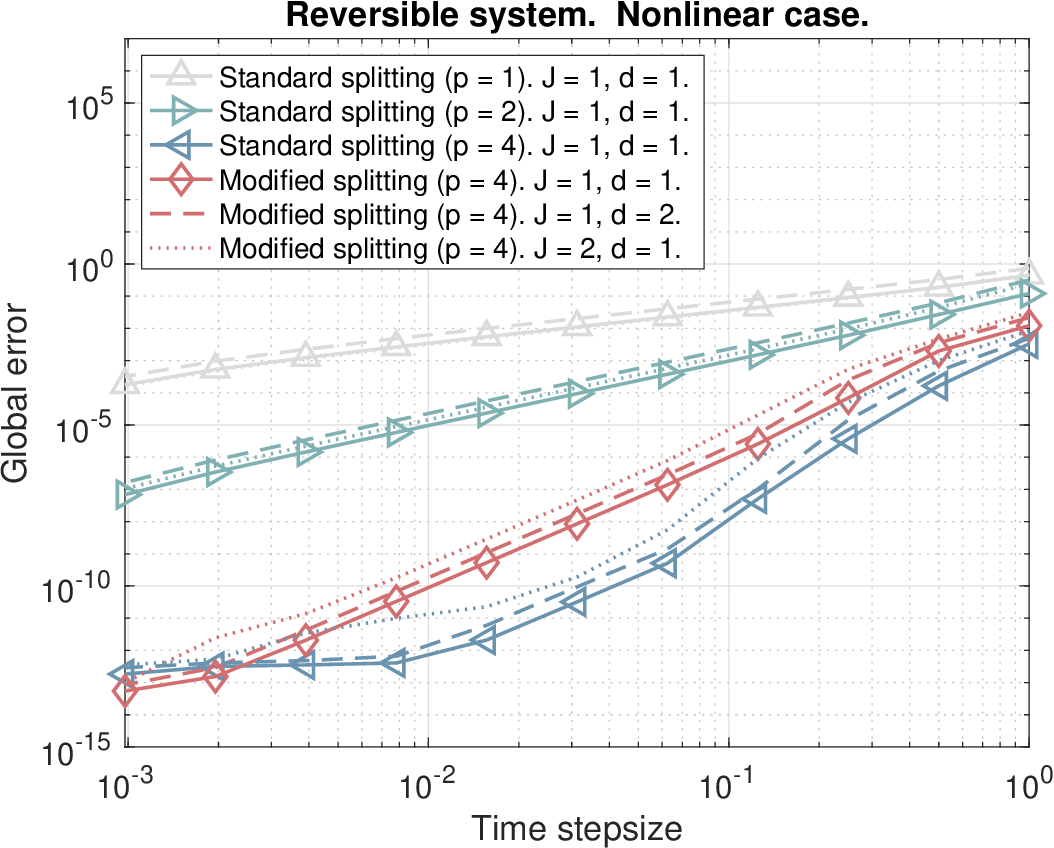} 
\end{center}
\caption{Comparisons of the results displayed in Figure~\ref{fig:FigureTO1} with related two-component systems (dotted lines).
The good agreement of the global errors for the cases $(J, d) \in \{(1,1), (1,2), (2,1)\}$ suggests a similar behaviour regarding stability and accuracy for more complex settings in three space dimensions, see Figure~\ref{fig:FigureTO5}.}
\label{fig:FigureTO2}
\end{figure}

\begin{figure}[t!]
\begin{center}
\includegraphics[width=6.5cm]{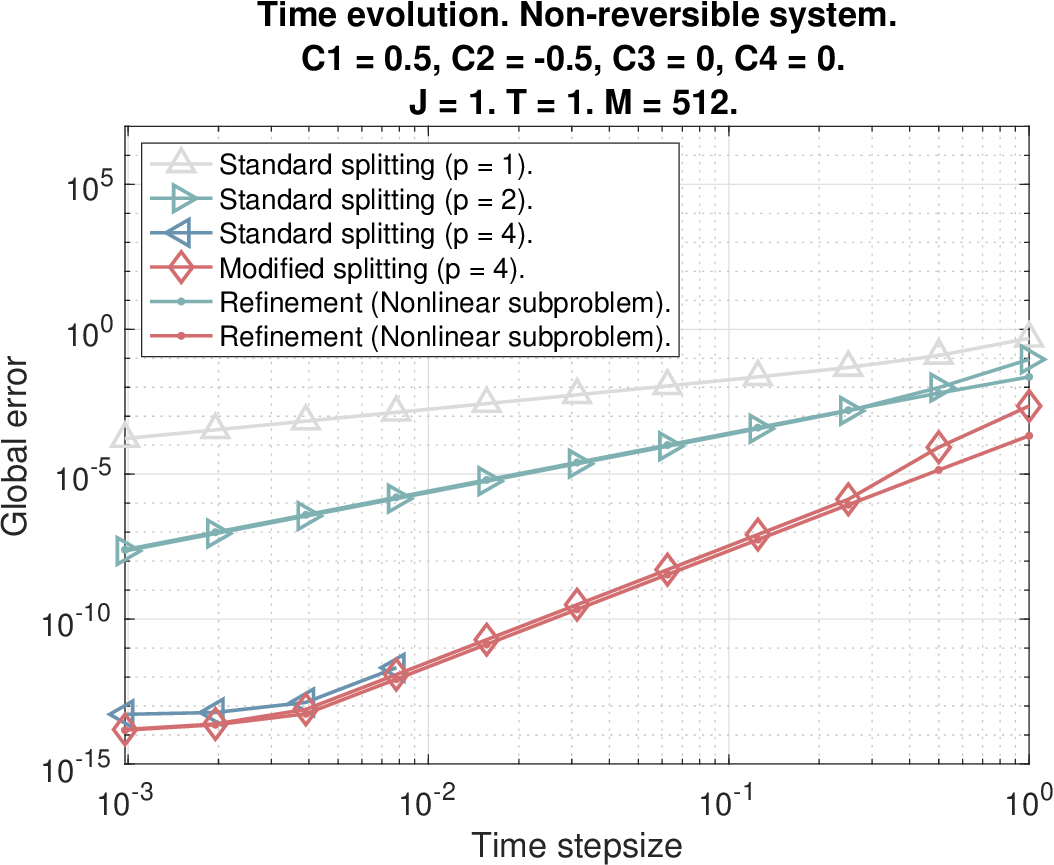} \quad
\includegraphics[width=6.5cm]{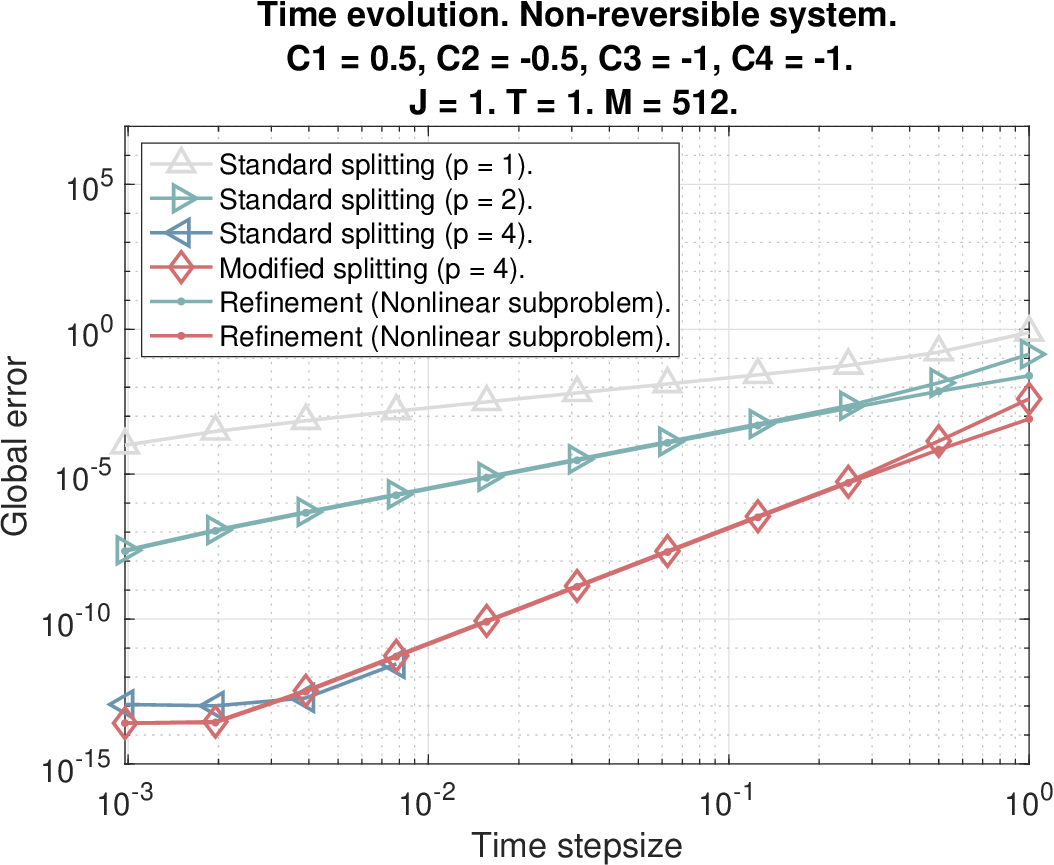} \\[2mm]
\includegraphics[width=6.5cm]{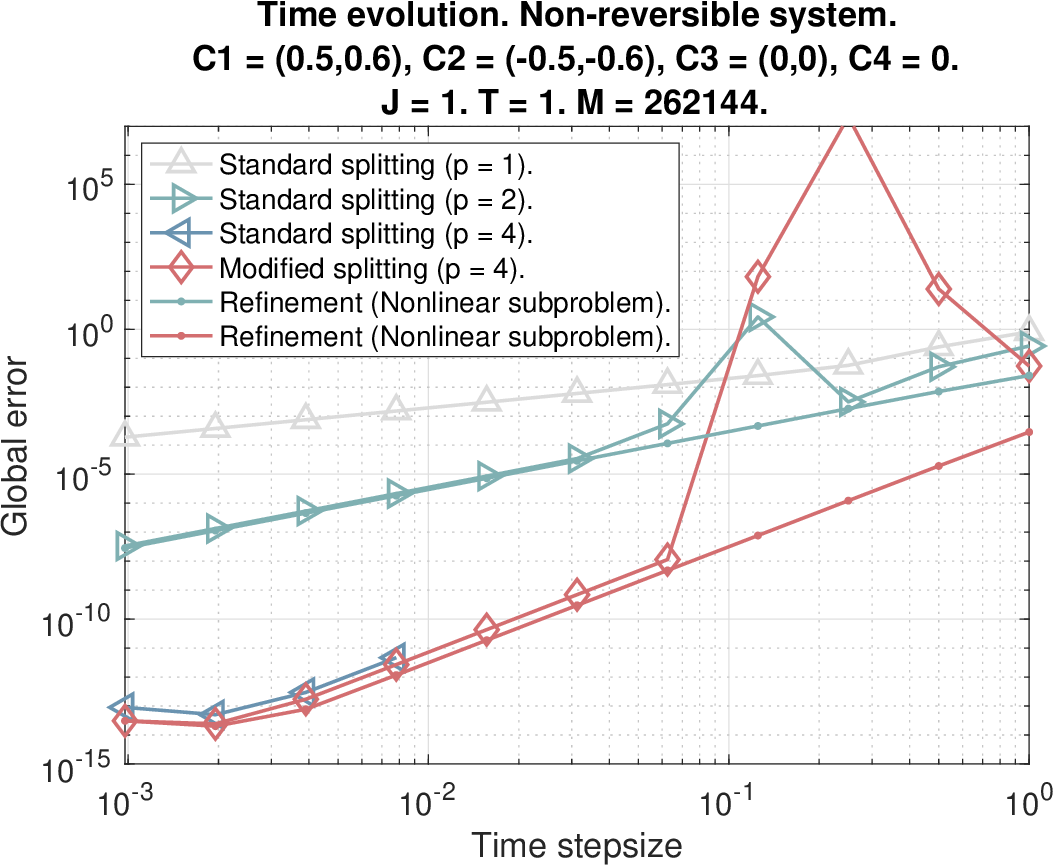} \quad
\includegraphics[width=6.5cm]{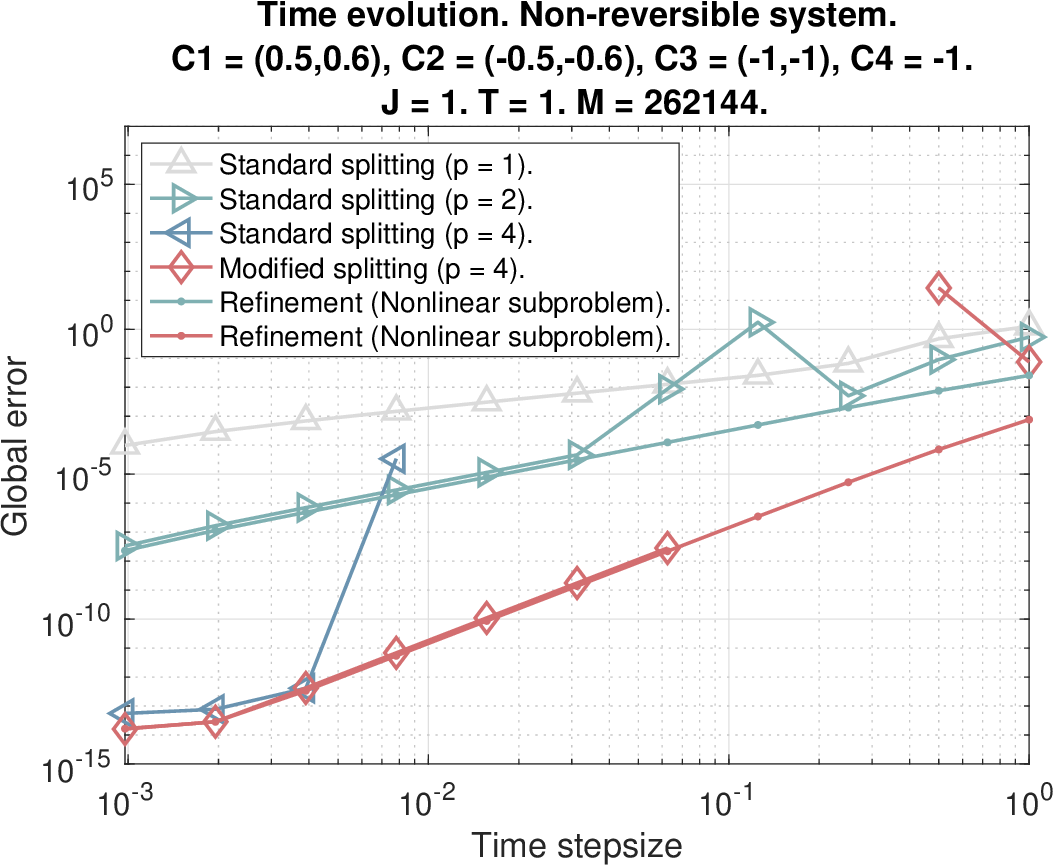} 
\end{center}
\caption{Time evolution of linear (first column) versus nonlinear (second column) non-reversible model problems in one (first row) and two (second row) space dimensions.
Application of standard and modified operator splitting methods, see Figure~\ref{fig:FigureTO1}.  
The slopes of the lines reflect the temporal orders of convergence and confirm the preservation of the nonstiff orders for problems with sufficiently regular solutions. 
For standard fourth-order splitting methods applied to non-reversible systems, severe stability issues due to the occurence of negative method coefficients and hence even failures for larger time stepsizes are observed.
On the contrary, modified splitting methods remain stable and yield highly accurate results for one space dimension (first row).
For problems in two space dimensions with increased stiffness (second row), the time integration of the subproblems comprising the potentials and the nonlinearites with refined time increments (factor $\tfrac{1}{4}$ for Strang splitting, factor $\tfrac{1}{8}$ for modified splitting) leads to significant improvements.}
\label{fig:FigureTO3}
\end{figure}

\begin{figure}[t!]
\begin{center}
\includegraphics[width=6.5cm]{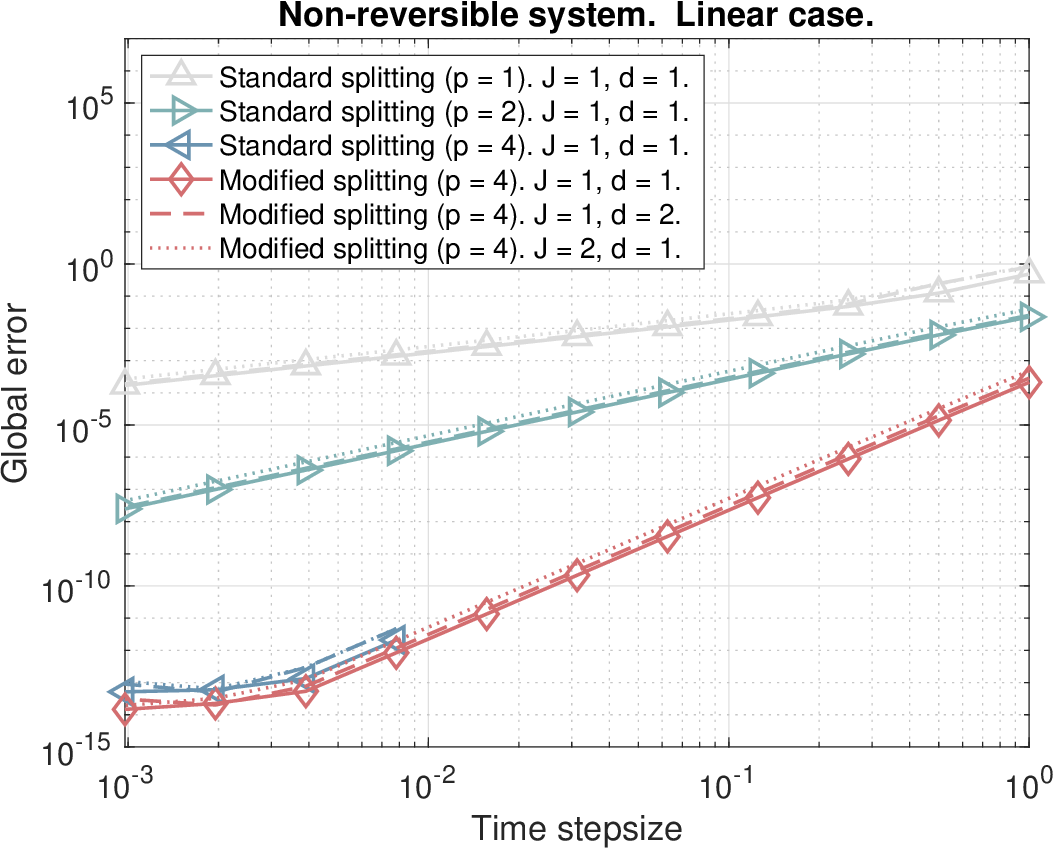} \quad
\includegraphics[width=6.5cm]{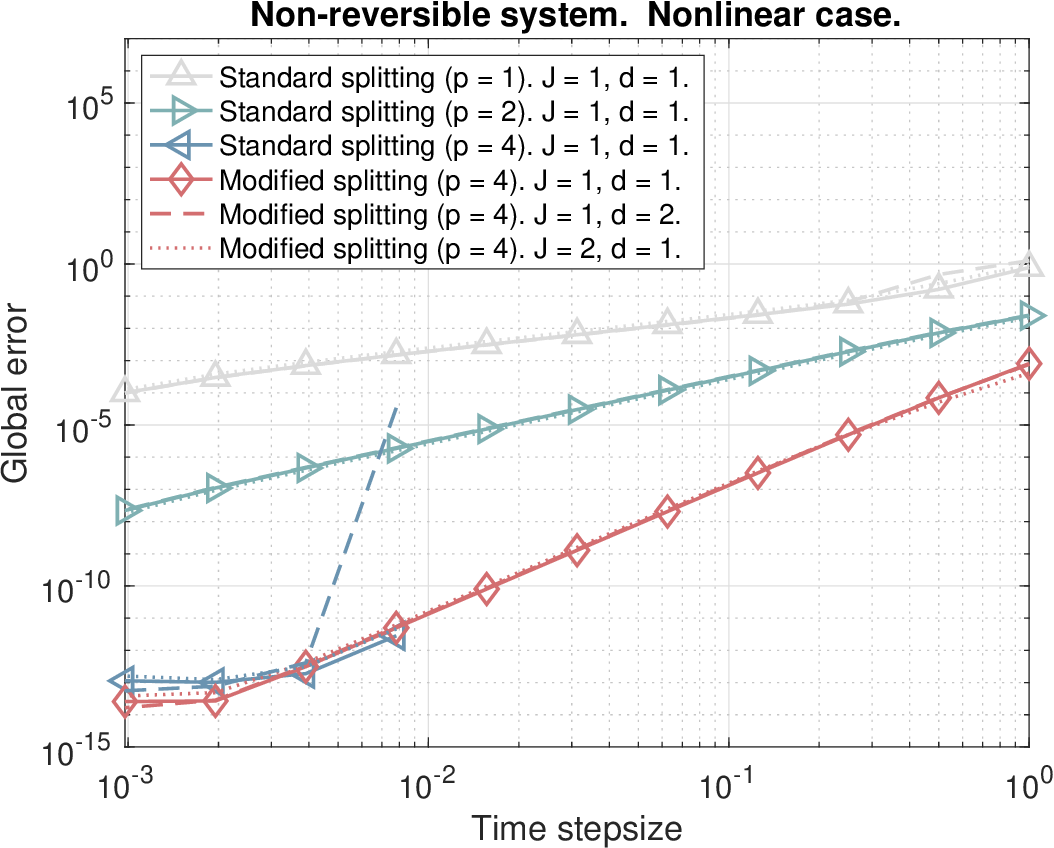} 
\end{center}
\caption{Comparisons of the results displayed in Figure~\ref{fig:FigureTO3} with related two-component systems (dotted lines).
The good agreement of the global errors for the cases $(J, d) \in \{(1,1), (1,2), (2,1)\}$ suggests a similar behaviour regarding stability and accuracy for more complex settings in three space dimensions, see Figure~\ref{fig:FigureTO5}.
Due to severe instabililies for non-reversible systems, standard higher-order operator splitting methods fail for larger time increments, whereas modified splitting methods with an appropriate resolution of the nonlinear subproblems remain stable and favourable in accuracy.}
\label{fig:FigureTO4}
\end{figure}

\begin{figure}[t!]
\begin{center}
\includegraphics[width=6.5cm]{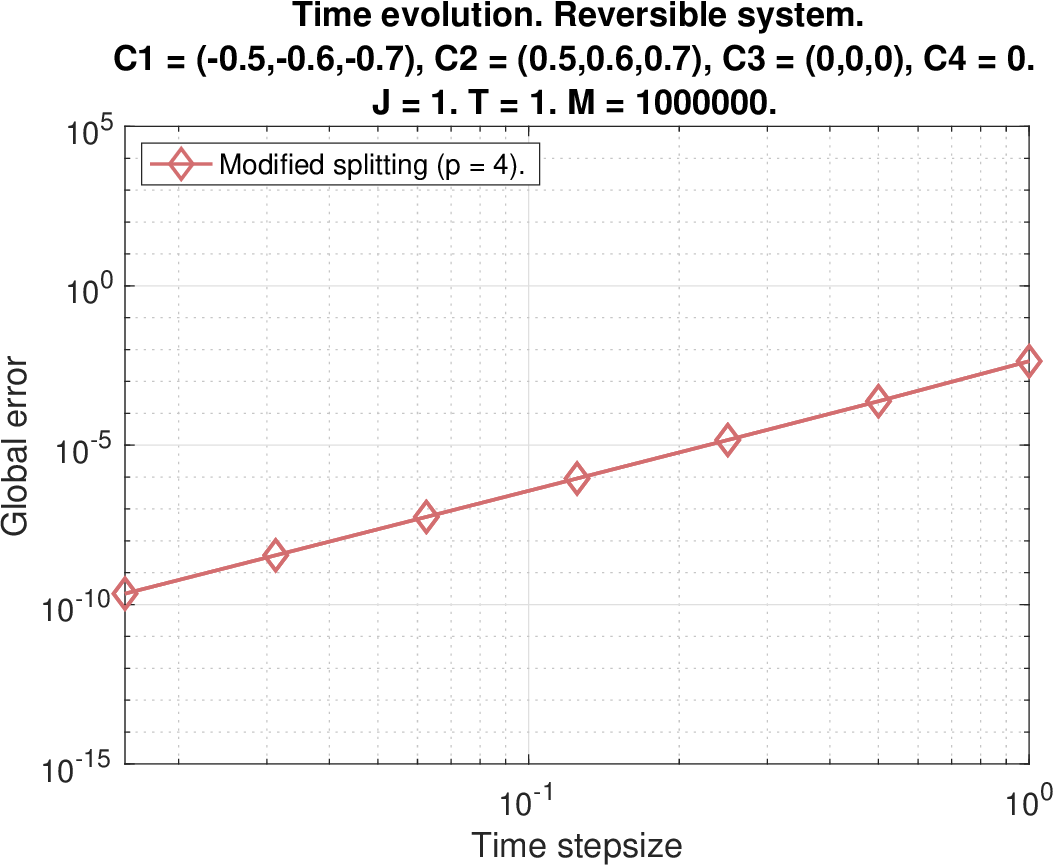} \quad
\includegraphics[width=6.5cm]{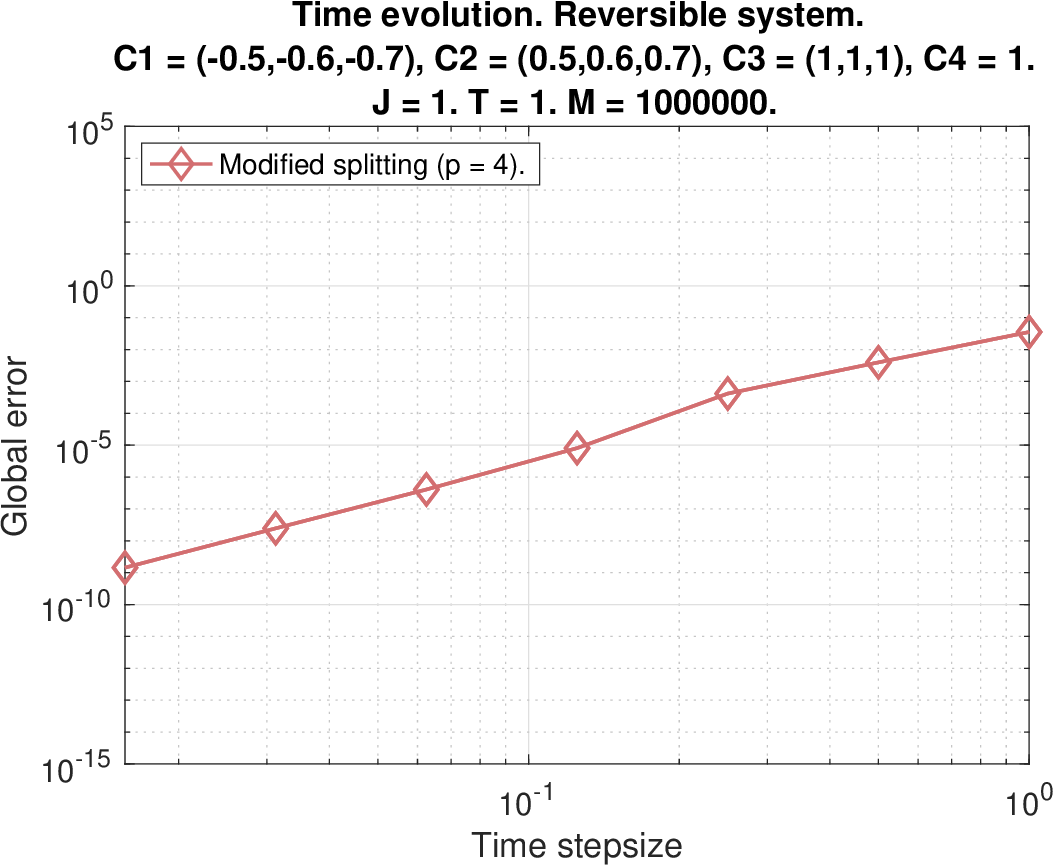} \\[2mm]
\includegraphics[width=6.5cm]{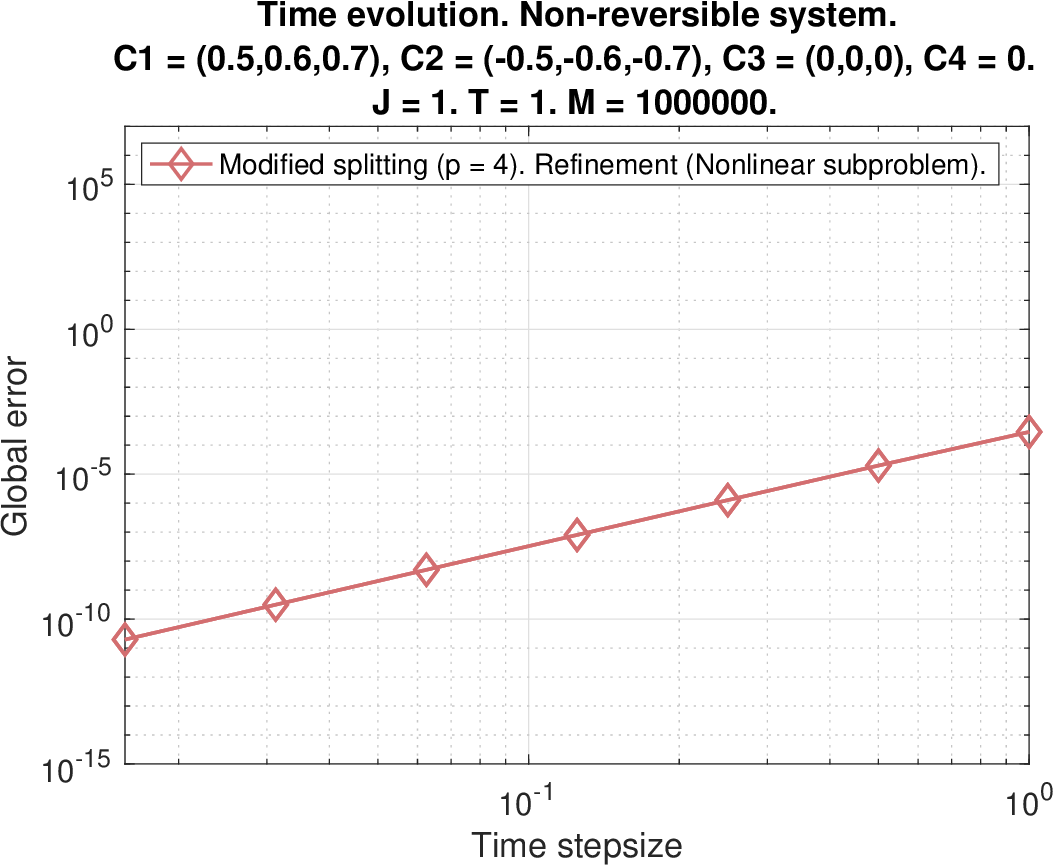} \quad
\includegraphics[width=6.5cm]{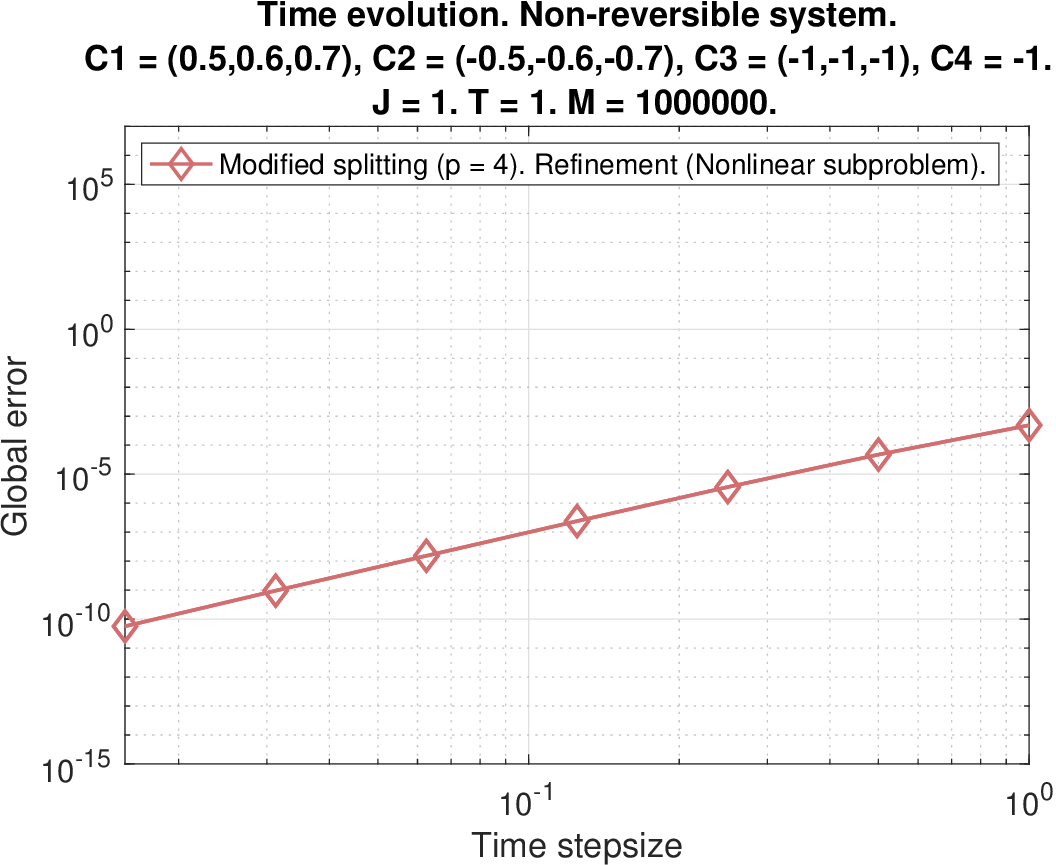} \\[2mm]
\end{center}
\caption{Confirmation of the favourable stability and accuracy behaviour of fourth-order modified splitting methods applied to linear (left column) versus nonlinear (right column) and reversible-in-time (first row) versus non-reversible (second row) model problems in three space dimensions.
In order to ensure a reliable numerical approximation of the subproblems comprising the potentials and the nonlinearites for larger time stepsizes, these stepsizes are refined by a factor $\tfrac{1}{16}$.}
\label{fig:FigureTO5}
\end{figure}

\MyParagraph{Stability and accuracy}
In a first numerical experiment, we study reversible-in-time and non-reversible model problems reflecting the type of problems that arise in the time evolution and imaginary time propagation of Gross--Pitaevskii systems.
For the purpose of verification, we make use of the fact that exact solutions are known in simplified linear cases. 
For the more complex nonlinear cases, numerical reference solutions are computed for refined time increments.
In Figures~\ref{fig:FigureTO1} to~\ref{fig:FigureTO5}, we display the global errors versus the time stepsizes in a logarithmic scale such that the slopes of the lines reflect the temporal orders of convergence.
Due to the fact that the chosen initial states and hence the solutions satisfy certain regularity and consistency requirements, the nonstiff orders of convergence are indeed retained.
Whenever highly accurate results are desirable, the fourth-order methods perform significantly better compared to the lower-order schemes, and the modified operator splitting method is competitive  in accuracy and efficiency with the optimised standard scheme.
Our main observation is that standard higher-order splitting methods applied to non-reversible systems encounter severe stability issues and even fail, whereas modified splitting methods remain stable also in higher space dimensions, provided that the time stepsizes satisfy moderate smallness assumptions.
As illustrated below, these moderate stability requirements can be fulfilled by means of an inexpensive local error control based on the second-order Strang splitting method to automatically adjust the time increments.

\begin{figure}
\begin{center}
\includegraphics[width=6.5cm]{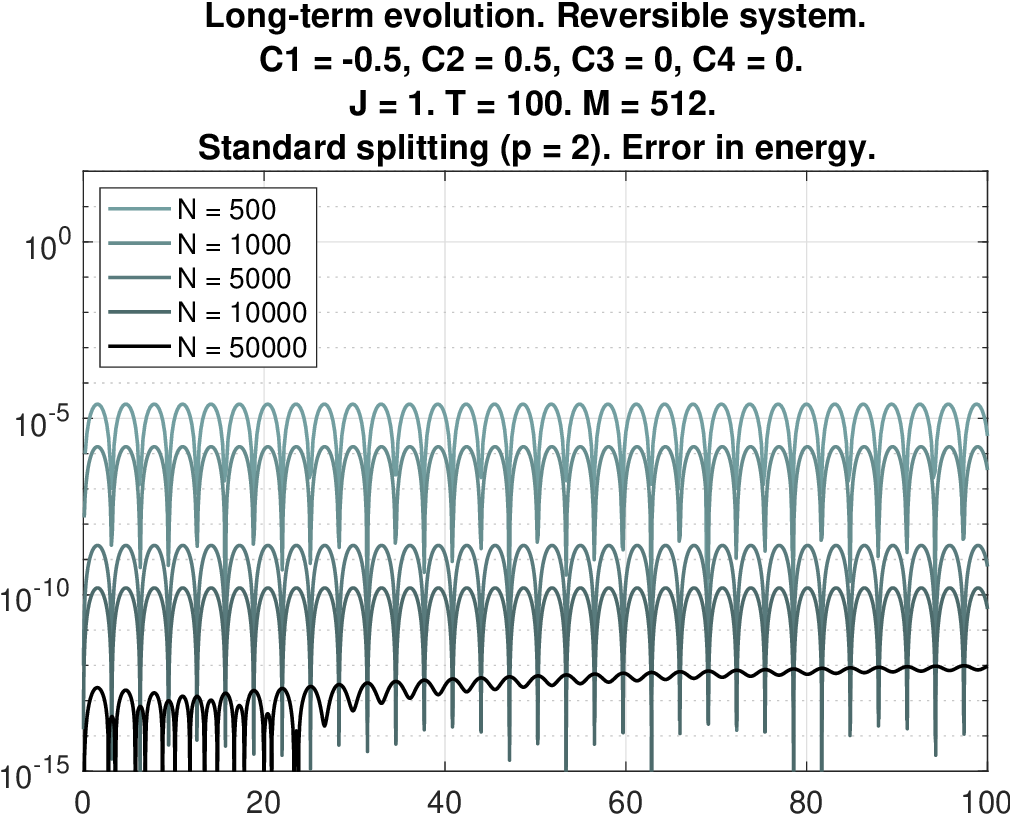} \quad
\includegraphics[width=6.5cm]{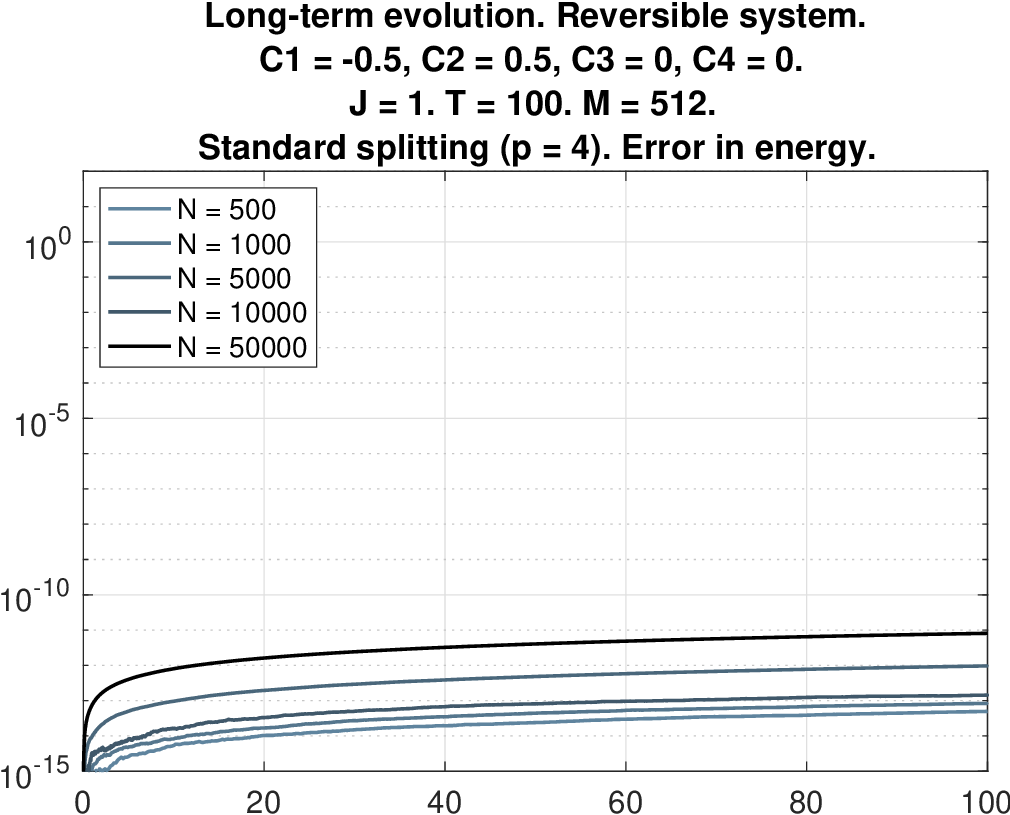} \quad
\includegraphics[width=6.5cm]{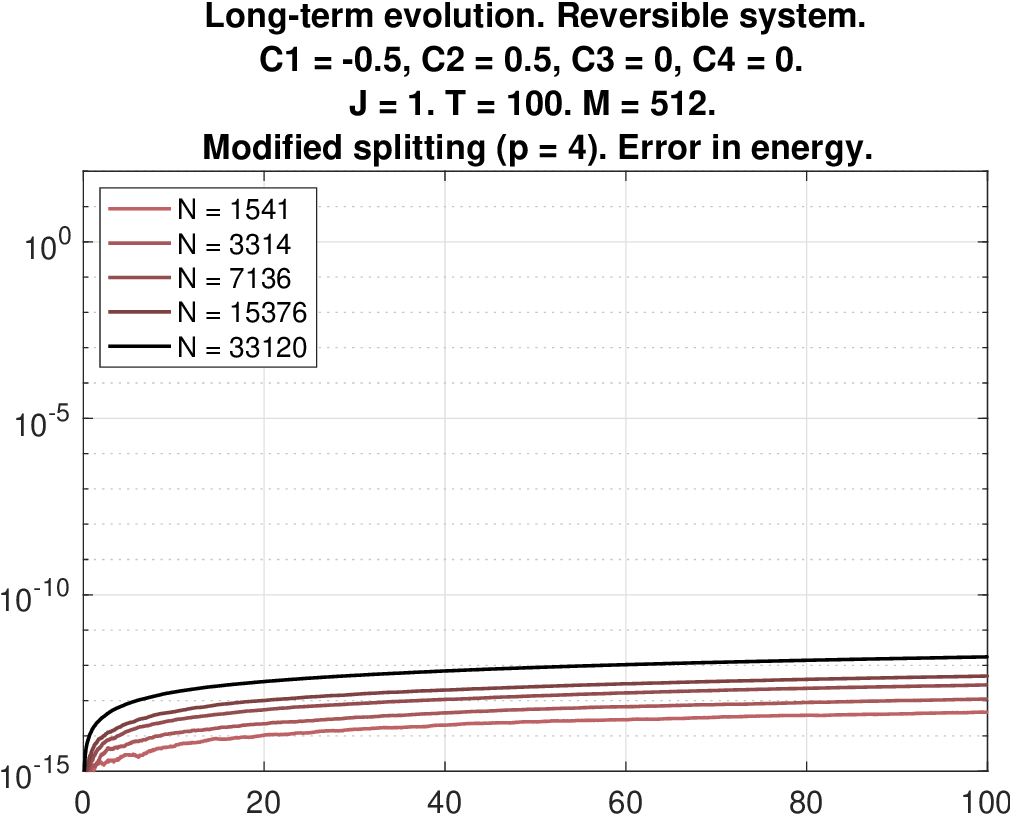} \quad
\includegraphics[width=6.5cm]{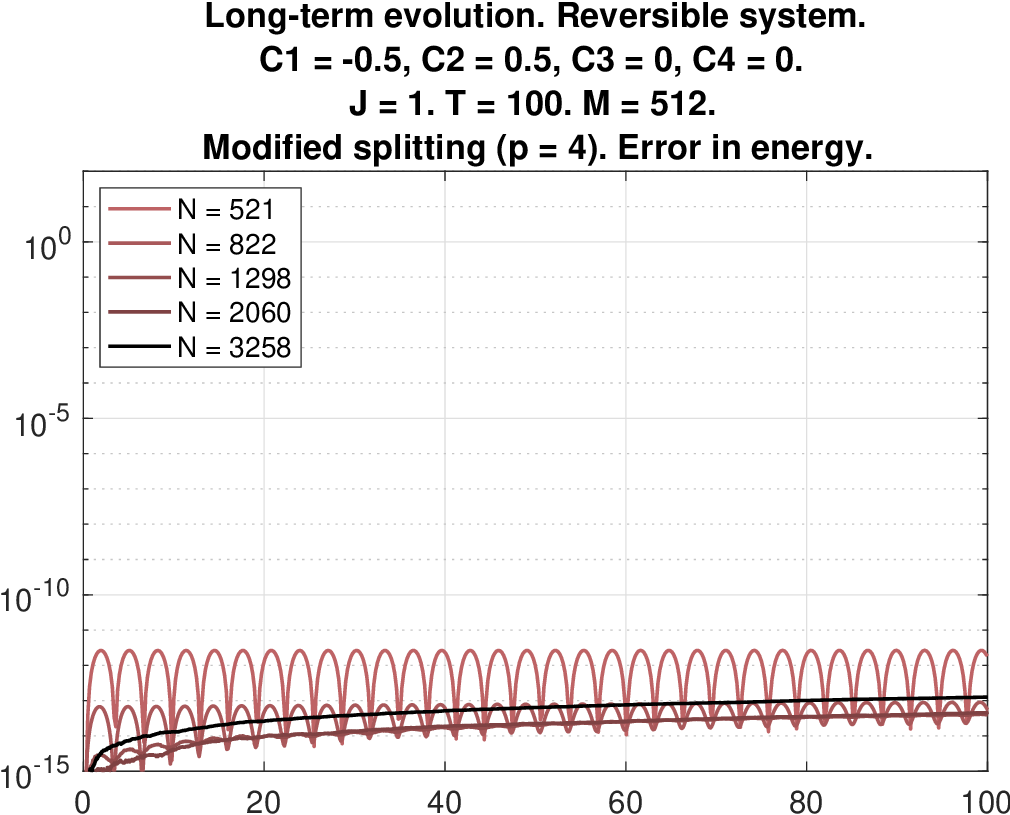} \quad
\end{center}
\caption{Energy preservation in a long-term integration of a linear Schr{\"o}dinger equation.
Application of standard splitting methods for prescribed numbers of time steps $N \in \{5 \cdot 10^2, 10^3, 5 \cdot 10^3, 10^4, 5 \cdot 10^4\}$ and resulting errors in the total energy (first row).
Application of an adaptive modified splitting method for prescribed tolerances $\text{Tol} \in \{10^{-4}, 10^{-5}, 10^{-6}, 10^{-7}, 10^{-8}\}$ with corresponding numbers of time stepsizes and resulting errors in the total energy.
On the one hand, the optimal time stepsize is determined by means of an estimate for the local error~$\text{Err}_{\text{Local}}$ based on the difference $\Psi_{\text{Modified}} - \Psi_{\text{Strang}}$, which corresponds to local order three (second row, left column).
On the other hand, the optimal time stepsize is determined by means of the local error estimate $\tau^2 \, \text{Err}_{\text{Local}}$, which corresponds to local order five and reduces the numbers of time steps (second row, right column).
The quotients of subsequent tolerances $Q_1 = 10$ and the quotients of subsequent numbers of time steps, $Q_2 = 2.15$ or $Q_2 = 1.58$, respectively, indeed confirm the local orders of convergence through $\ln(Q_1)/\ln(Q_2) \approx p + 1$.}
\label{fig:FigureTE1}
\end{figure}

\begin{figure}
\begin{center}
\includegraphics[width=6.5cm]{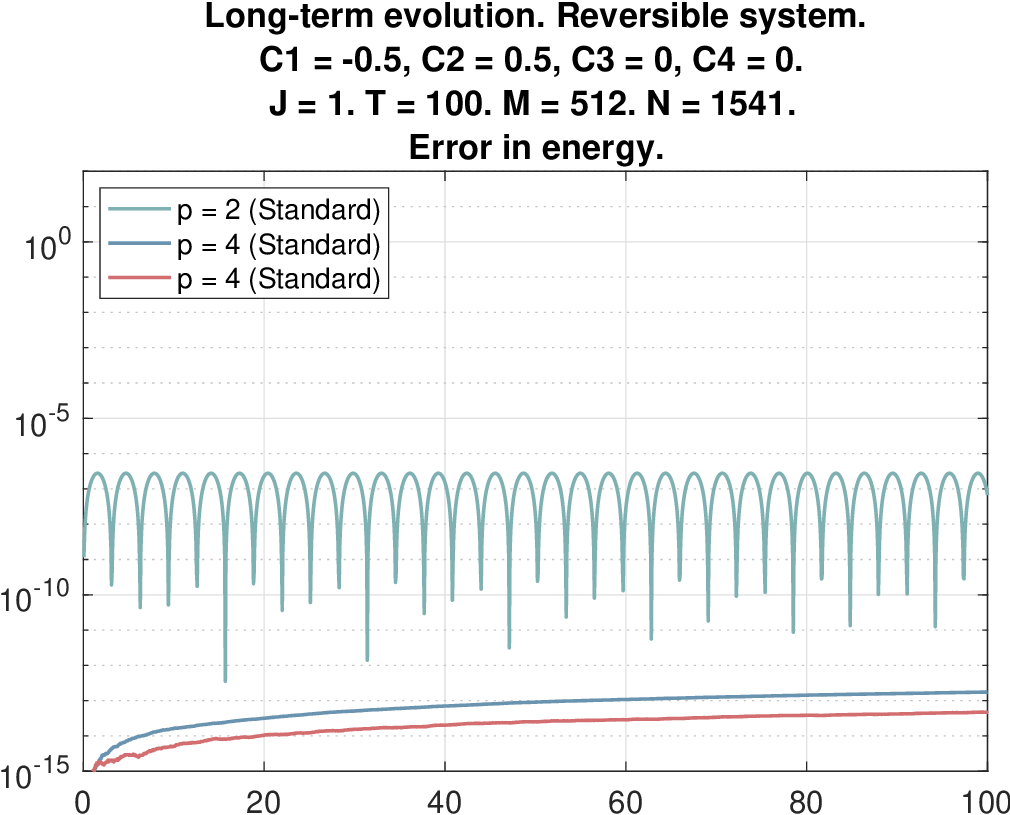} \quad
\includegraphics[width=6.5cm]{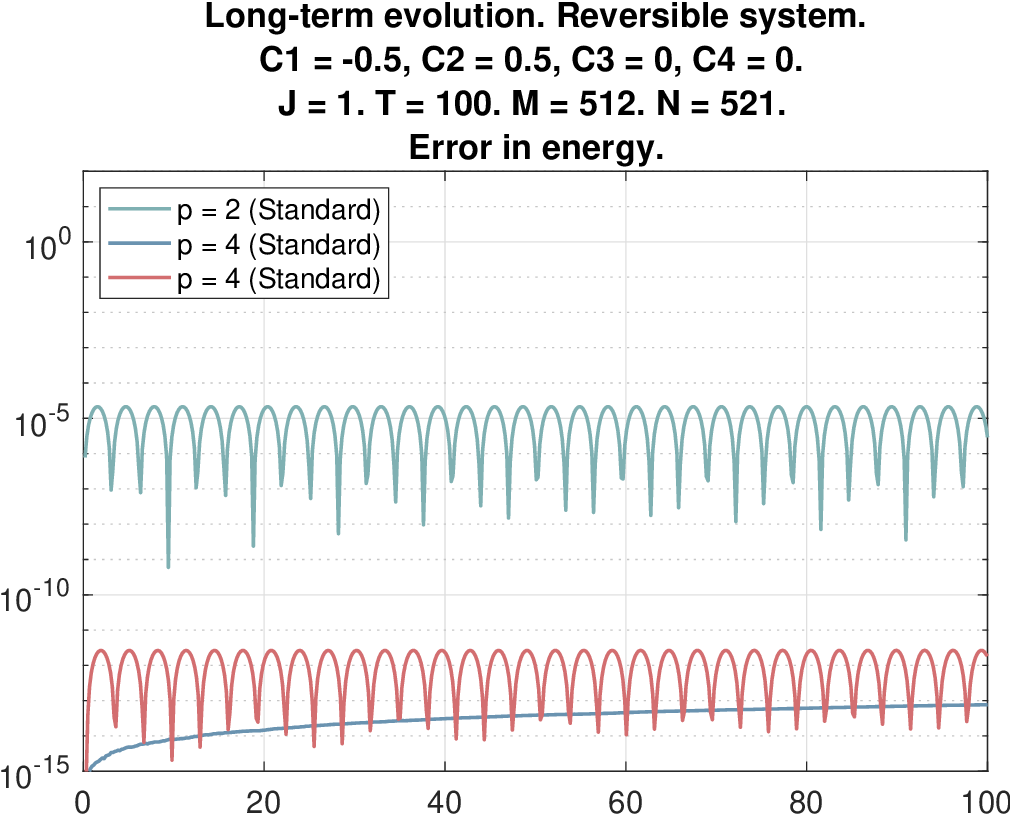} \\[2mm]
\includegraphics[width=6.5cm]{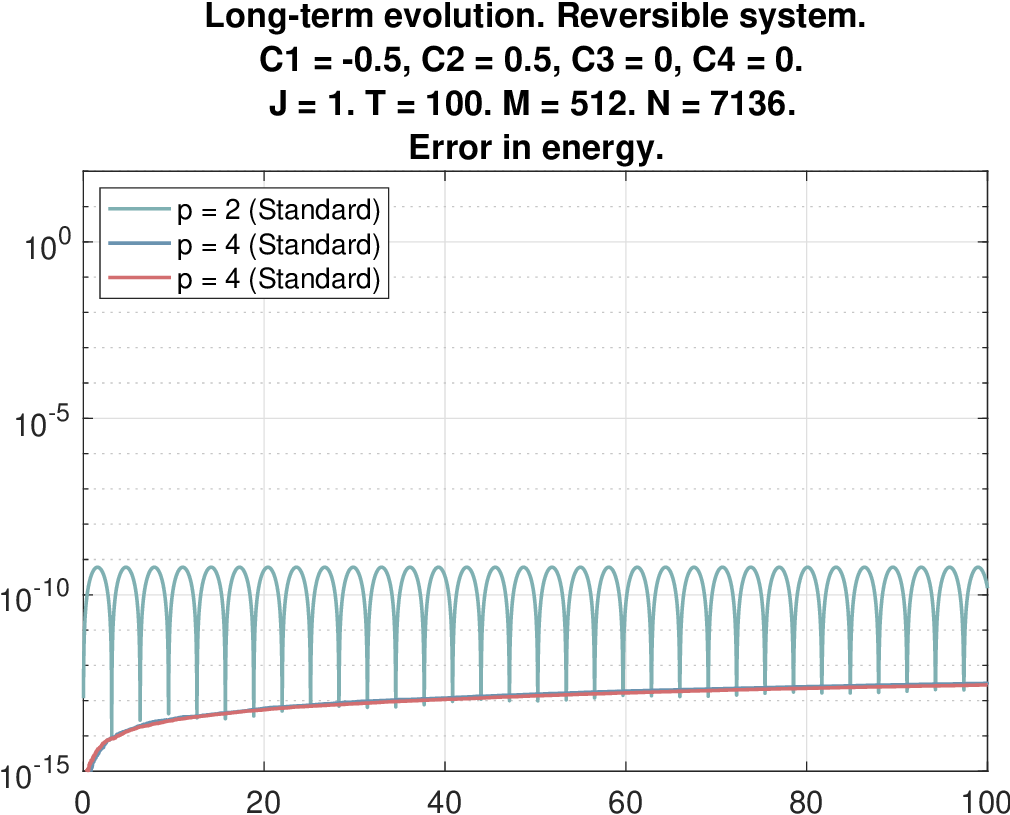} \quad
\includegraphics[width=6.5cm]{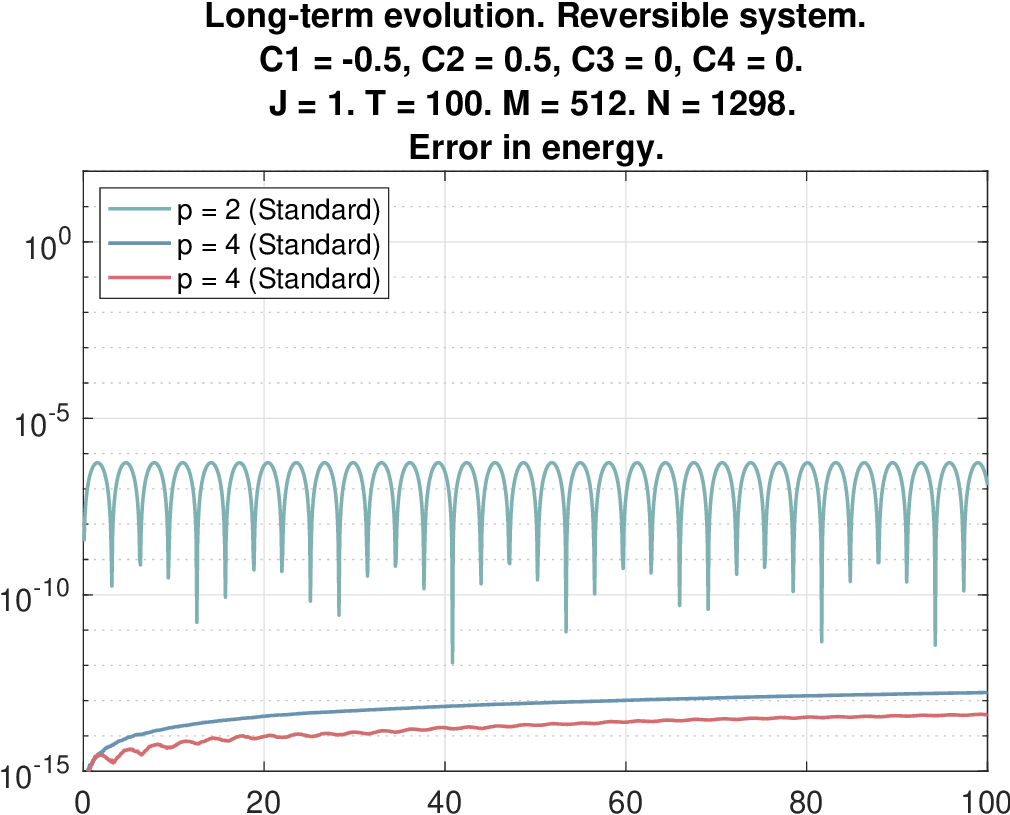} \\[2mm]
\includegraphics[width=6.5cm]{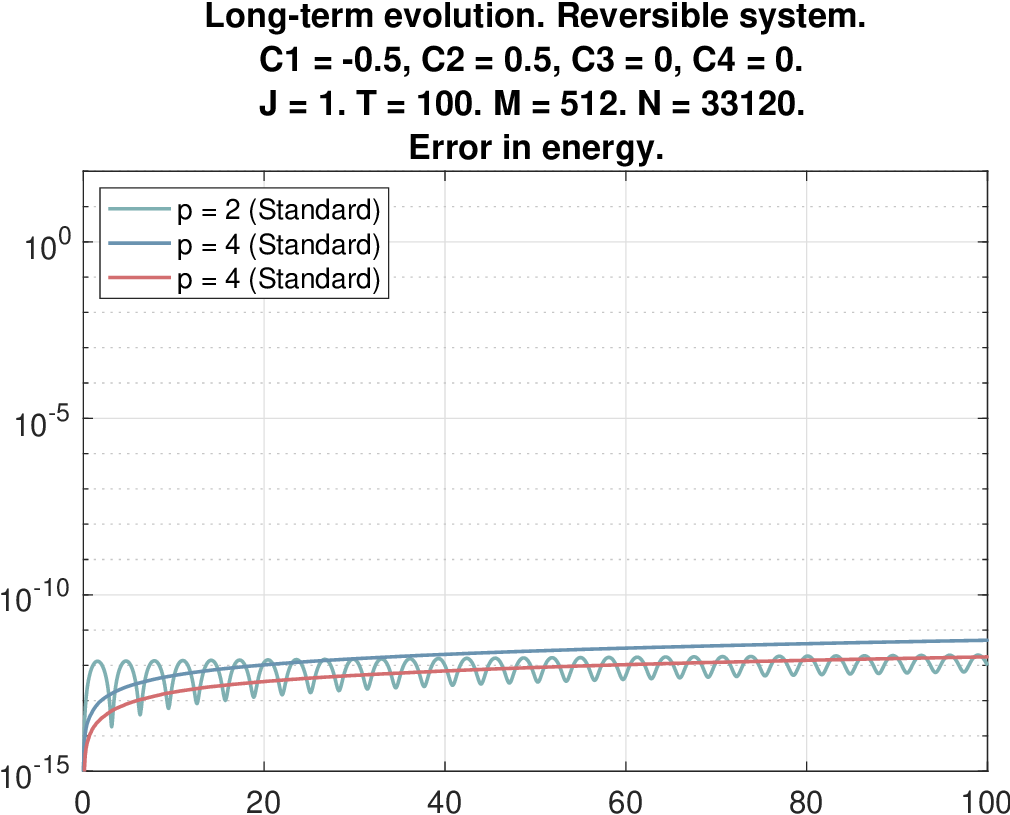} \quad
\includegraphics[width=6.5cm]{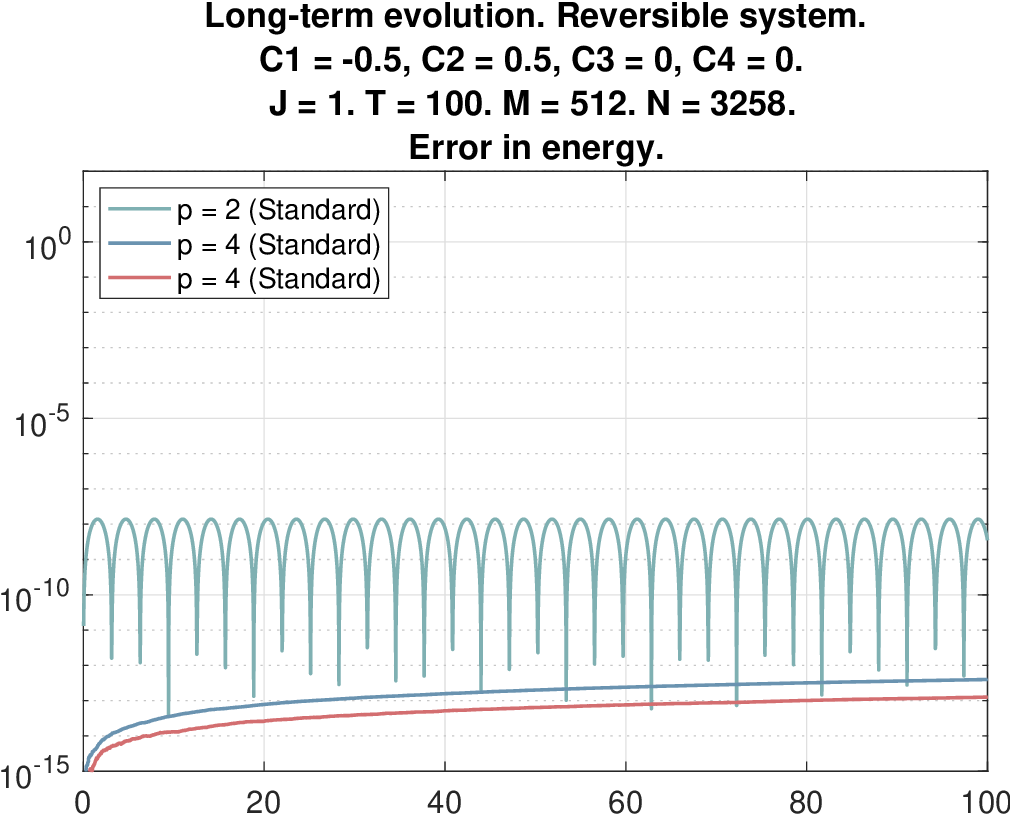} 
\end{center}
\caption{Energy preservation in a long-term integration of a linear Schr{\"o}dinger equation.
For prescribed tolerances, the time integration is performed by an adaptive fourth-order modified operator splitting method (local error estimate $\text{Err}_{\text{Local}}$ based on $\Psi_{\text{Modified}} - \Psi_{\text{Strang}}$ for left column,  local error estimate $\tau^2 \, \text{Err}_{\text{Local}}$ for right column).
Subsequently, for the resulting numbers of time steps and corresonding equidistant time increments, standard operator splitting methods of orders two and four are applied.
The results confirm the favourable performance of adaptive modified splitting methods, in particular when high accuracy is desirable.}
\label{fig:FigureTE2}
\end{figure}

\begin{figure}
\begin{center}
\includegraphics[width=6.5cm]{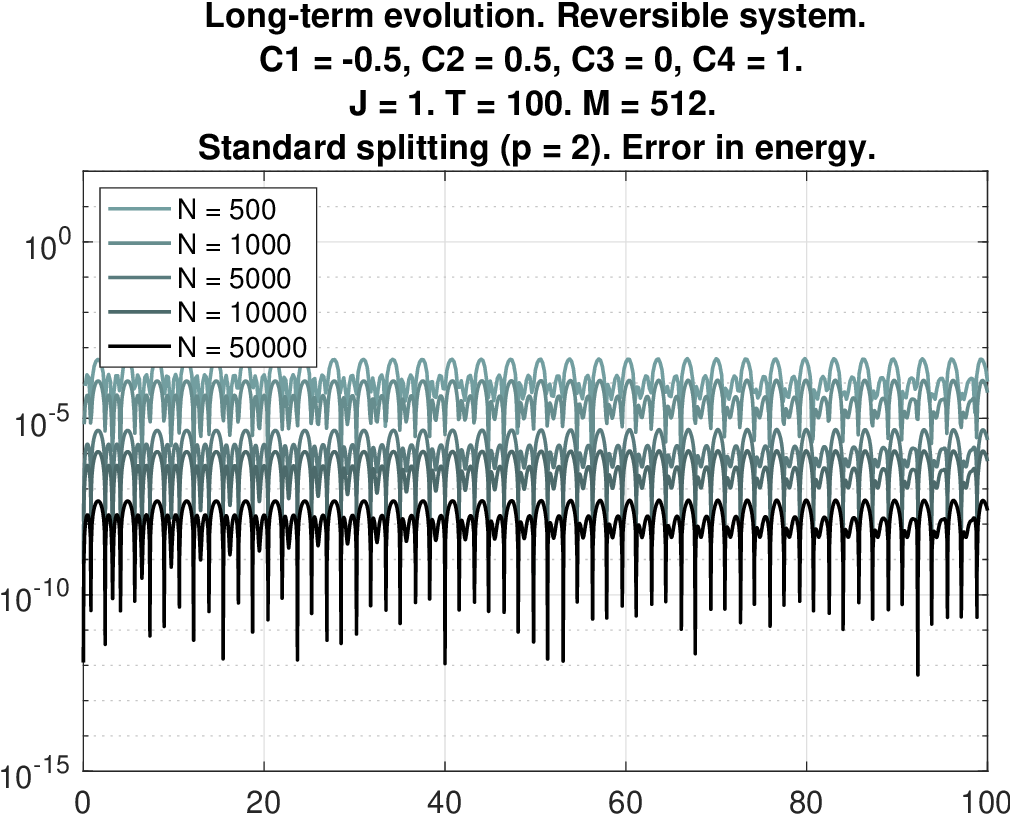} \quad
\includegraphics[width=6.5cm]{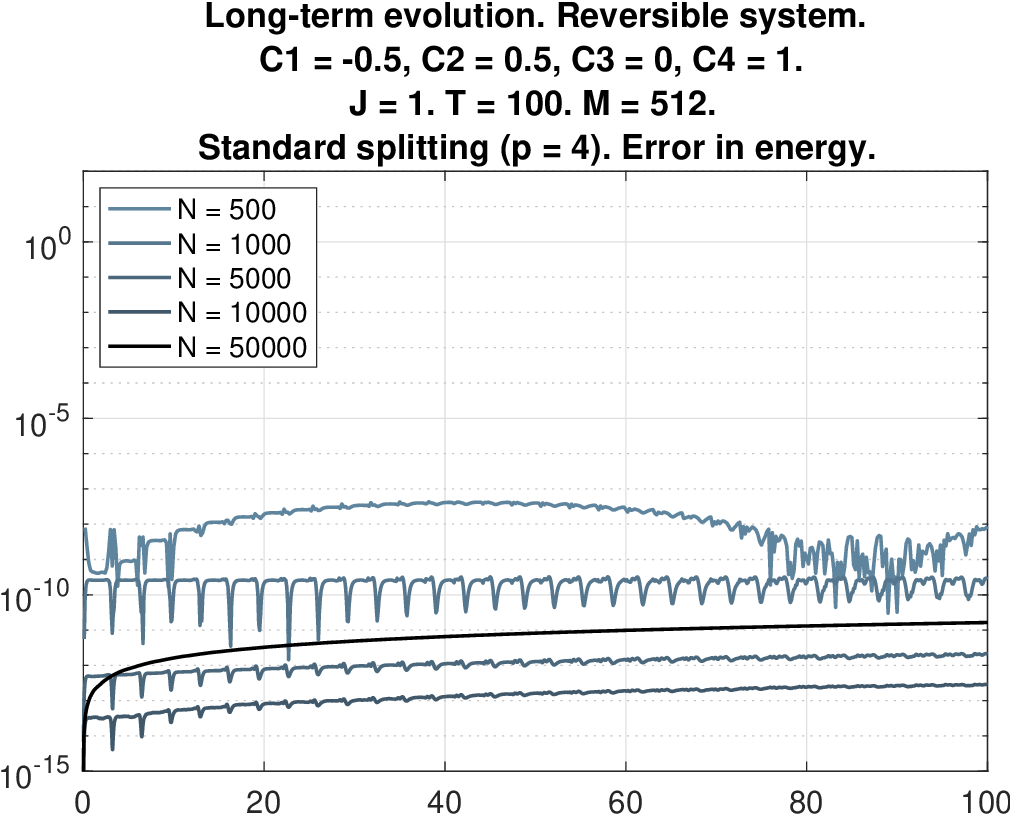} \\[2mm]
\includegraphics[width=6.5cm]{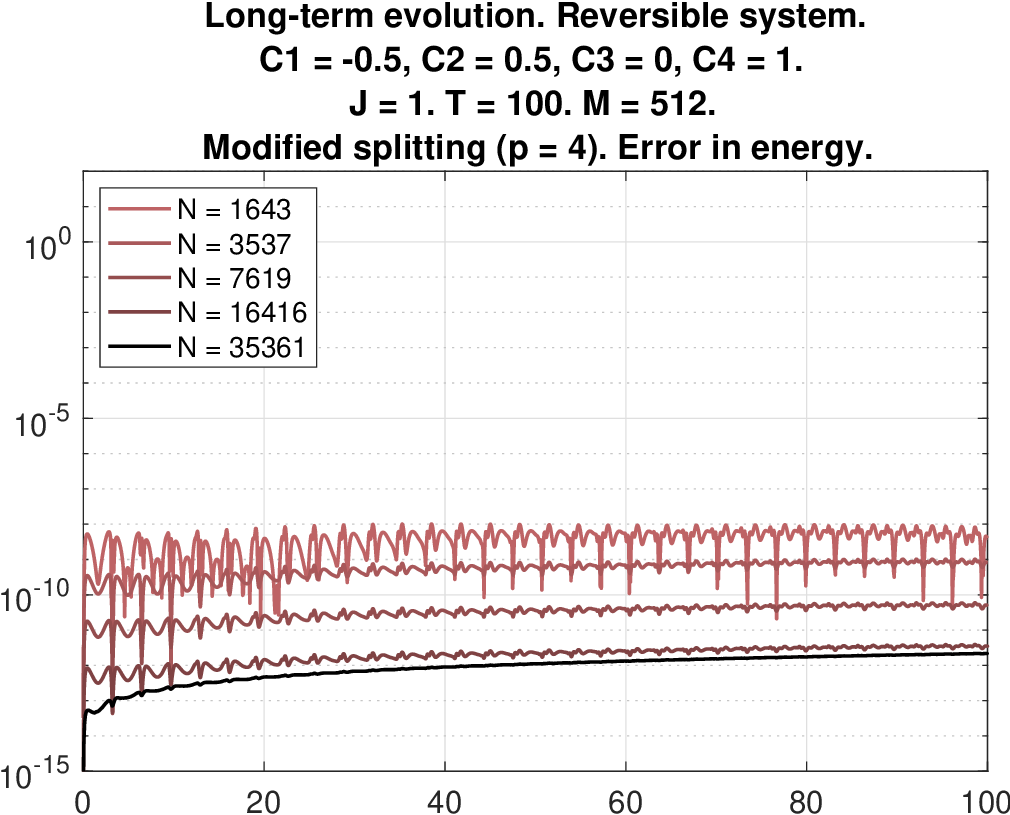} \quad
\includegraphics[width=6.5cm]{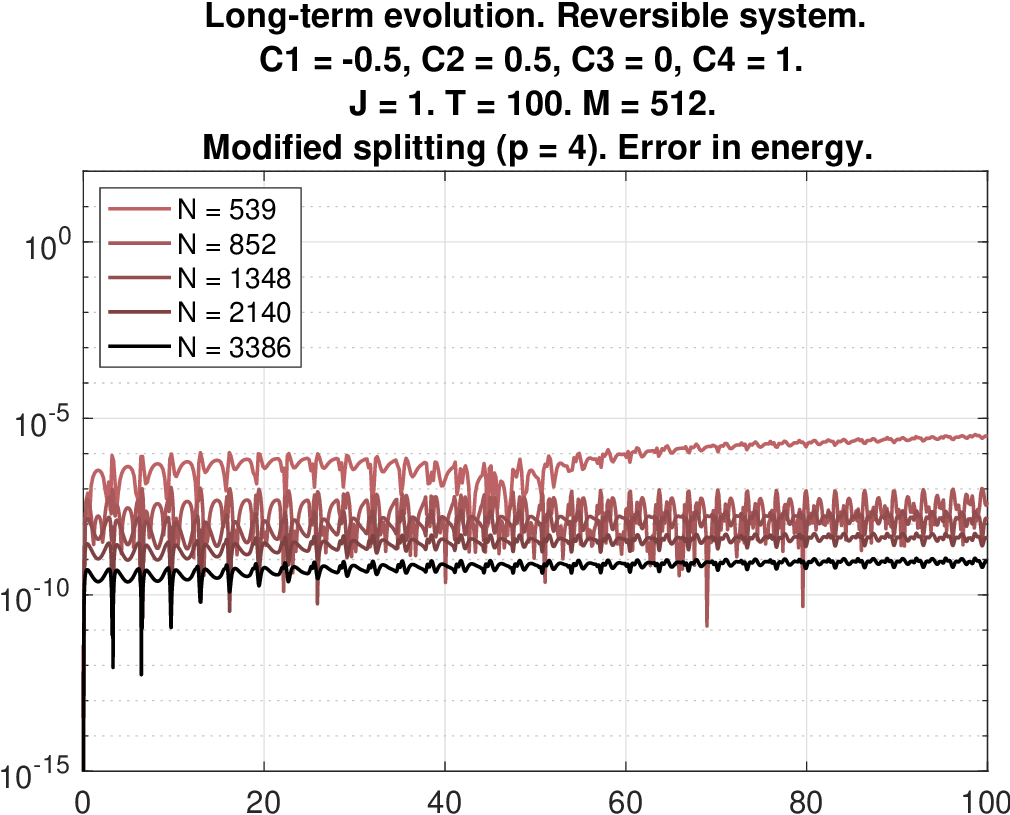} 
\end{center}
\caption{Corresponding results concerning the energy preservation in a long-term integration of a one-dimensional Gross--Pitaevskii equation, see Figure~\ref{fig:FigureTE1}.}
\label{fig:FigureTE3}
\end{figure}

\begin{figure}
\begin{center}
\includegraphics[width=6.5cm]{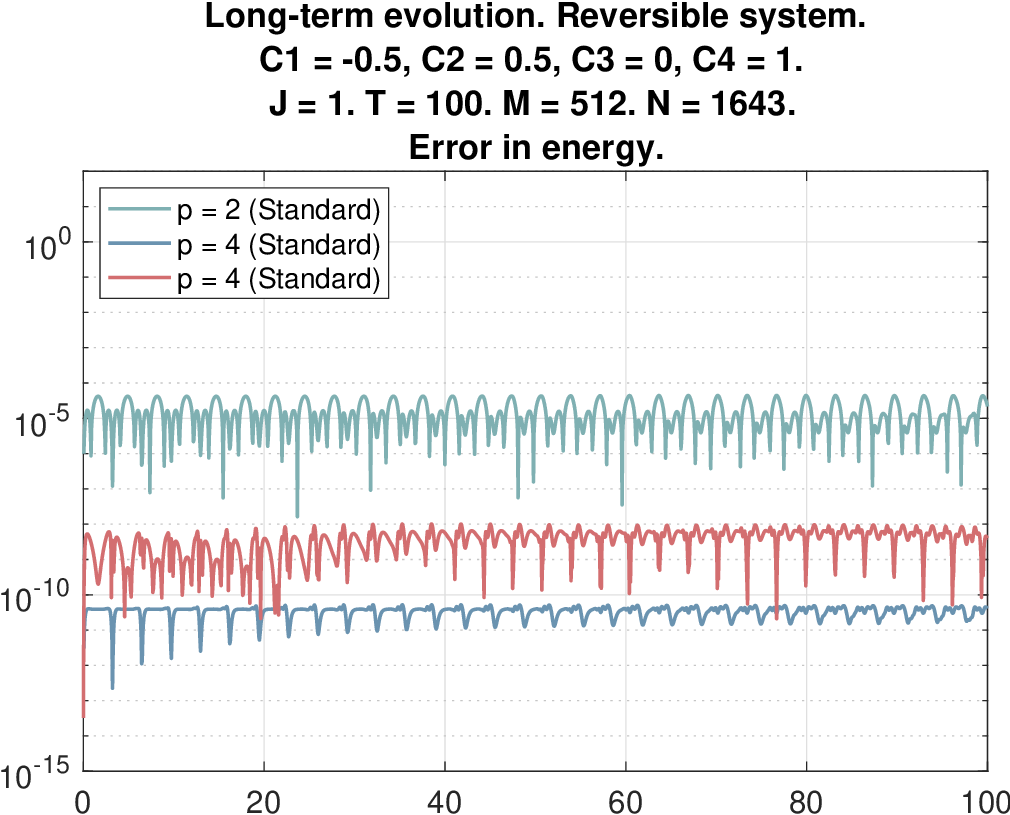} \quad
\includegraphics[width=6.5cm]{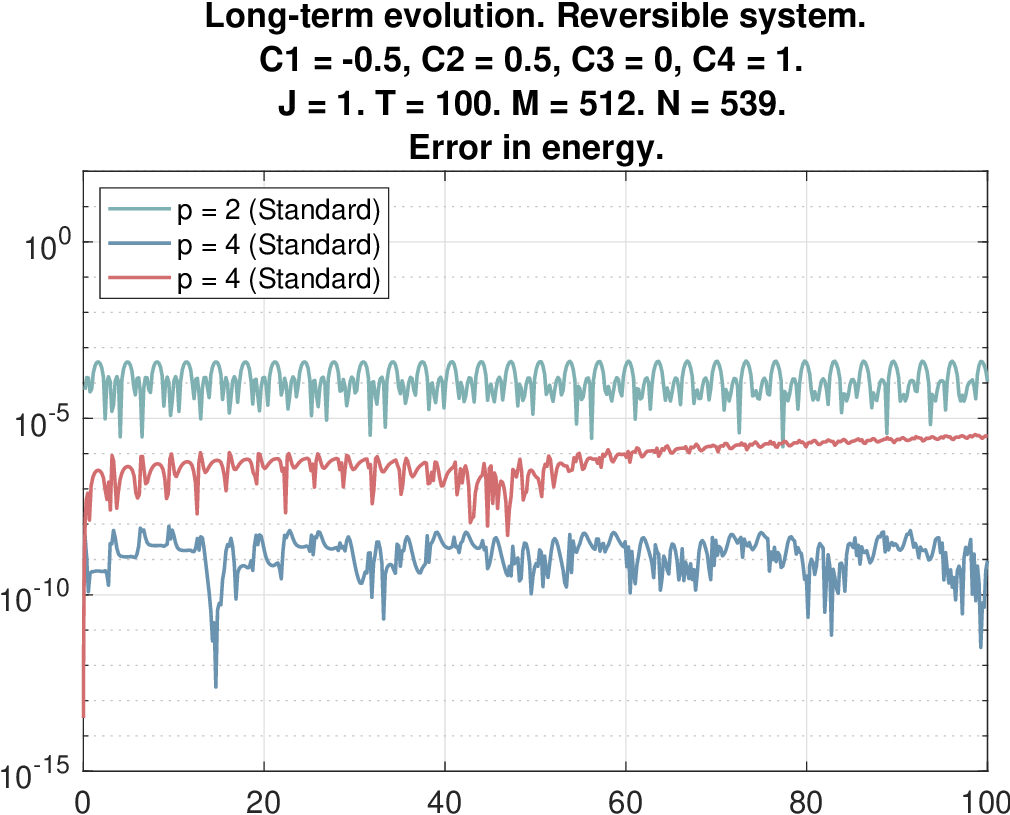} \\[2mm]
\includegraphics[width=6.5cm]{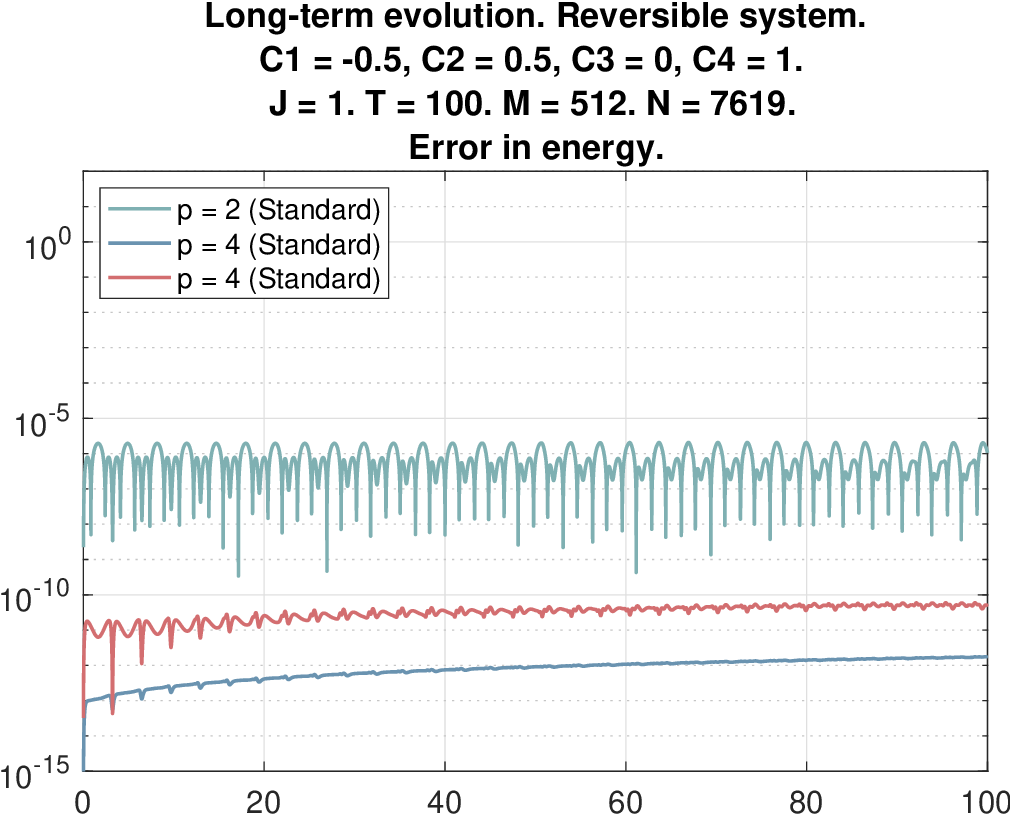} \quad
\includegraphics[width=6.5cm]{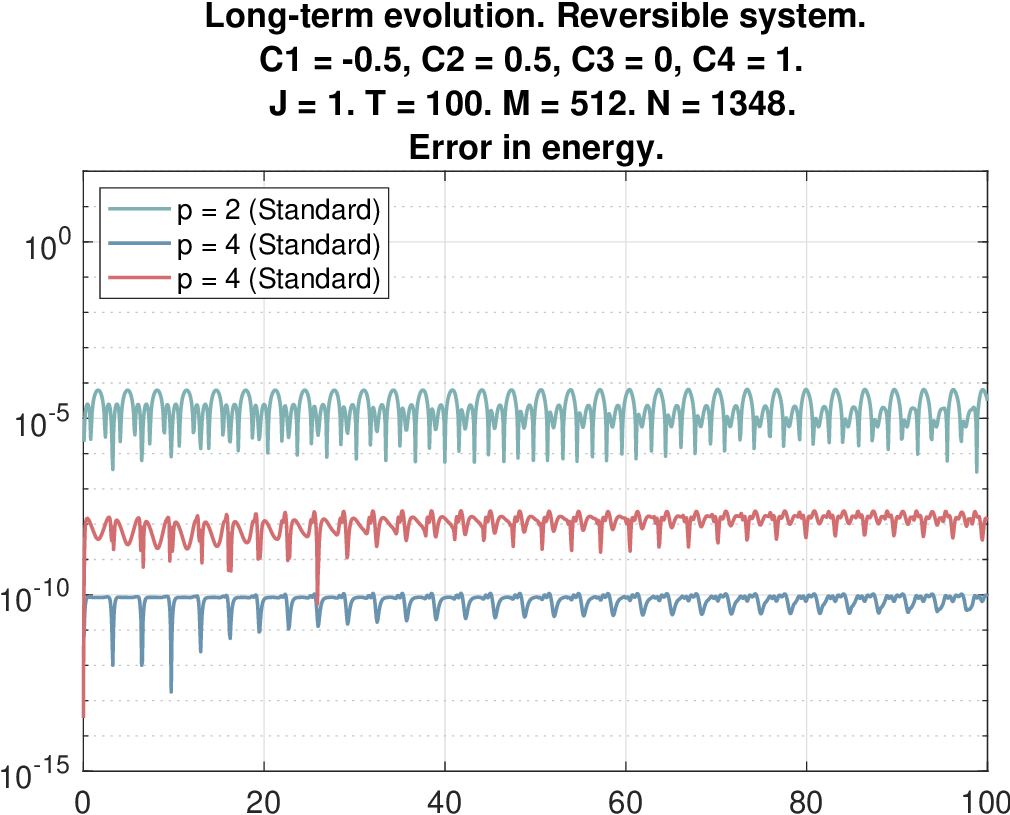} \\[2mm]
\includegraphics[width=6.5cm]{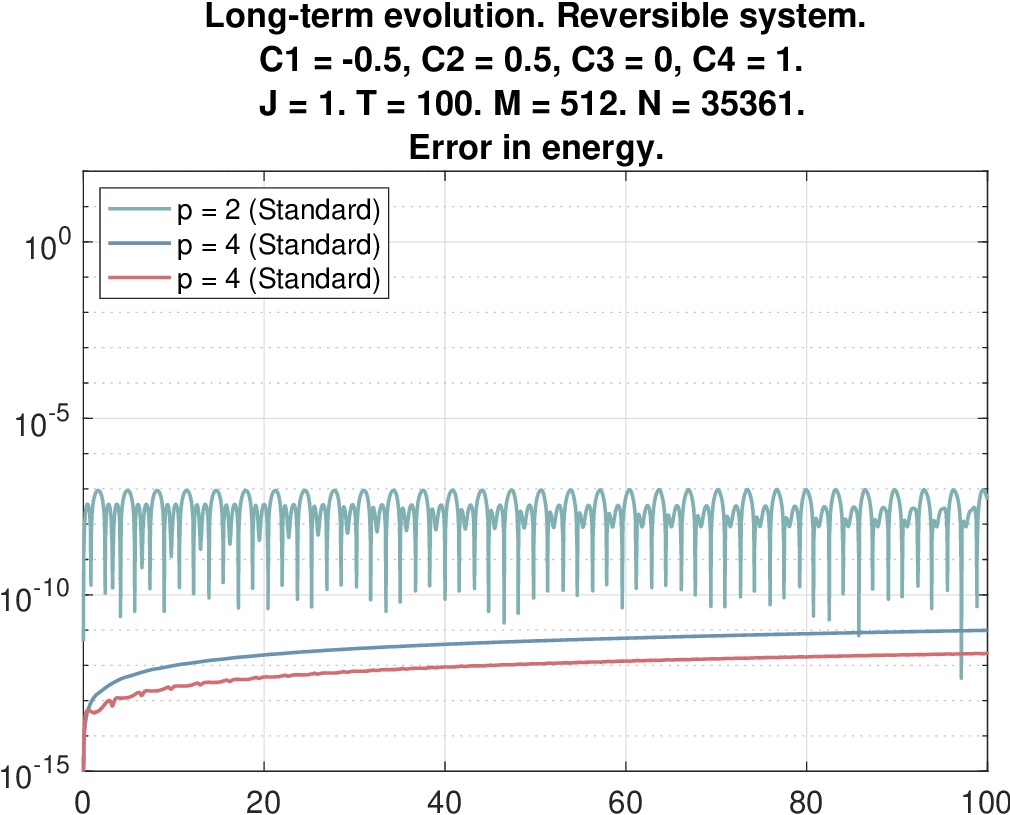} \quad
\includegraphics[width=6.5cm]{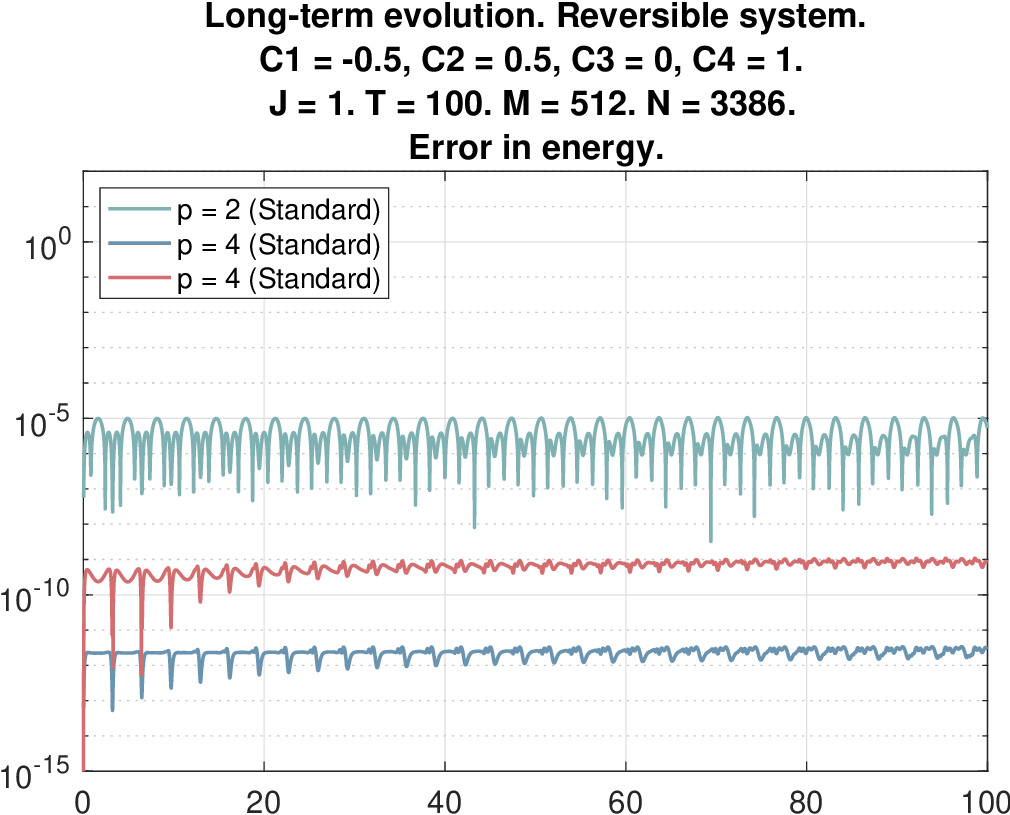} 
\end{center}
\caption{Corresponding results concerning the energy preservation in a long-term evolution of a one-dimensional Gross--Pitaevskii equation, see Figure~\ref{fig:FigureTE2}.
For low tolerances, the conservative strategy with local error estimate $\text{Err}_{\text{Local}}$ based on $\Psi_{\text{Modified}} - \Psi_{\text{Strang}}$ (left column) leads to highly accurate results.
For the adapted strategy with local error estimate $\tau^2 \, \text{Err}_{\text{Local}}$ and a reduced number of time steps (right column), which enhances efficiency, a connection between the errors in energy and the prescribed tolerances is observed.}
\label{fig:FigureTE4}
\end{figure}

\begin{figure}
\begin{center}
\includegraphics[width=6.5cm]{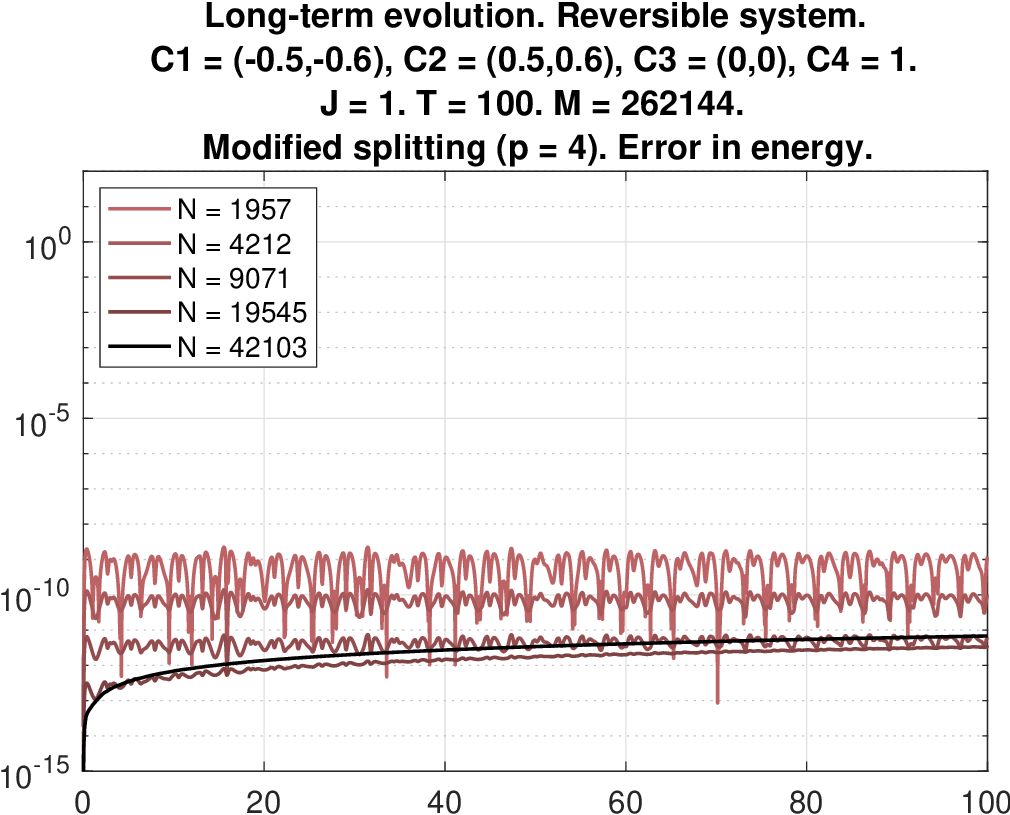} \quad
\includegraphics[width=6.5cm]{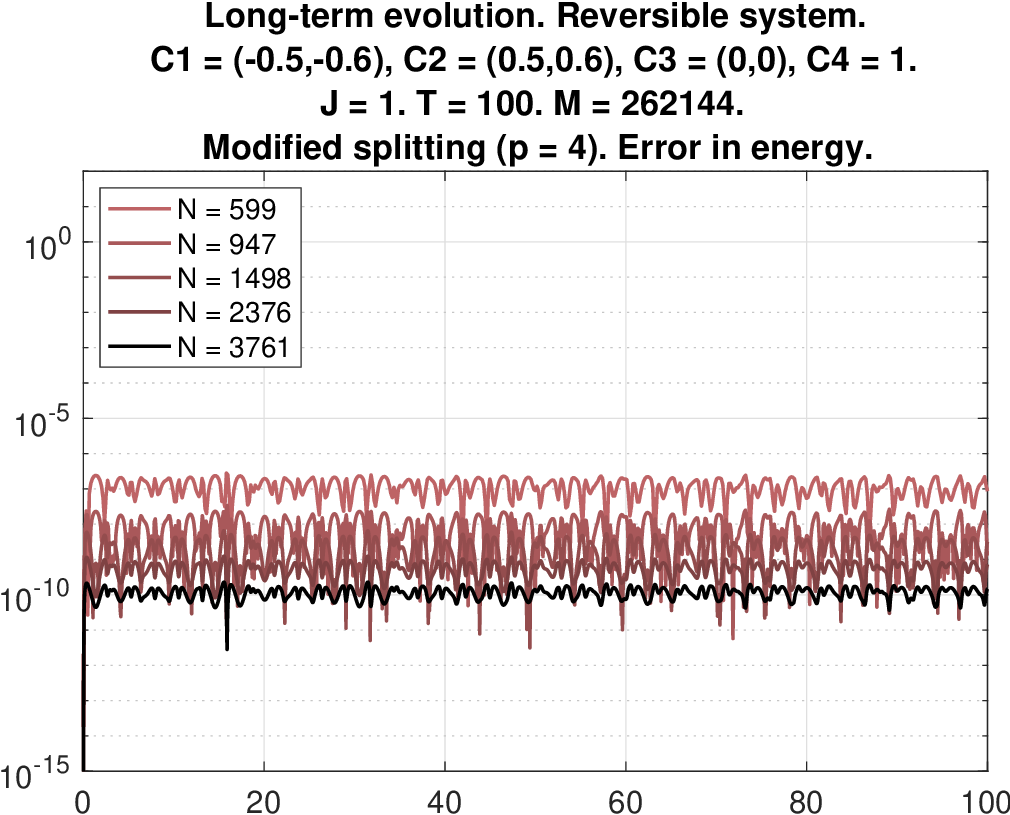} \\[2mm]
\includegraphics[width=6.5cm]{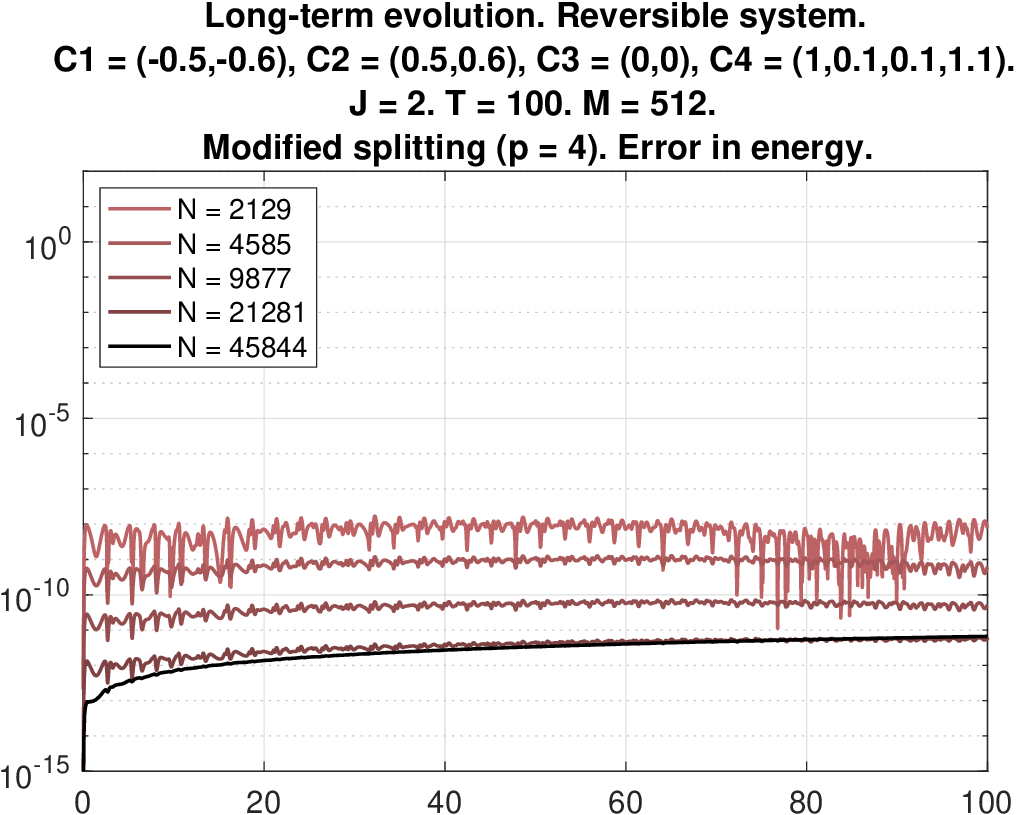} \quad
\includegraphics[width=6.5cm]{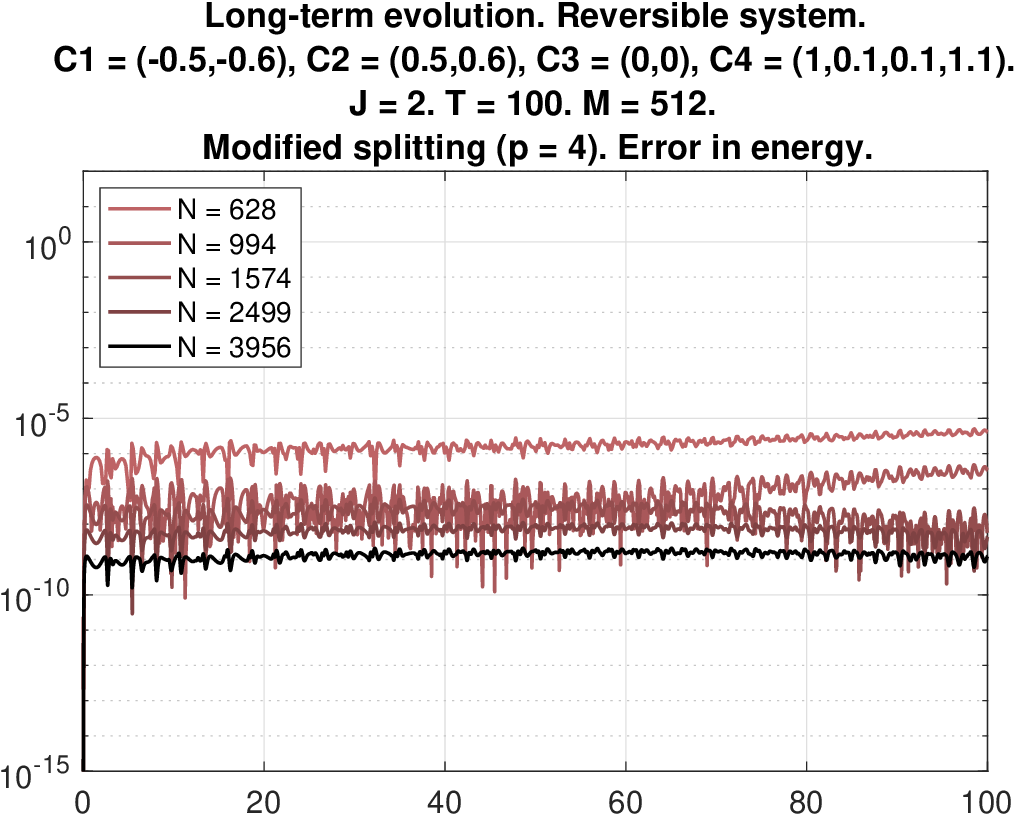} \quad
\end{center}
\caption{Corresponding results concerning the energy preservation in a long-term integration of a two-dimensional Gross--Pitaevskii equation (first row) and a one-dimensional two-component Gross--Pitaevskii system (second row), see Figure~\ref{fig:FigureTE3}.}
\label{fig:FigureTE5}
\end{figure}

\MyParagraph{Energy conservation}
In a second numerical experiment, we illustrate the favourable behaviour of adaptive modified operator splitting methods over longer times.
At first, we perform the time evolution of Gross--Pitaevskii systems by the standard splitting methods of orders two and four, respectively, for a prescribed sequence of time stepsizes. 
Then, we apply the fourth-order modified operator splitting method with a local error control based on the second-order Strang splitting method.
We point out that the use of a local error control to automatically adapt the time stepsizes releases the users from the task to find suitable time stepsize ranges.
That is, without an error control, suitable time increments had to be determined by several reruns, which involve high computational efforts, especially for the long-term evolution of Gross--Pitaevskii systems in three space dimensions. 
For the purpose of an initial test, we prescribe different tolerances $\text{Tol}$ ranging from~$10^{-4}$ to~$10^{-8}$ and maintain a tentative approach in the sense that we leave aside the option to enlarge the time increments by additional scaling factors. 
Subsequently, we replace the local error estimate~$\text{Err}_{\text{Local}}$ based on the difference of current approximations~$\Psi_{\text{Modified}}$ and~$\Psi_{\text{Strang}}$ by the local error estimate $\tau^2 \, \text{Err}_{\text{Local}}$ and repeat the long-term evolution. 
For the resulting numbers of time steps~$N$, we then once again apply the Strang splitting method and the optimised standard fourth-order splitting method with the corresonding equidistant time increments. 
We report that all splitting methods show an excellent mass conservation, which is a well-known intrinsic distinction to other classes of time integration methods such as Runge--Kutta and linear multistep methods. 
Concerning the conservation of total energy, also indicators for the achievable accuracy in the solution, the obtained results are displayed in Figures~\ref{fig:FigureTE1} to~\ref{fig:FigureTE4}.
We observe an excellent performance of higher-order standard and adaptive modified operator splitting methods compared to lower-order schemes. 
Specifically for more involved nonlinear settings, the results also suggest that adaptive modified splitting methods with the conservative strategy based on $\Psi_{\text{Modified}} - \Psi_{\text{Strang}}$ are preferable, whenever highly accurate results enclosing a favourable energy preservation are desirable, whereas efficiency is enhanced with the adapted stategy based on $\tau^2 \, \text{Err}_{\text{Local}}$. 
In order to confirm the reliability of our adaptive approach for sequences of prescribed tolerances, we include the outcomes for a Gross--Pitaevskii equation in two space dimensions and a one-dimensional two-component Gross--Pitaevskii system in Figure~\ref{fig:FigureTE5}.

\MyParagraph{Imaginary time propagation and time evolution}
In a final numerical experiment, we consider Gross--Pitaevskii equations in one, two, and three space dimensions for different parameter ranges.
We contrast the case of a moderate constant in the nonlinearity to the case of an additional multiple lattice potential and a large constant in the nonlinearity, which is of particular relevance for physical experimental set-ups. 
Simplified linear cases, where the exact ground state and time-dependent solutions are known, are once again included for the purpose of validation. 
Moreover, it is remarkable that the knowledge of the ground state solutions and the associated chemical potentials define the time-dependent solutions $\Psi(x, t) = \ee^{- \, \ii \, \mu \, t} \, \Phi(x)$ for $(x, t) \in \Omega \times [t_0, T]$, see~\eqref{eq:GPSGroundstate}.  
Consequently, in modulus, the time-dependent solutions coincide with the ground states.
Based on the adaptive fourth-order modified operator splitting method, we first perform the imaginary time propagation for the numerical computation of the ground state solutions and then evolve the Gross--Pitaevskii systems in time.
For the linear cases, we prescribe constant initial states, and for the nonlinear cases, we determine the Thomas--Fermi approximations~\eqref{eq:ThomasFermi}.
We point out that the second-order Strang splitting method and the optimised standard fourth-oder splitting method fail in the imaginary time propagation, when applied with time increment $\tau = \frac{1}{10}$, already in a single space dimension.
On the contrary, the adaptive fourth-order modified operator splitting method yields reliable results in the ground state computation for the initial time stepsize $\tau = \frac{1}{10}$ and the tolerance $\text{Tol} = 10^{-5}$.  
The profiles of the chosen real-valued initial states for the imaginary time propagation, the resulting ground state solutions at imaginary time~$T$, and the modulus of the time-dependent solutions evolved from $t_0 = 0$ to the final time~$T$ are illustrated in Figures~\ref{fig:FigureITTE1} to~\ref{fig:FigureITTE3}.
In space dimension three, we determine the section along $x_3 = 0$ and display the corresponding values $\Psi(x_1,x_2,0, t)$.
As additional information and indicators for the overall accuracy, we include the total energy and the corresponding errors in the headlines of the figures.
Alltogether, we conclude that the adaptive fourth-order modified operator splitting method leads to reliable and accurate results in the imaginary time propagation and the time evolution of Gross--Pitaevskii equations. 
\section{Conclusions}
\label{sec:Conclusions}
\MyParagraph{Summary}
The present work has been devoted to the introduction of adaptive modified operator splitting methods for the reliable imaginary time propagation and the efficient time evolution of two-component Gross--Pitaevskii systems.
In a series of numerical tests, we have demonstrated the excellent performance of a specific adaptive fourth-order modified operator splitting method involving positive coefficients in comparison with renowned standard splitting methods. 
We have implemented two strategies for an inexpensive expedient local error control that are advantageous for successful computations. 

\MyParagraph{Generalisations}
Straightforward extensions that require minor adaptations of our approach concern more general particle interaction terms such as quintic or fractional power nonlinearities.
Besides, it appears feasible to treat Gross--Pitaevskii equations with fractional Laplacians~\MyCite{AntoineEtAl2018} by means of a straightforward calculation of the iterated commutators, neglecting possible simplifications in the resulting relations. 
The numerical simulation of rotating Bose--Einstein condensates under the influence of external magnetic fields, modelled by Gross--Pitaevskii equations with additional rotation terms, is based on suitable reformulations within the rotating frames.
This yields non-autonomous problems of similar structures, which are resolved by the combination of Magnus-type integrators with modified operator splitting methods~\MyCite{BaderBlanesCasasThalhammer2019}.

\MyParagraph{Future investigations}
Our studies in the near future will be dedicated to the efficient numerical simulation of concrete settings that arise in actual laboratory set-ups.  
Amongst others, this will include the investigation of strongly interacting quantum systems driven by laser kicks and double species Bose--Einstein condensates.



\end{document}